\documentclass[12pt,leqno]{article}
\usepackage{amssymb}
\usepackage[mathscr]{eucal}
\usepackage{amsmath,amssymb,latexsym,theorem,bbm}
\usepackage{color}
\usepackage{appendix}

\newcommand{\NN}{\mathbb{N}}

\newcommand{\RR}{\mathbb{R}}

\newcommand{\ZZ}{\mathbb{Z}}

\newcommand{\bA}{{\boldsymbol{A}}}

\newcommand{\tA}{\widetilde{A}}

\newcommand{\tF}{\widetilde{F}}
\newcommand{\tcF}{\widetilde{\cF}}
\newcommand{\bF}{{\boldsymbol{F}}}

\newcommand{\tk}{\widetilde{k}}

\newcommand{\tp}{\widetilde{p}}

\newcommand{\tPP}{\widetilde{\PP}}

\newcommand{\bx}{{\boldsymbol{x}}}
\newcommand{\bX}{{\boldsymbol{X}}}
\newcommand{\by}{{\boldsymbol{y}}}
\newcommand{\bY}{{\boldsymbol{Y}}}
\newcommand{\tbX}{\widetilde{\bX}}

\newcommand{\tbW}{\widetilde{\bW}}

\newcommand{\bW}{{\boldsymbol{W}}}

\newcommand{\bxi}{{\boldsymbol{\xi}}}

\newcommand{\tOmega}{\widetilde{\Omega}}

\newcommand{\bzero}{{\boldsymbol{0}}}

\newcommand{\cB}{{\mathcal B}}
\newcommand{\cC}{{\mathcal C}}
\newcommand{\cD}{{\mathcal D}}

\newcommand{\cG}{{\mathcal G}}
\newcommand{\cF}{{\mathcal F}}
\newcommand{\cH}{{\mathcal H}}

\newcommand{\cM}{{\mathcal M}}

\newcommand{\cN}{{\mathcal N}}
\newcommand{\cP}{{\mathcal P}}

\newcommand{\cS}{{\mathcal S}}

\newcommand{\hcB}{\widehat{\cB}}

\newcommand{\cc}{\mathrm{c}}
\newcommand{\dd}{\mathrm{d}}
\newcommand{\ee}{\mathrm{e}}
\newcommand{\ii}{\mathrm{i}}

\newcommand{\EE}{\operatorname{\mathbb{E}}}
\newcommand{\PP}{\operatorname{\mathbb{P}}}

\newcommand{\tN}{\widetilde{N}}

\newcommand{\vare}{\varepsilon}

\renewcommand{\mid}{\,|\,}

\renewcommand{\leq}{\leqslant}
\renewcommand{\geq}{\geqslant}

\newcommand{\stochi}{\stackrel{\PP^{(i)}}{\longrightarrow}}

\newcommand{\as}{\stackrel{{\mathrm{a.s.}}}{\longrightarrow}}

\newcommand{\bbone}{\mathbbm{1}}

\newcommand{\nt}{{\lfloor nt\rfloor}}

\newcommand{\proofend}{\hfill\mbox{$\Box$}}

\numberwithin{equation}{section}

\theoremstyle{change} \theorembodyfont{\em}
\newtheorem{Lem}{Lemma.}[section]
\newtheorem{Thm}[Lem]{Theorem.}

\newtheorem{Cor}[Lem]{Corollary.}
\newtheorem{Def}[Lem]{Definition.}

\theorembodyfont{\rm}
\newtheorem{Rem}[Lem]{Remark.}

\begin{document}

\begin{center}
 {\bfseries\Large Yamada--Watanabe results for stochastic differential\\[2mm]
                   equations with jumps}
 \\[5mm]

 {\sc\large
  M\'aty\'as $\text{Barczy}^{*,\diamond}$,
  \ Zenghu $\text{Li}^{**}$,
  \ Gyula $\text{Pap}^{***}$}

\end{center}

\vskip0.2cm

\noindent
 * Faculty of Informatics, University of Debrecen,
   Pf.~12, H--4010 Debrecen, Hungary.

\noindent
 ** School of Mathematical Sciences, Beijing Normal University,
     Beijing 100875,  People's Republic of China.

\noindent
 *** Bolyai Institute, University of Szeged,
     Aradi v\'ertan\'uk tere 1, H--6720 Szeged, Hungary.

\noindent e--mails: barczy.matyas@inf.unideb.hu (M. Barczy),
                    lizh@bnu.edu.cn (Z. Li),
                    papgy@math.u-szeged.hu (G. Pap).

\noindent $\diamond$ Corresponding author.

%\vskip0.2cm

%\centerline{\sl December 8, 2013.}

\renewcommand{\thefootnote}{}
\vspace*{0.2cm}
\footnote{
The research of M. Barczy and G. Pap was realized in the frames of
 T\'AMOP 4.2.4.\ A/2-11-1-2012-0001 ,,National Excellence Program --
 Elaborating and operating an inland student and researcher personal support
 system''.
The project was subsidized by the European Union and co-financed by the
 European Social Fund.
Z. Li has been partially supported by NSFC under Grant No.\ 11131003
 and 973 Program under Grant No.\ 2011CB808001.}

\vspace*{-10mm}

\begin{abstract}
Recently, Kurtz (2007, 2014) obtained a general version of the
 Yamada--Watanabe and Engelbert theorems relating existence and uniqueness of
 weak and strong solutions of stochastic equations covering also the case of
 stochastic differential equations with jumps.
Following the original method of Yamada and Watanabe (1971), we give
 alternative proofs for the following two statements: pathwise uniqueness
 implies uniqueness in the sense of probability law, and weak existence
 together with pathwise uniqueness imply strong existence for stochastic
 differential equations with jumps.
\end{abstract}

\section{Introduction}
\label{section_intro}

In order to prove existence and pathwise uniqueness of a strong solution for
 stochastic differential equations, it is an important issue to clarify the
 connections between weak and strong solutions.
The first pioneering results are due to Yamada and Watanabe \cite{YamWat}
 for certain stochastic differential equations driven by Wiener processes.

We investigate stochastic differential equations with jumps.
Let \ $U$ \ be a second-countable locally compact Hausdorff space equipped
 with its Borel $\sigma$-algebra \ $\cB(U)$.
\ Let \ $m$ \ be a $\sigma$-finite Radon measure on \ $(U, \cB(U))$, \ meaning
 that the measure of compact sets is always finite.
Let \ $U_0, U_1 \in \cB(U)$ \ be disjoint subsets.
Let \ $d, r \in \NN$.
\ Let \ $b : [0, \infty) \times \RR^d \to \RR^d$,
 \ $\sigma : [0, \infty) \times \RR^d \to \RR^{d\times r}$,
 \ $f : [0, \infty) \times \RR^d \times U \to \RR^d$ \ and
 \ $g : [0, \infty) \times \RR^d \times U \to \RR^d$ \ be Borel measurable
 functions, where \ $[0, \infty) \times \RR^d \times U$ \ is equipped with its
 Borel $\sigma$-algebra
 \ $\cB([0, \infty) \times \RR^d \times U)
    = \cB([0, \infty)) \otimes \cB(\RR^d) \otimes \cB(U)$
 \ (see, e.g., Dudley \cite[Proposition 4.1.7]{Dud}).
Consider a stochastic differential equation (SDE)
 \begin{equation}\label{SDE_X_YW}
  \begin{aligned}
   \bX_t
   &= \bX_0 + \int_0^t \sigma(s, \bX_s) \, \dd \bW_s
            + \int_0^t \int_{U_0} f(s, \bX_{s-}, u) \, \tN(\dd s, \dd u) \\
   &\quad
      + \int_0^t b(s, \bX_s) \, \dd s
      + \int_0^t \int_{U_1} g(s, \bX_{s-}, u) \, N(\dd s, \dd u) , \qquad
   t \in [0, \infty)  ,
  \end{aligned}
 \end{equation}
 where \ $(\bW_t)_{t\geq0}$ \ is an $r$-dimensional standard Brownian motion,
 \ $N(\dd s, \dd u)$ \ is a Poisson random measure on
 \ $(0, \infty)  \times U$ \ with intensity measure \ $\dd s \, m(\dd u)$,
 \ $\tN(\dd s, \dd u) := N(\dd s, \dd u) - \dd s \, m(\dd u)$, \ and
 \ $(\bX_t)_{t\geq0}$ \ is a suitable process with values in \ $\RR^d$.

Yamada and Watanabe \cite{YamWat} proved that weak existence and pathwise
 uniqueness imply uniqueness in the sense of probability law and strong
 existence for the SDE  \eqref{SDE_X_YW} with \ $f = 0$ \ and \ $g = 0$.
\ Engelbert \cite{Eng} and Cherny \cite{Che2} extended this result to a
 somewhat more general class of equations and gave a converse in which the
 roles of existence and uniqueness are reversed, that is, joint uniqueness in
 the sense of probability law (see, Engelbert \cite[Definition 5]{Eng}) and
 strong existence imply pathwise uniqueness.
The original Yamada--Watanabe result arises naturally in the procedure of
 proving existence of solutions of a SDE; for a detailed discussion, see
 Kurtz \cite[pages 1--2]{Kur2}.

Jacod \cite{Jac} generalized the above mentioned result of Yamada and Watanabe
 for a SDE driven by a semimartingale, where the coefficient may depend on the
 paths both of the solution and of the driving process.
The Yamada--Watanabe result has been generalized by Ondrej\'at \cite{Ond} and
 R\"ockner et al.\ \cite{RocSchZha} for stochastic evolution equations in
 infinite dimensions, and by Tappe \cite{Tap} for semilinear stochastic
 partial differential equations with path-dependent coefficients.

Recently, there has been a renewed interest in generalizations of the results
 of  Yamada and Watanabe \cite{YamWat}.
Kurtz \cite{Kur1}, \cite{Kur2} continued the direction of Engelbert \cite{Eng}
 and Jacod \cite{Jac}.
He studied general stochastic models which relate stochastic inputs with
 stochastic outputs, and obtained a general version of the Yamada--Watanabe
 and Engelbert theorems relating existence and uniqueness of weak and strong
 solutions of stochastic models with the message that the original results
 are not limited to SDEs driven by Wiener processes.
In order to derive the original Yamada--Watanabe results from this general
 theory, proofs of pathwise uniqueness require appropriate adaptedness
 conditions, so two new notions, compatibility and partial compatibility
 between inputs and outputs have been introduced.
Due to Example 3.9 in Kurtz \cite{Kur1} and Page 7 in Kurtz \cite{Kur2}, the
 results are valid for SDEs driven by a Wiener process and Poisson random
 measures.

Following the ideas of Yamada and Watanabe \cite{YamWat}, we are going to give
 alternative proofs for the following two statements:

\begin{Thm}\label{Thm_pathwise_1}
Pathwise uniqueness for the SDE \eqref{SDE_X_YW} implies uniqueness in the
 sense of probability law.
\end{Thm}

\begin{Thm}\label{Thm_pathwise_2}
Weak existence and pathwise uniqueness for the SDE \eqref{SDE_X_YW} imply
 strong existence.
\end{Thm}

Note that Theorems \ref{Thm_pathwise_1} and \ref{Thm_pathwise_2} are
 generalizations of Proposition 1 and Corollary 1 in Yamada and Watanabe
 \cite{YamWat} (we do not intend to deal with generalization of their
 Corollary 3).
The definition of weak and strong solutions of the SDE \eqref{SDE_X_YW},
 pathwise uniqueness for the SDE \eqref{SDE_X_YW} and uniqueness in the sense
 of probability law, and a detailed, precise formulation of Theorem
 \ref{Thm_pathwise_2} will be given in the paper.
In the course of the proofs we developed a sequence of lemmas discussing
 several kinds of measurability, see Lemmas \ref{measurability} and \ref{k},
 and we also presented a key observation on the preservation of the joint
 distribution of the parts of the SDE \eqref{SDE_X_YW}, see Lemmas
 \ref{Lemma_distribution} and \ref{Lemma_distribution2}.

Our alternative proofs show the power of the original method of Yamada and
 Watanabe \cite{YamWat}, these proofs can be followed step by step and every
 technical detail is transparent in the paper.
This raises a question whether Kurtz's result could be proved via the
 walked-out path by Yamada and Watanabe.

Note that Situ \cite[Theorem 137]{Sit} also considered the SDE
 \eqref{SDE_X_YW} with \ $\RR^d \setminus \{\bzero\}$ \ instead of \ $U$ \ and
 with \ $g = 0$, \ and proved Theorems \ref{Thm_pathwise_1} and
 \ref{Thm_pathwise_2} under the resctrictive assumption
 \begin{equation}\label{Situ}
   \int_{\RR^d \setminus \{\bzero\}} \frac{\|u\|^2}{1 + \|u\|^2} \, m(\dd u)
   < \infty .
 \end{equation}
This assumption was needed for introducing an auxiliary c\`adl\`ag process in
 Lemma 139 in Situ \cite{Sit}.
In fact, one can get rid of condition \eqref{Situ} by using the space of point
 measures on \ $\RR_+ \times U$ \ as the space of trajectories of Poisson
 point processes instead of the space of c\`adl\`ag functions, see the proofs
 of Theorems \ref{Thm_pathwise_1} and \ref{Thm_pathwise_2}.
We call the attention that in the literature the result of Situ
 \cite[Theorem 137]{Sit} has been usually referred to without checking
 condition \eqref{Situ}, see, e.g., Li and Mytnik
 \cite[equation (3.1)]{LiMyt}, Dawson and Li
 \cite[equation (2.9)]{DawLi}, D\"oring and Barczy
 \cite[equation (3.23)]{DorBar} and Li and Pu
 \cite[equations (4.6) and (5.1)]{LiPu}, but Theorem \ref{Thm_pathwise_2}
 covers these situations as well.

We remark that Zhao \cite{Zhao} already adapted the original method of
 Yamada and Watanabe for the SDE \eqref{SDE_X_YW} driven only by a compensated
 Poisson random measure, i.e., with \ $\sigma = 0$ \ and \ $g = 0$, \ but for
 processes with values in a separable Hilbert space instead of \ $\RR^d$-valued
 processes.
Comparing with the results of the present paper, note that we explicitly stated
 and proved in Theorem \ref{Thm_pathwise_1} that pathwise uniqueness for the SDE
 \eqref{SDE_X_YW} implies uniqueness in the sense of probability law.

\section{Preliminaries}
\label{section_preliminaries}

Let \ $\ZZ_+$, \ $\NN$, \ $\RR$, \ $\RR_+$  \ and \ $\RR_{++}$ \ denote the set
 of non-negative integers, positive integers, real numbers, non-negative real
 numbers and positive real numbers, respectively.
For \ $x , y \in \RR$, \ we will use the notations
 \ $x \land y := \min \{x, y\} $.
\ By \ $\|\bx\|$ \ and \ $\|\bA\|$, \ we denote the Euclidean norm of a vector
 \ $\bx \in \RR^d$ \ and the induced matrix norm of a matrix
 \ $\bA \in \RR^{d\times d}$, \ respectively.
Throughout this paper, we make the conventions
 \ $\int_a^b := \int_{(a,b]}$ \ and \ $\int_a^\infty := \int_{(a,\infty)}$ \ for any
 \ $a, b \in \RR$ \ with \ $a < b$.
By \ $C(\RR_+, \RR^\ell)$ \ and \ $D(\RR_+, \RR^\ell)$ \ we denote the set of
 continuous and c\`adl\`ag \ $\RR^\ell$-valued functions defined on \ $\RR_+$,
 \ equipped with a metric inducing the local uniform topology (see, e.g.,
 Jacod and Shiryaev \cite[Section VI.1a]{JacShi}) and a metric inducing the
 so-called Skorokhod topology (see, e.g., Jacod and Shiryaev
 \cite[Theorem VI.1.14]{JacShi}), respectively.
Moreover, \ $\cC(\RR_+, \RR^\ell)$ \ and \ $\cD(\RR_+, \RR^\ell)$ \ denote the
 corresponding Borel \ $\sigma$-algebras on them.

Recall that \ $U$ \ is a second-countable locally compact Hausdorff space.
Note that \ $U$ \ is homeomorphic to a separable complete metric space, see,
 e.g., Kechris \cite[Theorem 5.3]{Kec}.
For our later purposes, we recall the notion of the space of point
 measures on \ $\RR_+ \times U$, \ of the space of simple point measures on
 \ $\RR_+ \times U$, \ and of the vague convergence.
We follow Resnick \cite[Chapter 3]{Res} and Ikeda and Watanabe
 \cite[Chapter I, Sections 8 and 9]{IkeWat}.

A point measure on \ $\RR_+ \times U$ \ is a measure \ $\pi$ \ of the
 following form: let \ $F \subset \NN$ \ and let
 \ $\{ (t_i, u_i) : i \in F \}$ \ be a countable collection of (not
 necessarily distinct) points of \ $\RR_+ \times U$, \ and let
 \[
   \pi := \sum_{i\in F} \delta_{(t_i,u_i)}
 \]
 assuming also that \ $\pi([0,t] \times B) <\infty$ \ for all \ $t \in \RR_+$
 \ and compact subsets \ $B \in \cB(U)$ \ (i.e., \ $\pi$ \ is a Radon measure
 meaning that the measure of compact sets is always finite, and consequently,
 it is locally finite), where \ $\delta_{(t_i,u_i)}$ \ denotes the Dirac measure
 concentrated on the point \ $(t_i, u_i)$.
\ Thus
 \[
   \pi([0,t] \times B) = \#\{i \in F : (t_i,u_i) \in [0,t]\times B\} , \qquad
   t \in \RR_+ , \quad B \in \cB(U) .
 \]

A point function (or point pattern) \ $p$ \ on \ $U$ \ is a mapping
 \ $p: D(p) \to U$, \ where the domain \ $D(p)$ \ is a countable subset of
 \ $\RR_{++}$ \ such that
 \ $\{s \in D(p) : \text{$s \in (0, t]$, \ $p(s) \in B$}\}$ \ is finite for
 all \ $t \in \RR_+$ \ and compact subsets \ $B \in \cB(U)$.
\ The counting measure \ $N_p$ \ on \ $\RR_{++} \times U$ \ corresponding to
 \ $p$ \ is defined by
 \[
   N_p((0,t] \times B)
   := \#\{s \in D(p) : \text{$s \in (0, t]$, \ $p(s) \in B$}\} , \qquad
   t \in \RR_{++} , \quad B \in \cB(U) .
 \]
Note that there is a (natural) bijection between the set of point functions
 on \ $U$ \ and the set of point measures \ $\pi$ \ on \ $\RR_+\times U$
 \ with \ $\pi(\{t\} \times U) \leq 1$, \ $t \in \RR_{++}$, \ and
 \ $\pi(\{0\} \times U) = 0$.
\ Namely, if \ $p: D(p) \to U$ \ is a point function, then the corresponding
 point measure is its counting measure
 \ $N_p = \sum_{t\in D(p)} \delta_{(t,p(t))}$.
\ The set of all point measures on \ $\RR_+ \times U$ \ will be denoted by
 \ $M(\RR_+ \times U)$, \ and define a \ $\sigma$-algebra
 \ $\cM(\RR_+ \times U)$ \ on it to be the smallest $\sigma$-algebra
 containing all sets of the form
 \[
   \{ \pi \in M(\RR_+\times U) : \pi([0, t] \times B) \in A \} \qquad
   \text{for \ $t \in \RR_+$, \ $B \in \cB(U)$, \ $A \in \cB([0, \infty])$.}
 \]
Alternatively, \ $\cM(\RR_+ \times U)$ \ is the smallest $\sigma$-algebra
 making all the mappings
 \ $M(\RR_+ \times U) \ni \pi \mapsto \pi([0, t] \times B) \in [0, \infty]$,
 \ $t \in \RR_+$, \ $B \in \cB(U)$, \ measurable.

Note that there is a (natural) bijection between the set of point processes
 (randomized point functions) \ $p$ \ defined on a probability space
 \ $(\Omega, \cF, \PP)$ \ with values in the space of point functions on \ $U$
 \ (in the sense of Ikeda and Watanabe
 \cite[Chapter I, Definition 9.1]{IkeWat}) and the set of
 \ $\cF / \cM(\RR_+ \times U)$-measurable mappings
 \ $p : \Omega \to M(\RR_+ \times U)$ \ with
 \ $p(\omega)(\{t\} \times U) \leq 1$
 \ for all \ $\omega \in \Omega$ \ and \ $t \in \RR_{++}$, \ and
 \ $p(\omega)(\{0\}\times U) = 0$ \ for all \ $\omega \in \Omega$ \ (which
 are (special) point processes in the sense of Resnick \cite[page 124]{Res}).

A point process \ $p$ \ on \ $U$ \ is called a Poisson point process if its
 counting measure \ $N_p$ \ is a Poisson random measure on \ $\RR_+ \times U$
 \ (for the definition of Poisson random measure see, e.g., Ikeda and Watanabe
 \cite[Chapter I, Definition 8.1]{IkeWat}).
A Poisson point process is stationary if and only if its intensity measure is
 of the form \ $\dd s \, \nu(\dd u)$ \ for some measure \ $\nu$ \ on
 \ $(U, \cB(U))$, \ which is called its charateristic measure.
\ If \ $\nu$ \ is a Radon measure, then \ $N_p((0, t] \times B)$ \ is Poisson
 distributed with parameter \ $t \nu(B) \in \RR_+$, \ hence
 \ $\{s \in D(p) : \text{$s \in (0, t]$, \ $p(s) \in B$}\}$ \ is finite with
 probability one for all \ $t \in \RR_+$ \ and compact subsets
 \ $B \in \cB(U)$.
\ Consequently, a stationary Poisson point process with a Radon charateristic
 measure is a stationary Poisson point process in the sense of Ikeda and
 Watanabe \cite[Chapter I, Definition 9.1]{IkeWat}.

Next we recall vague convergence.
Let \ $C_\cc(\RR_+ \times U, \RR_+)$ \ be the space of \ $\RR_+$-valued
 continuous functions defined on \ $\RR_+ \times U$ \ with compact support.
For \ $\pi, \pi_n \in M(\RR_+ \times U)$, \ $n\in\NN$, \ we say that \ $\pi_n$
 \ converges vaguely to \ $\pi$ \ as \ $n\to\infty$ \ if
 \[
   \lim_{n\to\infty} \int_{\RR_+\times U} f \, \dd \pi_n
   = \int_{\RR_+\times U} f \, \dd \pi
 \]
 for all \ $f \in C_\cc(\RR_+ \times U, \RR_+)$.
\ For a topology on \ $M(\RR_+ \times U)$ \ giving this notion of convergence,
 see page 140 in Resnick \cite{Res}.
Recall that \ $\cM(\RR_+ \times U)$ \ coincides with the Borel
 \ $\sigma$-algebra generated by the open sets with respect to the vague
 topology on \ $M(\RR_+ \times U)$, \ see, e.g., Resnick
 \cite[Exercises 3.4.2(b) and 3.4.5]{Res}.

In what follows we equip the spaces \ $C(\RR_+, \RR^\ell)$,
 \ $D(\RR_+, \RR^\ell)$, \ $\ell\in\NN$, \ and \ $M(\RR_+\times U)$ \ with some
 \ $\sigma$-algebras that will be used later on.
For each \ $\ell \in \NN$, \ let us equip \ $C(\RR_+, \RR^\ell)$ \ and
 \ $D(\RR_+, \RR^\ell)$ \ with the \ $\sigma$-algebras
 \[
   \cC_t(\RR_+, \RR^\ell) := \varphi_t^{-1}(\cC(\RR_+, \RR^\ell))
    \qquad \text{and} \qquad
   \cD_t(\RR_+, \RR^\ell) := \varphi_t^{-1}(\cD(\RR_+, \RR^\ell)) ,
   \qquad t\in\RR_+,
 \]
 respectively, where \ $\varphi_t : D(\RR_+, \RR^\ell) \to D(\RR_+, \RR^\ell)$
 \ is the mapping
 \begin{equation}\label{varphi_t}
  (\varphi_t(z))(s) := z(t \land s) ,
  \qquad z \in D(\RR_+, \RR^\ell), \quad s \in \RR_+ ,
 \end{equation}
 which stops the function \ $z$ \ at \ $t$.
\ It is easy to check that for all \ $t \in \RR_+$, \ $\cC_t(\RR_+, \RR^\ell)$
 \ coincides with the smallest \ $\sigma$-algebra containing all the
 finite-dimensional cylinder sets of the form
 \[
   \big\{ w \in C(\RR_+,\RR^\ell) : (w(t_1), \ldots, w(t_n)) \in A \big\} ,
   \qquad n \in \NN , \quad A \in \cB(\RR^{n\ell}) , \quad
   t_1, \ldots, t_n \in [0,t] ,
 \]
 and then
 \begin{align}\label{C_infty}
  \cC(\RR_+, \RR^\ell)
  = \sigma\Biggl(\bigcup_{t\in\RR_+} \cC_t(\RR_+, \RR^\ell)\Biggr) ,
 \end{align}
 see, e.g., Problem 2.4.2 in Karatzas and Shreve \cite{KarShr}.
Similarly, for all \ $t\in\RR_+$, \ $\cD_t(\RR_+, \RR^\ell)$ \ coincides with
 the smallest \ $\sigma$-algebra containing all the finite-dimensional
 cylinder sets of the form
 \[
   \big\{ y \in D(\RR_+,\RR^\ell) : (y(t_1), \ldots, y(t_n)) \in A \big\} ,
   \qquad n \in \NN , \quad A \in \cB(\RR^{n\ell}) , \quad
   t_1, \ldots, t_n \in [0,t] ,
 \]
 and then
 \begin{align*}
  \cD(\RR_+, \RR^\ell)
  = \sigma\Biggl(\bigcup_{t\in\RR_+} \cD_t(\RR_+, \RR^\ell)\Biggr) ,
 \end{align*}
 hence \ $\cD_t(\RR_+, \RR^\ell)$ \ coincides with \ $\cD_t^0(\RR^\ell)$ \ in
 Definition VI.1.1 in Jacod and Shiryaev \cite{JacShi}.
Finally, let us equip \ $M(\RR_+ \times U)$ \ with the $\sigma$-algebras
 \ $\cM_t(\RR_+ \times U)$, \ $t \in \RR_+$, \ being the smallest
 $\sigma$-algebra containing all sets of the form
 \[
   \{ \pi \in M(\RR_+ \times U) : \pi([0, s] \times B) \in A \} \quad
   \text{with \ $s \in [0, t]$, \ $B \in \cB(U)$, \ $A \in \cB([0, \infty])$} .
 \]
Note that
 \begin{align}\label{M_infty}
  \cM(\RR_+ \times U)
  = \sigma\Biggl(\bigcup_{t\in\RR_+} \cM_t(\RR_+ \times U)\Biggr) ,
 \end{align}
 since the union of the generator system of the $\sigma$-algebras
 \ $\cM_t(\RR_+ \times U)$, \ $t \in \RR_+$, \ forms a generator system of
 \ $\cM(\RR_+ \times U)$.

\section{Notions of weak and strong solutions}
\label{section_weak_strong}

If \ $(\Omega, \cF, \PP)$ \ is a probability space, then, by \ $\PP$-null sets
 from a sub $\sigma$-algebra \ $\cH\subset \cF$, \ we mean the elements of the
 set
 \[
   \{ A \subset \Omega \,
      : \, \exists B \in \cH \;
        \text{\ such that \ $A \subset B$ \ and \ $\PP(B) = 0$ } \} .
 \]

\begin{Def}\label{Def_weak_solution}
Let \ $n$ \ be a probability measure on \ $(\RR^d, \cB(\RR^d))$.
\ A weak solution of the SDE \eqref{SDE_X_YW} with initial distribution \ $n$
 \ is a tuple
 \ $\bigl( \Omega, \cF, (\cF_t)_{t\in\RR_+}, \PP, \bW, p, \bX \bigr)$, \ where
 \begin{enumerate}
  \item[\textup{(D1)}]
   $(\Omega, \cF, (\cF_t)_{t\in\RR_+}, \PP)$ \ is a filtered probability space
    satisfying the usual hypotheses (i.e., \ $(\cF_t)_{t\in\RR_+}$ \ is right
    continuous and \ $\cF_0$ \ contains all the $\PP$-null sets in \ $\cF$);
  \item[\textup{(D2)}]
   $(\bW_t)_{t\in\RR_+}$ \ is an $r$-dimensional standard
    \ $(\cF_t)_{t\in\RR_+}$-Brownian motion;
  \item[\textup{(D3)}]
   $p$ \ is a stationary $(\cF_t)_{t\in\RR_+}$-Poisson point process on \ $U$
    \ with characteristic measure \ $m$;
  \item[\textup{(D4)}]
   $(\bX_t)_{t\in\RR_+}$ \ is an \ $\RR^d$-valued \ $(\cF_t)_{t\in\RR_+}$-adapted
    c\`adl\`ag process such that
     \begin{itemize}
      \item[\textup{(a)}]
       the distribution of \ $\bX_0$ \ is \ $n$, \
      \item[\textup{(b)}]
       $\PP\left(\int_0^t
                  \big(\|b(s, \bX_s)\| + \|\sigma(s, \bX_s)\|^2 \big)\, \dd s
                 < \infty\right) = 1$,
        \ $t\in\RR_+$,
      \item[\textup{(c)}]
       $\PP\left( \int_0^t
                   \int_{U_0} \|f(s, \bX_s, u)\|^2 \, \dd s \, m(\dd u)
                  < \infty \right)
        = 1$,
        \ $t \in \RR_+$,
      \item[\textup{(d)}]
       $\PP\left( \int_0^t \int_{U_1} \|g(s, \bX_{s-}, u)\| \, N(\dd s, \dd u)
                   < \infty \right)
         = 1$,
        \ $t \in \RR_+$, \ where \ $N(\dd s, \dd u)$ \ is the counting measure
        of \ $p$ \ on \ $\RR_{++} \times U$,
      \item[\textup{(e)}]
       equation \eqref{SDE_X_YW} holds \ $\PP$-a.s., where
        \ $\tN(\dd s, \dd u) := N(\dd s, \dd u) - \dd s \, m(\dd u)$.
    \end{itemize}
 \end{enumerate}
\end{Def}

For the definitions of an $(\cF_t)_{t\in\RR_+}$-Brownian motion and an
 $(\cF_t)_{t\in\RR_+}$-Poisson point process, see, e.g., Ikeda and Watanabe
 \cite[Chapter I, Definition 7.2 and Chapter II, Definition 3.2]{IkeWat}.

In the next remark we point out that the integrals in the SDE \eqref{SDE_X_YW}
 are well-defined under the conditions of Definition \ref{Def_weak_solution}
 and have c\`adl\`ag modifications as functions of \ $t$.

\begin{Rem}\label{intexist}
If conditions (D1), (D2) and (D4)(b) are satisfied, then
 \ $\bigl(\int_0^t \sigma(s, \bX_s) \, \dd \bW_s\bigr)_{t\in\RR_+}$ \ is
 well-defined and has continuous sample paths almost surely, see, Ikeda and
 Watanabe \cite[Chapter II, Definition 1.9]{IkeWat}.
Indeed, \ $(\sigma(t, \bX_t))_{t\in\RR_+}$ \ is $(\cF_t)_{t\in\RR_+}$-adapted
 (since \ $\bX$ \ is $(\cF_t)_{t\in\RR_+}$-adapted and \ $\sigma$ \ is
 measurable), \ $(\sigma(t, \bX_t))_{t\in\RR_+}$ \ is measurable (since \ $\bX$
 \ is measurable, because it has right-continuous paths, see Karatzas and
 Shreve \cite[Remark 1.1.14]{KarShr}, and \ $\sigma$ \ is measurable), and
 \ $\PP\bigl(\int_0^t \|\sigma(s, \bX_s)\|^2 \, \dd s < \infty\bigr) = 1$,
 \ $t \in \RR_+$.

Concerning conditions (D4)(c) and (d), note that the mappings
 \ $\RR_+ \times U_0 \times \Omega \ni (s, u, \omega)
    \mapsto f(s, \bX_{s-}(\omega), u) \in \RR^d$
 \ and \ $\RR_+ \times U_1 \times \Omega \ni (s, u, \omega)
    \mapsto g(s, \bX_{s-}(\omega), u) \in \RR^d$
 \ are \ $(\cF_t)_{t\in\RR_+}$-predictable, see Lemma \ref{Lemma_pred}.

Hence condition (D4)(c) is satisfied if and only if the mapping
 \ $\RR_+ \times U_0 \times \Omega \ni (s, u, \omega)
    \mapsto f(s, \bX_{s-}(\omega), u) \in \RR^d$
 \ is in the (multidimensional version of the) class \ $\bF_p^{2,loc}$ \ defined
 on page 62 in Ikeda and Watanabe \cite{IkeWat}, i.e., if it is
 \ $(\cF_t)_{t\in\RR_+}$-predictable and there exists a sequence
 \ $(\tau_n)_{n\in\NN}$ \ of \ $(\cF_t)_{t\in\RR_+}$-stopping times such that
 \ $\tau_n \uparrow \infty$ \ almost surely as \ $n \to \infty$ \ and
 \begin{equation}\label{fcondm}
   \EE\left( \int_0^{t\wedge\tau_n} \int_{U_0}
              \|f(s, \bX_s, u)\|^2 \, \dd s \, m(\dd u) \right)
   < \infty , \qquad t \in\RR_+ , \quad n \in \NN .
 \end{equation}
Indeed, if (D4)(c) holds then \eqref{fcondm} is satisfied for
 \[
   \tau_n
   := \inf\left\{t \in \RR_+
                 : \int_0^t \int_{U_0}
                    \|f(s, \bX_s, u)\|^2 \, \dd s \, m(\dd u) \geq n\right\}
      \land n ,
   \qquad n \in \NN ,
 \]
 where \ $\tau_n \uparrow \infty$ \ almost surely as \ $n \to \infty$.
\ On the other hand, \eqref{fcondm} implies
 \ $\PP\left( \int_0^{t\wedge\tau_n} \int_{U_0}
               \|f(s, \bX_s, u)\|^2 \, \dd s \, m(\dd u) < \infty \right)
    = 1$
 \ for all \ $t \in \RR_+$ \ and \ $n \in \NN$, \ and hence (D4)(c),
 because \ $\tau_n \uparrow \infty$ \ almost surely as \ $n \to \infty$.

Moreover, if conditions (D1), (D3) and (D4)(c) are satisfied, then the process
 \[
   \left(\int_0^t \int_{U_0}
          f(s, \bX_{s-}, u) \, \tN(\dd s, \dd u)\right)_{t\in\RR_+}
 \]
 is well-defined and has c\`adl\`ag sample paths almost surely.
Indeed, for each \ $n \in \NN$,
 \[
   \left(\int_0^{t\land\tau_n} \int_{U_0}
          f(s, \bX_{s-}, u) \, \tN(\dd s, \dd u)\right)_{t\in\RR_+}
   = \left(\int_0^t \int_{U_0}
            \bbone_{[0,\tau_n]}(s)
            f(s, \bX_{s-}, u) \, \tN(\dd s, \dd u)\right)_{t\in\RR_+} ,
 \]
 see page 63 in Ikeda and Watanabe \cite{IkeWat}.
The integrand
 \ $\RR_+ \times U_0 \times \Omega \ni (s, u, \omega)
          \mapsto \bbone_{[0,\tau_n]}(s) f(s, \bX_{s-}(\omega), u) \in \RR^d$
 \ belongs to the (multidimensional version of the) class \ $\bF_p^2$ \ defined
 on page 62 in Ikeda and Watanabe \cite{IkeWat}, hence the process on the
 right hand side is a square integrable \ $(\cF_t)_{t\in\RR_+}$-martingale, see
 page 63 in Ikeda and Watanabe \cite{IkeWat}.
By Theorem 1.3.13 in Karatzas and Shreve \cite{KarShr}, this process has a
 c\`adl\`ag modification.
Here we point out that for using this theorem, we need completeness and right
 continuity of the filtration \ $(\cF_t)_{t\in\RR_+}$.
\ Further, we also obtain
 \[
   \int_0^{t\land\tau_n} \int_{U_0} f(s, \bX_{s-}, u) \, \tN(\dd s, \dd u)
   \as \int_0^t \int_{U_0} f(s, \bX_{s-}, u) \, \tN(\dd s, \dd u) \qquad
   \text{as \ $n \to \infty$}
 \]
 for all \ $t \in \RR_+$, \ since \ $\tau_n \uparrow \infty$ \ almost surely
 as \ $n \to \infty$.

Recalling that the mapping
 \ $\RR_+ \times U_1 \times \Omega \ni (s, u, \omega)
    \mapsto g(s, \bX_{s-}(\omega), u) \in \RR^d$
 \ is \ $(\cF_t)_{t\in\RR_+}$-predictable, condition (D4)(d) is satisfied if and
 only if the mapping
 \ $\RR_+ \times U_1 \times \Omega \ni (s, u, \omega)
    \mapsto g(s, \bX_{s-}(\omega), u) \in \RR^d$
 \ is in the (multidimensional version of the) class \ $\bF_p$ \ defined on
 page 61 in Ikeda and Watanabe \cite{IkeWat}.

Further, if conditions (D1), (D3) and (D4)(d) are satisfied, then, by
 definition, the process
 \[
   \left(\int_0^t \int_{U_1}
          g(s, \bX_{s-}, u) \, N(\dd s, \dd u)\right)_{t\in\RR_+}
   = \left(\sum_{s\in(0,t]\cap D(p)}
            g(s, \bX_{s-}, p(s)) \bbone_{U_1}(p(s)) \right)_{t\in\RR_+}
 \]
 is well-defined and has c\`adl\`ag sample paths, where \ $D(p)$ \ is the
 domain of \ $p$ \ (being a countable subset of \ $\RR_{++}$).
Indeed, for each \ $\omega \in \Omega$, \ by definition, the mappings
 \begin{align*}
  &\RR_+ \ni t \mapsto
    \sum_{s\in(0,t]\cap D(p)(\omega)}
     g(s, \bX_{s-}(\omega), p(s)(\omega)) \bbone_{U_1}(p(s)(\omega)) ,  \\
  &\RR_+ \ni t \mapsto
    \sum_{s\in(0,t)\cap D(p)(\omega)}
     g(s, \bX_{s-}(\omega), p(s)(\omega)) \bbone_{U_1}(p(s)(\omega))
 \end{align*}
 are right and left continuous, respectively.
\proofend
\end{Rem}

\begin{Rem}\label{Rem_sum_finite}
If \ $m(U_1) < \infty$, \ then condition (D4)(d) is satisfied
 automatically, since then \ $\EE(N((0,t] \times U_1) = t m(U_1) < \infty$
 \ implies \ $\PP(N((0,t] \times U_1) < \infty) = 1$, \ and hence
 \ $\int_0^t \int_{U_1} \|g(s, \bX_{s-}, u)\| \, N(\dd s, \dd u)
    = \sum_{s\in(0,t]\cap D(p)} \|g(s, \bX_{s-}, p(s))\| \bbone_{U_1}(p(s))$
 \ is a finite sum with probability one.
\proofend
\end{Rem}

\begin{Rem}\label{dRM_YW}
Note that if conditions (D1)--(D3) are satisfied, then \ $\bW$ \ and \ $p$
 \ are automatically independent according to Theorem 6.3 in Chapter II of
 Ikeda and Watanabe \cite{IkeWat}, since the intensity measure
 \ $\dd s \, m(\dd u)$ \ of \ $p$ \ is deterministic.

Moreover, if \ $\bigl( \Omega, \cF, (\cF_t)_{t\in\RR_+}, \PP, \bW, p, \bX \bigr)$
 \ is a weak solution of the SDE \eqref{SDE_X_YW}, then \ $\cF_0$, \ $\bW$
 \ and \ $p$ \ are mutually independent, and hence \ $\bX_0$, \ $\bW$ \ and
 \ $p$ \ are mutually independent as well.
Indeed, the conditional joint charateristic function of \ $\bW$ \ and the
 counting measure of \ $p$ \ with respect to \ $\cF_0$ \ equals to the product
 of the (unconditional) charateristic functions of \ $\bW$ \ and the counting
 measure of \ $p$, \ see equation (6.12) in Chapter II of Ikeda and Watanabe
 \cite{IkeWat} applied with \ $X = \bW$ \ and \ $s = 0$, \ and then one can
 use Lemma 2.6.13 in Karatzas and Shreve \cite{KarShr}.
Since \ $\bX_0$ \ is measurable with respect to \ $\cF_0$ \ due to (D4), we
 have the mutual independence of \ $\bX_0$, \ $\bW$ \ and \ $p$.

The thinnings \ $p_0$ \ and \ $p_1$ \ of \ $p$ \ onto \ $U_0$ \ and \ $U_1$
 \ are again stationary $(\cF_t)_{t\in\RR_+}$-Poisson point processes on \ $U_0$
 \ and \ $U_1$, \ respectively, and their characteristic measures are the
 restrictions \ $m|_{U_0}$ \ and \ $m|_{U_1}$ \ of \ $m$ \ onto \ $U_0$ \ and
 \ $U_1$, \ respectively (this can be checked calculating their conditional
 Laplace transforms, see Ikeda and Watanabe \cite[page 44]{IkeWat}).

Remark that for any weak solution of the SDE \eqref{SDE_X_YW}, \ $\bX_0$,
 \ the Brownian motion \ $\bW$ \ and the stationary Poisson point processes
 \ $p_0$ \ and \ $p_1$ \ are mutually independent according again to Theorem
 6.3 in Chapter II of Ikeda and Watanabe \cite{IkeWat}.
Indeed, one can argue as before taking into account also that the intensity
 measures of \ $p_0$ \ and \ $p_1$ \ are deterministic, and condition (6.11) of
 this theorem is satisfied, because \ $p_0$ \ and \ $p_1$ \ live on disjoint
 subsets of \ $U$.
\proofend
\end{Rem}

\begin{Def}\label{Def_pathwise_uniqueness}
We say that pathwise uniqueness holds for the SDE \eqref{SDE_X_YW} if whenever
 \ $\bigl( \Omega, \cF, (\cF_t)_{t\in\RR_+}, \PP, \bW, p, \bX \bigr)$ \ and
 \ $\bigl( \Omega, \cF, (\cF_t)_{t\in\RR_+}, \PP, \bW, p, \tbX \bigr)$ \ are weak
 solutions of the SDE \eqref{SDE_X_YW} such that \ $\PP(\bX_0 = \tbX_0) = 1$,
 \ then \ $\PP(\text{$\bX_t = \tbX_t$ \ for all \ $t \in \RR_+$}) = 1$.
\end{Def}

\begin{Rem}
One may also consider the following more strict definition of pathwise
 uniqueness.
Namely, one could say that pathwise uniqueness holds for the SDE
 \eqref{SDE_X_YW} if whenever
 \ $\bigl( \Omega, \cF, (\cF_t)_{t\in\RR_+}, \PP, \bW, p, \bX \bigr)$ \ and
 \ $\bigl( \Omega, \cF, (\tcF_t)_{t\in\RR_+}, \PP, \bW, p, \tbX \bigr)$ \ are
 weak solutions of the SDE \eqref{SDE_X_YW} such that
 \ $\PP(\bX_0 = \tbX_0) = 1$, \ then
 \ $\PP(\bX_t = \tbX_t \; \text{for all} \; t \in \RR_+) = 1$.
\ Note that in this definition we require that \ $\bW$ \ is an
\ $(\cF_t)_{t\in\RR_+}$-Brownian motion and an \ $(\tcF_t)_{t\in\RR_+}$-Brownian
 motion as well, and since it is not necessarily true that \ $\bW$ \ is an
 \ $(\sigma (\cF_t\cup\tcF_t))_{t\in\RR_+}$-Brownian motion, it is not clear
 whether this more strict definition of pathwise uniqueness and the one given
 in \ref{Def_pathwise_uniqueness} are equivalent.
According to Ikeda and Watanabe \cite[Chapter IV, Remark 1.3]{IkeWat}, they
 are equivalent.
We also point out that in our statements and proofs we use pathwise uniqueness
 in the sense of Definition \ref{Def_pathwise_uniqueness}, and we do not use
 the above mentioned equivalence of the two kinds of definitions.
\proofend
\end{Rem}

\begin{Def}\label{Def_uniqueness_in_law}
We say that uniqueness in the sense of probability law holds for the SDE
 \eqref{SDE_X_YW} if whenever
 \ $\bigl( \Omega, \cF, (\cF_t)_{t\in\RR_+}, \PP, \bW, p, \bX \bigr)$ \ and
 \ $\bigl( \tOmega, \tcF, (\tcF_t)_{t\in\RR_+}, \tPP, \tbW, \tp, \tbX \bigr)$
 \ are weak solutions of the SDE \eqref{SDE_X_YW} with the same initial
 distribution, i.e., \ $\PP(\bX_0 \in B) = \tPP(\tbX_0 \in B)$ \ for all
 \ $B \in \cB(\RR^d)$, \ then \ $\PP(\bX \in C) = \tPP(\tbX \in C)$ \ for all
 \ $C \in \cD(\RR_+, \RR^d)$.
\end{Def}

Now we define strong solutions.
Consider the following objects:
 \begin{enumerate}
  \item[(E1)]
    a probability space \ $(\Omega, \cF, \PP)$;
  \item[(E2)]
   an $r$-dimensional standard Brownian motion \ $(\bW_t)_{t\in\RR_+}$;
  \item[(E3)]
   a stationary Poisson point process \ $p$ \ on \ $U$
    \ with characteristic measure \ $m$;
  \item[(E4)]
   a random vector \ $\bxi$ \ with values in \ $\RR^d$, \ independent of
    \ $\bW$ \ and \ $p$.
 \end{enumerate}

\begin{Rem}\label{dRM_strong}
Note that if conditions (E1)--(E4) are satisfied, then \ $\bxi$, \ $\bW$ \ and
 \ $p$ \ are automatically mutually independent according to Remark
 \ref{dRM_YW}.
\proofend
\end{Rem}

Provided that the objects (E1)--(E4) are given, let
 \ $(\cF^{\bxi,\bW\!,\,p}_t)_{t\in\RR_+}$ \ be the augmented filtration generated by
 \ $\bxi$, \ $\bW$ \ and \ $p$, \ i.e., for each \ $t \in \RR_+$,
 \ $\cF^{\bxi,\bW\!,\,p}_t$ \ is the $\sigma$-field generated by
 \ $\sigma(\bxi; \, \bW_s, s\in[0,t]; \, p(s), s\in(0,t]\cap D(p))$ \ and by
 the \ $\PP$-null sets from
 \ $\sigma(\bxi; \, \bW_s, s\in\RR_+; \, p(s), s\in\RR_{++}\cap D(p))$ \ (which
 is similar to the definition in Karatzas and Shreve \cite[page 285]{KarShr}).
One can check that
 \begin{itemize}
  \item
   $(\cF^{\bxi,\bW\!,\,p}_t)_{t\in\RR_+}$ \ satisfies the usual hypotheses;
  \item
   $(\bW_t)_{t\in\RR_+}$ \ is a standard
    $(\cF^{\bxi,\bW\!,\,p}_t)_{t\in\RR_+}$-Brownian motion;
  \item
   $p$ is a stationary $(\cF^{\bxi,\bW\!,\,p}_t)_{t\in\RR_+}$-Poisson point process
    on $U$ with characteristic measure $m$.
 \end{itemize}
Indeed, by Remark \ref{dRM_strong}, \ $\bW$ \ is a standard
 \ $(\sigma(\bxi; \, \bW_s, s\in[0,t]; \,
            p(s), s\in(0,t]\cap D(p)))_{t\in\RR_+}$-Brownian motion, and \ $p$
 \ is a stationary
 $(\sigma(\bxi; \, \bW_s, s\in[0,t]; \,
          p(s), s\in(0,t]\cap D(p)))_{t\in\RR_+}$-Poisson
 point process on \ $U$ \ with characteristic measure \ $m$.
\ Hence, by Theorems 6.4 and 6.5 in Chapter II in Ikeda and Watanabe
 \cite{IkeWat}, \ $(\bW,p)$ \ has the strong Markov property with respect to
 the filtration
 \ $(\sigma(\bxi; \, \bW_s, s\in[0,t]; \,
     p(s), s\in(0,t]\cap D(p)))_{t\in\RR_+}$.
\ Then Proposition 2.7.7 in Karatzas and Shreve \cite{KarShr} yields that the
 augmented filtration \ $(\cF^{\bxi,\bW\!,\,p}_t)_{t\in\RR_+}$ \ satisfies the usual
 hypotheses.
Moreover, the augmentation of \ $\sigma$-fields does not disturb the
 definition of a standard Wiener process and a stationary Poisson point
 process, hence \ $(\bW_t)_{t\in\RR_+}$ \ is a standard
 $(\cF^{\bxi,\bW\!,\,p}_t)_{t\in\RR_+}$-Brownian motion, and \ $p$ \ is a
 stationary $(\cF^{\bxi,\bW\!,\,p}_t)_{t\in\RR_+}$-Poisson point process on \ $U$
 \ with characteristic measure \ $m$.
For the standard Wiener process, see, e.g., Karatzas and Shreve
 \cite[Theorem 2.7.9]{KarShr}.
The main point is to show that \ $\bW_t - \bW_s$ \ is independent of
 \ $\cF^{\bxi,\bW,p}_s$ \ for all \ $s,t\in\RR_+$ \ with \ $s < t$, \ and
 \ $p(t) - p(s)$ \ is independent of \ $\cF^{\bxi,\bW,p}_s$ \ for all
 \ $s, t \in D(p)$ \ with \ $s < t$, \ detailed as follows (in order to shed
 some light what is going on behind).
\ Let \ $s, t \in \RR_+$ \ with \ $s < t$, \ and \ $F \in \cF^{\bxi,\bW\!,\,p}_s$.
\ Then, by Problem 2.7.3 in Karatzas and Shreve \cite{KarShr}, there exists
 \ $\tF \in \sigma(\bxi; \, \bW_u, u\in[0,s]; \, p(u), u\in(0,s]\cap D(p))$
 \ such that \ $F \Delta \tF$ \ is a \ $\PP$-null set from
 \ $\sigma(\bxi; \, \bW_u, u\in\RR_+; \, p(u), u\in\RR_{++}\cap D(p))$,
 \ where \ $F \Delta \tF$ \ denotes the symmetric difference of \ $F$ \ and
 \ $\tF$.
\ Using that
 \[
   \PP(A)
   = \PP(B) + \PP(A \cap(\Omega\setminus B))
     - \PP((\Omega\setminus A) \cap B) ,
   \qquad A, B \in \cF ,
 \]
 we get for all \ $K \in \cB(\RR^r)$,
 \begin{align*}
  &\PP(\{ \bW_t - \bW_s \in K \} \cap F)\\
  &\quad
   = \PP(\{ \bW_t - \bW_s \in K \} \cap \tF)\\
  &\quad\quad
     + \PP( \{ \bW_t - \bW_s \in K \} \cap F
            \cap (\{ \bW_t - \bW_s \not\in K \}
                  \cup (\Omega\setminus \tF)) ) \\
  &\quad\quad
     - \PP( (\{ \bW_t - \bW_s \not\in K \} \cup (\Omega\setminus F) )
            \cap \{ \bW_t - \bW_s \in K \} \cap \tF ) \\
  &\quad
   = \PP(\{ \bW_t - \bW_s\in K \} \cap \tF)
     + \PP(\{ \bW_t - \bW_s\in K \} \cap F \cap (\Omega\setminus \tF) ) \\
  &\quad\quad
     - \PP(\{ \bW_t - \bW_s \in K \} \cap (\Omega\setminus F) \cap \tF)  \\
  &\quad
   = \PP(\{ \bW_t - \bW_s\in K \} \cap \tF)
   = \PP(\bW_t - \bW_s\in K) \PP(\tF)
    = \PP(\bW_t - \bW_s\in K) \PP(F) ,
 \end{align*}
 where the last but one step follows from the independence of
 \ $\bW_t - \bW_s$ \ and \ $\tF$.
\ A similar argument shows the independence of \ $p(t) - p(s)$ \ and \ $F$.

\begin{Def}\label{Def_strong_solution2}
Suppose that the objects \textup{(E1)--(E4)} are given.
A strong solution of the SDE \eqref{SDE_X_YW} on \ $(\Omega, \cF, \PP)$ \ and
 with respect to the standard Brownian motion \ $\bW$, \ the stationary
 Poisson point process \ $p$ \ and initial value \ $\bxi$, \ is an
 \ $\RR^d$-valued \ $(\cF^{\bxi,\bW\!,\,p}_t)_{t\in\RR_+}$-adapted c\`{a}dl\`{a}g
 process \ $(\bX_t)_{t \in \RR_+}$ \ with \ $\PP(\bX_0 = \bxi) = 1$ \ satisfying
 \textup{(D4)(b)--(d)}.
\end{Def}

Clearly, if \ $(\bX_t)_{t \in \RR_+}$ \ is a strong solution, then
 \ $\bigl( \Omega, \cF, (\cF^{\bxi,\bW\!,\,p}_t)_{t\in\RR_+}, \PP, \bW, p,
           \bX \bigr)$
 \ is a weak solution with initial distribution being the distribution of
 \ $\bxi$.

\section{Proof of Theorem \ref{Thm_pathwise_1}}
\label{section_proof1}

Our presentation as follows is a generalization of the one given in Section
 5.3.D in Karatzas and Shreve \cite{KarShr}.

Let us consider a weak solution
 \ $\bigl( \Omega, \cF, (\cF_t)_{t\in\RR_+}, \PP, \bW, p, \bX \bigr)$ \ of
 the SDE \eqref{SDE_X_YW} with initial distribution \ $n$ \ on
 \ $(\RR^d, \cB(\RR^d))$.
\ Then \ $\PP(\bX_0 \in B) = n(B)$, \ $B \in \cB(\RR^d)$.
\ We put \ $\bY_t:= \bX_t - \bX_0$ \ for \ $t \in \RR_+$, \ and we regard the
 solution \ $\bX$ \ as consisting of four parts: \ $\bX_0$, \ $\bW$, \ $p$
 \ and \ $\bY$.
\ Let us consider the product space
 \begin{align}\label{space_theta}
   \Theta
   := \RR^d \times C(\RR_+, \RR^r) \times M(\RR_+\times U)
      \times D(\RR_+, \RR^d)
 \end{align}
 equipped with the Borel $\sigma$-algebra
 \[
   \cB(\Theta)
   = \cB(\RR^d) \otimes \cC(\RR_+, \RR^r) \otimes \cM(\RR_+\times U)
     \otimes \cD(\RR_+, \RR^d) ,
 \]
 see, e.g., Dudley \cite[Proposition 4.1.7]{Dud}.
The quadruplet \ $(\bX_0, \bW, p, \bY)$ \ induce the probability measure \ $P$
 \ on \ $(\Theta, \cB(\Theta))$ \ according to the prescription
 \begin{align}\label{help10}
  P(A) := \PP[(\bX_0, \bW, p, \bY) \in A] , \qquad A \in \cB(\Theta) .
 \end{align}
We denote by \ $\theta = (\bx, w, \pi, y)$ \ a generic element of \ $\Theta$.
\ The marginal of \ $P$ \ on the \ $\bx$-coordinate of \ $\theta$ \ is the
 probability measure \ $n$ \ on \ $(\RR^d, \cB(\RR^d))$, \ the marginal on the
 \ $w$-coordinate is an $r$-dimensional Wiener measure \ $P_{\bW\!,\,r}$ \ on
 \ $(C(\RR_+, \RR^r),\cC(\RR_+, \RR^r))$, \ the marginal on the
 \ $\pi$-coordinate is the distribution \ $P_{U,m}$ \ on
 \ $(M(\RR_+ \times U), \cM(\RR_+ \times U))$ \ of a stationary Poisson point
 process \ $p$ \ on \ $U$ \ with characteristic measure \ $m$.
\ Moreover, the distribution of the triplet \ $(\bx, w, \pi)$ \ under \ $P$
 \ is the product measure \ $n \times P_{\bW\!,\,r} \times P_{U,m}$ \ because
 \ $\bX_0$ \ is $\cF_0$-measurable and \ $\bW$, \ $p$ \ and \ $\cF_0$ \ are
 independent, see Remark \ref{dRM_YW}.
Furthermore, \ $\PP(\bY_0 = \bzero) = 1$.

The product space \ $\Theta$ \ defined in \eqref{space_theta} is a complete,
 separable metric space, since \ $\RR^d$ \ is a complete, separable metric
 space with the usual Euclidean metric, \ $C(\RR_+,\RR^r)$ \ is a complete,
 separable metric space with a metric inducing the local uniform topology
 (see, e.g., Jacod and Shiryaev \cite[Section VI.1a]{JacShi}),
 \ $D(\RR_+, \RR^d)$ \ is a complete, separable metric space with a metric
 inducing the so-called Skorokhod topology (see, e.g., Jacod and Shiryaev
 \cite[Theorem VI.1.14]{JacShi}), and the vague topology on the space
 \ $M(\RR_+ \times U)$ \ of all point measures on \ $\RR_+ \times U$ \ is
 metrizable as a complete, separable metric space (see, e.g., Resnick
 \cite[Proposition 3.17, page 147]{Res}).
Hence there exists a regular conditional probability for \ $\cB(\Theta)$
 \ given \ $(\bx, w, \pi)$, \ by an application of Karatzas and Shreve
 \cite[Chapter 5, Theorem 3.19]{KarShr} with the random variable
 \ $\Theta \ni (\bx, w, \pi, y) \mapsto (\bx, w, \pi)$.
\ We shall be interested in conditional probabilities of sets in
 \ $\cB(\Theta)$ \ only of the form
 \ $\RR^d \times C(\RR_+,\RR^r) \times M(\RR_+ \times U) \times F$, \ where
 \ $F\in \cD(\RR_+,\RR^d)$.
\ Consequently, with a slight abuse of notation, there exists a function
 \begin{align}\label{help11}
  Q : \RR^d \times C(\RR_+, \RR^r) \times M(\RR_+\times U)
      \times \cD(\RR_+, \RR^d)
      \to [0, 1]
 \end{align}
 enjoying the following properties:
 \begin{enumerate}
  \item[(R1)]
   for each \ $\bx \in \RR^d$, \ $w \in C(\RR_+, \RR^r)$ \ and
    \ $\pi \in M(\RR_+ \times U)$, \ the set function
    \ $\cD(\RR_+, \RR^d) \ni F \mapsto Q(\bx, w, \pi, F)$ \ is a
    probability measure on \ $(D(\RR_+, \RR^d), \cD(\RR_+, \RR^d))$;
  \item[(R2)]
   for each \ $F \in \cD(\RR_+, \RR^d)$, \ the mapping
    \ $\RR^d \times C(\RR_+, \RR^r) \times M(\RR_+ \times U) \ni (\bx, w, \pi)
       \mapsto Q(\bx, w, \pi, F)$
    \ is
    \ $\cB(\RR^d) \otimes \cC(\RR_+, \RR^r)
       \otimes \cM(\RR_+\times U) / \cB([0, 1])$-measurable;
  \item[(R3)]
   for each
    \ $G \in \cB(\RR^d) \otimes \cC(\RR_+, \RR^r)
             \otimes \cM(\RR_+\times U)$
    \ and \ $F \in \cD(\RR_+, \RR^d)$, \ we have
   \[
     P(G \times F)
     = \int_G
        Q(\bx, w, \pi, F) \, n(\dd\bx) \, P_{\bW\!,\,r}(\dd w) \,
        P_{U,m}(\dd\pi) .
   \]
 \end{enumerate}
We can call \ $Q(\bx, w, \pi, \cdot)$ \ as the regular conditional
 probability for \ $\cD(\RR_+, \RR^d)$ \ given \ $(\bx, w, \pi)$.

Let us now consider two weak solutions
 \ $\bigl( \Omega^{(i)}, \cF^{(i)}, (\cF_t^{(i)})_{t\in\RR_+}, \PP^{(i)}, \bW^{(i)},
           p^{(i)}, \bX^{(i)} \bigr)$,
 \ $i \in \{1, 2\}$ \ of the SDE \eqref{SDE_X_YW} with the same initial
 distribution \ $n$ \ on \ $(\RR^d, \cB(\RR^d))$, \ thus
 \[
   \PP^{(1)}[\bX^{(1)}_0 \in B] = \PP^{(2)}[\bX^{(2)}_0 \in B] = n(B), \qquad
   B \in \cB(\RR^d) .
 \]
According to \eqref{help10}, let
 \begin{align*}
  P_i(A) := \PP^{(i)}[(\bX^{(i)}_0, \bW^{(i)}, p^{(i)}, \bY^{(i)}) \in A] , \qquad
  A \in \cB(\Theta) , \qquad i \in \{1, 2\} ,
 \end{align*}
 and, as explained before, there exist functions
 \begin{align}\label{help16}
  Q_i : \RR^d \times C(\RR_+, \RR^r) \times M(\RR_+\times U)
        \times \cD(\RR_+, \RR^d)
        \to [0, 1], \qquad i \in \{1, 2\} ,
 \end{align}
 enjoying the properties (R1)--(R3).

First, we bring the two triplets \ $(\bX^{(i)}, \bW^{(i)}, p^{(i)})$,
 \ $i \in \{1, 2\}$, \ together on the same, canonical space, while preserving
 the joint distribution of the coordinates within each triplet.
Let \ $\Omega := \Theta \times D(\RR_+, \RR^d)$
 \ equipped with the $\sigma$-algebra \ $\cF$, \ which is the completion of
 the product $\sigma$-algebra \ $\cB(\Theta) \otimes \cD(\RR_+, \RR^d)$
 \ by the collection \ $\cN$ \ of null sets under the probability measure
 \begin{align}\label{help12}
  \begin{split}
   \PP_{1,2}(A) := &\int_{\RR^d \times C(\RR_+, \RR^r) \times M(\RR_+ \times U)}
               \Bigg(\int_{D(\RR_+, \RR^d) \times D(\RR_+, \RR^d)}
                      \bbone_A (\bx, w, \pi, y^{(1)}, y^{(2)}) \\
                    &\qquad\qquad\qquad Q_1(\bx, w, \pi, \dd y^{(1)}) \,
                            Q_2(\bx, w, \pi, \dd y^{(2)})\Bigg)
                     n(\dd\bx) \, P_{\bW\!,\,r}(\dd w) \, P_{U,m}(\dd\pi)
  \end{split}
 \end{align}
 for \ $A \in \cB(\Theta) \otimes \cD(\RR_+, \RR^d)$, \ where we have
 denoted by \ $(\bx, w, \pi, y^{(1)}, y^{(2)})$ \ a generic element of
 \ $\Omega$, \ and then we extend \ $\PP_{1,2}$ \ to \ $\cF$.
\ Especially, for all
 \ $G \in \cB(\RR^d) \otimes \cC(\RR_+, \RR^r)\otimes \cM(\RR_+\times U)$
 \ and \ $F_1, F_2 \in \cD(\RR_+,\RR^d)$,
 \[
   \PP_{1,2}(G\times F_1\times F_2)
     = \int_G
        Q_1(\bx, w, \pi, F_1) \, Q_2(\bx, w, \pi, F_2) \, n(\dd\bx) \,
        P_{\bW\!,\,r}(\dd w) \, P_{U,m}(\dd\pi) .
 \]
In order to endow \ $(\Omega, \cF, \PP_{1,2})$ \ with a filtration that
 satisfies the usual conditions, for each \ $t\in\RR_+$, \ we take
 \ $\cG_t := \sigma(f_{s,B} : s \in [0, t], \, B \in \cB(U))$, \ where the
 mapping
 \ $f_{s,B}: \Omega \to \RR^d\times\RR^r\times[0,\infty]\times\RR^d\times\RR^d$
 \ is defined by
 \begin{align*}
  f_{s,B}(\bx, w, \pi, y^{(1)},y^{(2)})
  := \bigr(\bx, w_s, \pi( [0,s]\times B), y^{(1)}_s , y^{(2)}_s\bigr) ,
  \qquad (\bx, w, \pi, y^{(1)},y^{(2)}) \in \Omega ,
 \end{align*}
 and put
 \begin{align*}
  \widetilde{\cG}_t := \sigma(\cG_t \cup \cN), \qquad
   \cF_t := \widetilde{\cG}_{t+}
         := \bigcap_{\vare>0} \widetilde \cG_{t+\vare} , \qquad
  t \in \RR_+ .
 \end{align*}
We note that for each \ $t \in \RR_+$,
 \[
   \cG_t = \widehat\cG_t
    = \cB(\RR^d) \otimes \cC_t(\RR_+,\RR^r) \otimes \cM_t(\RR_+\times U)
      \otimes \cD_t(\RR_+,\RR^d)\otimes \cD_t(\RR_+,\RR^d) ,
 \]
 where
 \ $\widehat{\cG}_t
    := \sigma(\widehat{f}_{s,B} : s \in [0, t], \, B \in \cB(U))$,
 \ and the mapping \ $\widehat{f}_{s,B}: \Omega \to \Omega$ \ is defined by
 \begin{align*}
  \widehat{f}_{s,B}(\bx, w, \pi, y^{(1)},y^{(2)})
  :=\bigr(\bx, (w_{t\wedge s})_{t\in\RR_+}, \pi|_{[0,s]\times B},
          (y^{(1)}_{t\wedge s})_{t\in\RR_+} , (y^{(2)}_{t\wedge s})_{t\in\RR_+} \bigr)
 \end{align*}
 for \ $(\bx, w, \pi, y^{(1)},y^{(2)}) \in \Omega$.
\ Indeed, for all \ $t\in\RR_+$, \ by definition,
 the \ $\sigma$-algebra \ $\cG_t$ \ coincides with the
 \ $\sigma$-algebra generated by the sets
 \begin{align*}
  E_1&\times \{w\in C(\RR_+,\RR^r) : w(s)\in E_2\}
      \times \{\pi\in M(\RR_+\times U): \pi([0,s]\times B)\in E_3\}\\
     & \times \{y^{(1)}\in D(\RR_+,\RR^d) : y^{(1)}(s)\in E_4\}
       \times \{y^{(2)}\in D(\RR_+,\RR^d) : y^{(2)}(s)\in E_5\}
 \end{align*}
 for \ $s \in [0,t]$, \ $B \in \cB(U)$, \ $E_1 \in \cB(\RR^d)$,
 \ $E_2 \in \cB(\RR^{r})$, \ $E_3 \in \cB([0,\infty])$ \ and
 \ $E_4, E_5 \in \cB(\RR^{d})$.
\ Moreover, as in Problem 2.4.2 in Karatzas and Shreve \cite{KarShr}, the
 \ $\sigma$-algebra \ $\widehat\cG_t$ \ coincides with the \ $\sigma$-algebra
 generated by the sets
 \begin{align*}
  E_1&\times \{w \in C(\RR_+,\RR^r)
               : (w(t_{1,1}\wedge s), \ldots, w(t_{1,n_1}\wedge s))\in E_2\} \\
     &\times \{\pi \in M(\RR_+\times U)
               : (\pi([0,t_{2,1}\wedge s]\times B_1), \ldots,
                  \pi([0,t_{2,n_2}\wedge s]\times B_{n_2})) \in E_3\} \\
     &\times \{y^{(1)} \in D(\RR_+,\RR^d)
               : (y^{(1)}(t_{3,1}\wedge s), \ldots, y^{(1)}(t_{3,n_3}\wedge s))
                 \in E_4\} \\
     &\times \{y^{(2)} \in D(\RR_+,\RR^d)
               : (y^{(2)}(t_{4,1}\wedge s), \ldots, y^{(2)}(t_{4,n_4}\wedge s))
                 \in E_5\}
 \end{align*}
 for \ $s \in [0,t]$, \ $t_{i,j} \in \RR_+$, \ $i \in \{1,2,3,4\}$,
 \ $j \in \{1, \ldots, n_i\}$, \ $B_1, \ldots, B_{n_2} \in \cB(U)$,
 \ $E_1 \in \cB(\RR^d)$, \ $E_2 \in \cB(\RR^{rn_1})$,
 \ $E_3 \in \cB([0,\infty]^{n_2})$, \ $E_4 \in \cB(\RR^{dn_3})$ \ and
 \ $E_5 \in \cB(\RR^{dn_4})$.
\ Since for any stochastic process \ $(\xi_t)_{t\in\RR_+}$,
 \begin{align}\label{help_merhetoseg}
  \sigma (\xi_t: \, t\in[0,s])
  = \sigma((\xi_{t_1},\ldots,\xi_{t_n})
           : t_i\in[0,s], \, i\in\{1,\ldots,n\}, \, n\in\NN) ,
  \qquad s \in \RR_+ ,
 \end{align}
 we get \ $\widehat\cG_t = \cG_t$, \ $t\in\RR_+$.

The \ $\pi$-coordinate process on \ $\Omega$ \ induces a point process
 \ $p_\pi$ \ on \ $U$ \ with characteristic measure \ $m$ \ in a natural way,
 since, as it was recalled, there is a bijection between the set of point
 functions on \ $U$ \ and the set of point measures \ $\pi$ \ on
 \ $\RR_+ \times U$ \ with \ $\pi(\{0\} \times U) = 0$ \ and
 \ $\pi(\{t\} \times U) \leq 1$, \ $t \in \RR_{++}$, \ and
 \[
   \PP_{1,2}\big(\big\{(\bx,w,\pi,y^{(1)},y^{(2)}) \in \Omega
                 : \pi(\{0\}\times U)=0, \;
                   \pi(\{t\} \times U)\leq1, \;
                   t\in\RR_{++}\big\}\big)
   = 1 ,
 \]
 which follows from \eqref{help12} using that \ $P_{U,m}$ \ is the distribution
 on \ $(M(\RR_+ \times U), \cM(\RR_+ \times U))$ \ of a stationary Poisson
 point process on \ $U$ \ with characteristic measure \ $m$ \ implying that
 \begin{align*}
  P_{U,m}\big(\big\{ \pi \in M(\RR_+ \times U)
                    : \pi(\{0\}\times U)=0, \;
                      \pi(\{t\} \times U)\leq1, \;
                      t\in\RR_{++} \big\}\big)
  = 1 .
 \end{align*}

Next we check that
 \ $\bigl( \Omega, \cF, (\cF_t)_{t\in\RR_+}, \PP_{1,2}, w, p_\pi,
           (\bx + y^{(i)}_t)_{t\in\RR_+} \bigr)$,
 \ $i \in\{1, 2\}$, \ are weak solutions of the SDE \eqref{SDE_X_YW} with the
 same initial distribution \ $n$.
\ Using the definitions of \ $P_i$, \ $i \in \{1, 2\}$, \ $\PP_{1,2}$, \ (R1)
 and (R3) we get
 \begin{equation}\label{PP}
  \PP_{1,2}[\omega = (\bx, w, \pi, y^{(1)}, y^{(2)}) \in \Omega
      : (\bx, w, \pi, y^{(i)}) \in A]
  = \PP^{(i)}[(\bX_0^{(i)}, \bW^{(i)}, p^{(i)}, \bY^{(i)}) \in A]
 \end{equation}
 for all \ $A \in \cB(\Theta)$ \ and \ $i \in \{1, 2\}$.
\ Indeed, with \ $i = 1$,
 \ $G \in \cB(\RR^d) \otimes \cC(\RR_+, \RR^r)
          \otimes \cM(\RR_+\times U)$
 \ and \ $F \in \cD(\RR_+, \RR^d)$, \ by Fubini theorem,
 \begin{multline*}
  \PP_{1,2}[\omega = (\bx, w, \pi, y^{(1)}, y^{(2)}) \in \Omega
       : (\bx, w, \pi, y^{(1)}) \in G \times F] \\
  \begin{aligned}
   &= \int_{\{\omega\in\Omega:(\bx,w,\pi,y^{(1)})\in G\times F\}}
       Q_1(\bx, w, \pi, \dd y^{(1)}) \, Q_2(\bx, w, \pi, \dd y^{(2)})
       \, n(\dd\bx) \, P_{\bW,r}(\dd w) \, P_{U,m}(\dd\pi) \\
   &= \int_G
       Q_1(\bx, w, \pi, F) \, Q_2(\bx, w, \pi, D(\RR_+, \RR^d))
       \, n(\dd\bx) \, P_{\bW,r}(\dd w) \, P_{U,m}(\dd\pi) \\
   &= \int_G
       Q_1(\bx, w, \pi, F) \, n(\dd\bx) \, P_{\bW,r}(\dd w)
       \, P_{U,m}(\dd\pi)
    = P_1(G\times F)  \\
   &= \PP^{(1)}[(\bX_0^{(1)}, \bW^{(1)}, p^{(1)}, \bY^{(1)}) \in G \times F] .
  \end{aligned}
 \end{multline*}
So the distribution of \ $(\bx + y^{(i)}, w, p_\pi)$ \ under \ $\PP_{1,2}$ \ is
 the same as the distribution of
 \ $(\bX^{(i)}_0 + \bY^{(i)}, \bW^{(i)}, p^{(i)}) = (\bX^{(i)}, \bW^{(i)}, p^{(i)})$
 \ under \ $\PP^{(i)}$.
\ Due to the definition of a weak solution, under \ $\PP^{(i)}$, \  $\bW^{(i)}$
 \ is an $r$-dimensional standard \ $(\cF_t^{(i)})_{t\in\RR_+}$-Brownian motion,
 and \ $p^{(i)}$ \ is a stationary $(\cF_t^{(i)})_{t\in\RR_+}$-Poisson point
 process on \ $U$ \ with characteristic measure \ $m$.
\ Consequently, by the definition of \ $(\cG_t)_{t\in\RR_+}$ \ (which is nothing
 else but the natural filtration corresponding to the coordinate processes),
 under \ $\PP_{1,2}$, \ the $w$-coordinate process is an $r$-dimensional
 standard \ $(\cG_t)_{t\in\RR_+}$-Brownian motion, the process \ $p_\pi$ \ is a
 stationary $(\cG_t)_{t\in\RR_+}$-Poisson point process on \ $U$ \ with
 characteristic measure \ $m$, \ and \ $(\bx+y^{(i)}_t)_{t\in\RR_+}$ \ is
 \ $(\cG_t)_{t\in\RR_+}$-adapted, \ $i\in\{1, 2\}$.
\ Further, the same is true if we replace the filtration \ $(\cG_t)_{t\in\RR_+}$
 \ by \ $(\cF_t)_{t\in\RR_+}$, \ see, Lemma \ref{Lemma_filtration}.
Note also that the filtration \ $(\cF_t)_{t\in\RR_+}$ \ satisfies the usual
 conditions.
All in all, for each \ $i\in\{1,2\}$, \ the tuple
 \ $\bigl( \Omega, \cF, (\cF_t)_{t\in\RR_+}, \PP_{1,2}, w, p_\pi,
           (\bx + y^{(i)}_t)_{t\in\RR_+} \bigr)$ \ satisfies (D1)--(D3).

Hence it remains to check that, for each \ $i\in\{1,2\}$, \ the tuple
 \ $\bigl( \Omega, \cF, (\cF_t)_{t\in\RR_+}, \PP_{1,2}, w, p_\pi,
           (\bx + y^{(i)}_t)_{t\in\RR_+} \bigr)$ \ satisfies (D4).
For each \ $i \in \{1, 2\}$, \ let us apply Lemma \ref{Lemma_distribution2}
 with the following choices
  \[
    \bigl( \Omega^{(i)}, \cF^{(i)}, (\cF_t^{(i)})_{t\in\RR_+}, \PP^{(i)}, \bW^{(i)},
           p^{(i)}, \bX^{(i)} \bigr)
 \]
 and
 \[
   \bigl( \Omega, \cF, (\cF_t)_{t\in\RR_+}, \PP_{1,2}, w, p_\pi,
          (\bx + y^{(i)}_t)_{t\in\RR_+} \bigr) .
 \]
Since
 \ $\bigl( \Omega^{(i)}, \cF^{(i)}, (\cF_t^{(i)})_{t\in\RR_+}, \PP^{(i)}, \bW^{(i)},
           p^{(i)}, \bX^{(i)} \bigr)$
 \ is a weak solution of the SDE \eqref{SDE_X_YW} with initial distribution
 \ $n$, \ the tuple
 \ $\bigl( \Omega^{(i)}, \cF^{(i)}, (\cF_t^{(i)})_{t\in\RR_+}, \PP^{(i)}, \bW^{(i)},
           p^{(i)}, \bX^{(i)} \bigr)$
 \ satisfies \textup{(D1)}--\textup{(D4)}.
Further, as it was explained before, the tuple
 \ $\bigl( \Omega, \cF, (\cF_t)_{t\in\RR_+}, \PP_{1,2}, w,
           p_\pi, (\bx + y^{(i)}_t)_{t\in\RR_+} \bigr)$
 \ satisfies \textup{(D1)}--\textup{(D3)}, the process
 \ $(\bx + y^{(i)}_t)_{t\in\RR_+}$ \ is adapted to the filtration
 \ $(\cF_t)_{t\in\RR_+}$, \ and the distribution of
 \ $(\bX^{(i)}, \bW^{(i)}, p^{(i)})$ \ under \ $\PP^{(i)}$ \ is the same as the
 distribution of \ $(\bx + y^{(i)}, w, p_\pi)$ \ under \ $\PP_{1,2}$.
\ Then Lemma  \ref{Lemma_distribution2} yields that
 the tuple
 \ $\bigl( \Omega, \cF, (\cF_t)_{t\in\RR_+}, \PP_{1,2}, w,
           p_\pi, (\bx + y^{(i)}_t)_{t\in\RR_+} \bigr)$
 \ satisfies  \textup{(D4)(a)--(d)} and the distribution of
 \begin{align*}
  \biggl(&\bX_t^{(i)} - \bX_0^{(i)}
          - \int_0^t b(s, \bX_s^{(i)}) \, \dd s
          - \int_0^t \sigma(s, \bX_s^{(i)}) \, \dd \bW_s^{(i)} \\
         &- \int_0^t \int_{U_0} f(s, \bX_{s-}^{(i)}, u) \, \tN^{(i)}(\dd s, \dd u)
          - \int_0^t \int_{U_1}
             g(s, \bX_{s-}^{(i)}, u) \, N^{(i)}(\dd s, \dd u) \biggr)_{t\in\RR_+}
 \end{align*}
 on \ $(D(\RR_+, \RR^d), \cD(\RR_+, \RR^d))$ \ under \ $\PP^{(i)}$ \ is
 the same as the distribution of
 \begin{align*}
  \biggl(&y^{(i)}_t - y^{(i)}_0
          - \int_0^t b(s, \bx+y^{(i)}_s) \, \dd s
          - \int_0^t \sigma(s, \bx+y^{(i)}_s) \, \dd w_s \\
         &- \int_0^t \int_{U_0} f(s, \bx+y^{(i)}_{s-}, u) \, \tN_\pi(\dd s, \dd u)
          - \int_0^t \int_{U_1}
             g(s, \bx+y^{(i)}_{s-}, u) \, N_\pi(\dd s, \dd u) \biggr)_{t\in\RR_+}
 \end{align*}
 on \ $(D(\RR_+, \RR^d), \cD(\RR_+, \RR^d))$ \ under \ $\PP_{1,2}$, \
 where \ $N_\pi(\dd s, \dd u)$ \ is the counting measure
 of \ $p_\pi$ \ on \ $\RR_+ \times U$, \ and
 \ $\tN_\pi(\dd s, \dd u) := N_\pi(\dd s, \dd u) - \dd s \, m(\dd u)$.
\ Using also that for each \ $i \in \{1, 2\}$, \ the first process and the
 identically 0 process are indistinguishable (since the SDE \eqref{SDE_X_YW}
 holds \ $\PP^{(i)}$-a.s. for \ $(\bX^{(i)}_t)_{t\in\RR_+}$), \ we obtain that the
 tuple
 \ $\bigl( \Omega, \cF, (\cF_t)_{t\in\RR_+}, \PP_{1,2}, w, p_\pi,
           (\bx + y^{(i)}_t)_{t\in\RR_+} \bigr)$ \ satisfies (D4), as desired.
It is worth mentioning that this is the place where we use that the filtration
 \ $(\cF_t)_{t\in\RR_+}$ \ satisfies the usual conditions in order to ensure
 that the second process above has a c\`adl\`ag modification, see
 Remark \ref{intexist}.
The filtrations \ $(\cG_t)_{t\in\RR_+}$ \ and \ $(\widetilde\cG_t)_{t\in\RR_+}$
 \ do not necessarily satisfy the usual conditions, this is the reason for
 introducing the filtration \ $(\cF_t)_{t\in\RR_+}$.

We have \ $\PP_{1,2}(\bx + y^{(1)}_0 = \bx + y^{(2)}_0) = 1$, \ because, by
 \eqref{PP},
 \ $\PP_{1,2}(y^{(i)}_0 = \bzero) = \PP^{(i)}(\bY_0^{(i)} = \bzero) = 1$,
 \ $i \in \{1, 2\}$.
\ Since
 \ $\bigl( \Omega, \cF, (\cF_t)_{t\in\RR_+}, \PP_{1,2}, w, p_\pi,
           (\bx + y^{(i)}_t)_{t\in\RR_+} \bigr)$,
 \ $i \in \{1, 2\}$, \ are weak solutions of the SDE \eqref{SDE_X_YW} with the
 same initial distribution \ $n$, \ and
 \ $\PP_{1,2}(\bx + y^{(1)}_0 = \bx + y^{(2)}_0) = 1$, \ pathwise uniqueness
 implies
 \ $\PP_{1,2}(\text{$\bx + y^{(1)}_t = \bx + y^{(2)}_t$ for all $t \in \RR_+$})
    = 1$,
 \ or equivalently,
 \begin{equation}\label{Pu}
  \PP_{1,2}[\omega = (\bx, w, \pi, y^{(1)}, y^{(2)}) \in \Omega : y^{(1)} = y^{(2)}]
  = 1 ,
 \end{equation}
 hence, applying \eqref{PP},
 \begin{align*}
  \PP^{(1)}[(\bX_0^{(1)}, \bW^{(1)}, p^{(1)}, \bY^{(1)}) \in A]
  &= \PP_{1,2}[\omega = (\bx, w, \pi, y^{(1)}, y^{(2)}) \in \Omega
         : (\bx, w, \pi, y^{(1)}) \in A] \\
  &= \PP_{1,2}[\omega = (\bx, w, \pi, y^{(1)}, y^{(2)}) \in \Omega
         : (\bx, w, \pi, y^{(2)}) \in A] \\
  &= \PP^{(2)}[(\bX_0^{(2)}, \bW^{(2)}, p^{(2)}, \bY^{(2)}) \in A]
 \end{align*}
 for all \ $A \in \cB(\Theta)$.
\ Since \ $\bX^{(i)} = \bX^{(i)}_0 + \bY^{(i)}$, \ $i \in \{1, 2\}$, \ and the
 mapping
 \ $\RR^d \times D(\RR_+, \RR^d) \ni (\bx_0, \by)
    \mapsto \bx_0+\by \in D(\RR_+, \RR^d)$
 \ is continuous (see, e.g., Jacod and Shiryaev
 \cite[Proposition VI.1.23]{JacShi}), we have
 \[
   \PP^{(1)}[\bX^{(1)} \in \tA] = \PP^{(2)}[\bX^{(2)} \in \tA] ,
   \qquad \tA \in \cD(\RR_+, \RR^d) ,
 \]
 and then we obtain uniqueness in the sense of probability law.
\proofend

\section{Precise formulation and proof of Theorem
            \ref{Thm_pathwise_2}}
\label{section_proof2}

Our first result is a counterpart of Lemma 1.1 in Chapter IV in Ikeda and
 Watanabe \cite{IkeWat} for stochastic differential equations with jumps,
 compare also with Situ \cite[page 106, Fact A]{Sit}.

\begin{Lem}\label{measurability}
If \ $\bigl( \Omega, \cF, (\cF_t)_{t\in\RR_+}, \PP, \bW, p, \bX\bigr)$ \ is a
 weak solution of the SDE \eqref{SDE_X_YW} with initial distribution \ $n$
 \ on \ $(\RR^d, \cB(\RR^d))$, \ then for every fixed \ $t \in \RR_+$ \ and
 \ $F \in \cD_t(\RR_+, \RR^d)$, \ the mapping
 \[
   \RR^d \times C(\RR_+, \RR^r) \times M(\RR_+\times U) \ni (\bx, w, \pi)
   \mapsto Q(\bx, w, \pi, F)
 \]
 is \ $\widehat{\cB}_t/\cB([0,1])$-measurable, where \ $\widehat{\cB}_t$
 \ denotes the completion of
 \ $\cB(\RR^d) \otimes \cC_t(\RR_+, \RR^r) \otimes \cM_t(\RR_+ \times U)$
 by the null sets of \ $n \times P_{\bW,r} \times P_{U,m}$ \ from
 \ $\cB(\RR^d) \otimes \cC(\RR_+, \RR^r)
    \otimes \cM(\RR_+ \times U)$.
\end{Lem}

\noindent
\textbf{Proof.}
Consider the regular conditional probability
 \[
   Q_t : \RR^d \times C(\RR_+, \RR^r) \times M(\RR_+ \times U)
            \times \cD_t(\RR_+, \RR^d)
            \to [0, 1]
 \]
 for \ $\cD_t(\RR_+, \RR^d)$ \ given \ $(x, \varphi_t(w), \psi_t(\pi))$,
 \ where, for each \ $t \in \RR_+$, \ the stopped mapping
 \ $\varphi_t : C(\RR_+, \RR^r) \to C(\RR_+, \RR^r)$ \ is defined in
 \eqref{varphi_t}, and \ $\psi_t : M(\RR_+ \times U) \to M(\RR_+ \times U)$,
 \ $\psi_t(\pi) := \pi\big|_{[0,t]\times U}$, \ $\pi \in M(\RR_+ \times U)$,
 \ i.e., \ $\psi_t(\pi)$ \ denotes the restriction of \ $\pi$ \ onto
 \ $[0, t] \times U$.
\ The mapping \ $Q_t$ \ enjoy properties analogous to (R1)--(R3).
Namely,
\begin{enumerate}
  \item[$(\widetilde{\mathrm R1})$]
   for each \ $\bx \in \RR^d$, \ $w \in C(\RR_+, \RR^r)$ \ and
    \ $\pi \in M(\RR_+ \times U)$, \ the set function
    \ $\cD_t(\RR_+, \RR^d) \ni F \mapsto Q_t(\bx, w, \pi, F)$ \ is a
    probability measure on \ $(D(\RR_+, \RR^d), \cD_t(\RR_+, \RR^d))$;
  \item[$(\widetilde{\mathrm R2})$]
   for each \ $F \in \cD_t(\RR_+, \RR^d)$, \ the mapping
   \ $\RR^d \times C(\RR_+, \RR^r) \times M(\RR_+ \times U) \ni (\bx, w, \pi)
    \mapsto Q_t(\bx, w, \pi, F)$ \ is
   \ $\cB(\RR^d) \otimes \cC_t(\RR_+, \RR^r)
      \otimes \cM_t(\RR_+ \times U)/\cB([0,1])$-measurable;
  \item[$(\widetilde{\mathrm R3})$]
   for every
    \ $G \in \cB(\RR^d) \otimes \cC_t(\RR_+, \RR^r)
       \otimes \cM_t(\RR_+ \times U)$
    \ and \ $F \in \cD_t(\RR_+, \RR^d)$,
   \begin{equation*}%\label{PGF}
    P(G \times F)
    = \int_G
       Q_t(\bx, w, \pi, F) \, n(\dd\bx) \, P_{\bW\!,\,r}(\dd w)
       \, P_{U,m}(\dd\pi) ,
   \end{equation*}
   where the probability measure \ $P$ \ is defined in \eqref{help10}.
\end{enumerate}
In order to prove the statement, it suffices to check that
 \begin{equation}\label{QiQit}
  \text{$Q(\bx, w, \pi, F) = Q_t(\bx, w, \pi, F)$ \quad
        for \ $n \times P_{\bW\!,\,r} \times P_{U,m}$-a.e. \ $(\bx, w, \pi)$.}
 \end{equation}
Indeed, then \ $(n \times P_{\bW\!,\,r} \times P_{U,m})(N) = 0$ \ for
 \begin{align*}
  N &:= \Big\{ (\bx,w,\pi) \in \RR^d \times C(\RR_+, \RR^r)
                               \times M(\RR_+ \times U)
               : Q(\bx, w, \pi, F) \ne Q_t(\bx, w, \pi, F) \Big\} \\
    &\in \cB(\RR^d) \otimes \cC(\RR_+, \RR^r)
          \otimes \cM(\RR_+ \times U) ,
 \end{align*}
 and what is more, \ $N \in \widehat{\cB}_t$, \ since
 \begin{align*}
  \widehat{\cB}_t
   = \sigma\big(\cB(\RR^d) \otimes \cC_t(\RR_+, \RR^r) \otimes \cM_t(\RR_+ \times U)
            \cup \cN\big) ,
 \end{align*}
 where
 \begin{align*}
  \cN:=\Big\{ A\subset \RR^d \times C(\RR_+, \RR^r) \times M(\RR_+ \times U)
              & : \exists \, B \in \cB(\RR^d) \otimes \cC(\RR_+, \RR^r)
                                   \otimes \cM(\RR_+ \times U) \\
              &\text{\;\; with \ $A \subset B$,
                     \ $(n \times P_{\bW,r} \times P_{U,m})(B)=0$} \Big\} ,
 \end{align*}
 and \ $N \in \cN$.
\ Hence for all \ $E \in \cB([0, 1])$,
 \begin{align*}
  \Big\{(\bx, w, \pi) \in \RR^d \times C(\RR_+, \RR^r) \times M(\RR_+ \times U)
        : Q(\bx, w, \pi, F) \in E \Big\}
  = A_1 \cup A_2 ,
 \end{align*}
 where
 \begin{align*}
  A_1
  &:=\Big\{(\bx,w,\pi)\in \RR^d \times C(\RR_+, \RR^r) \times M(\RR_+ \times U)
           : Q(\bx,w,\pi,F)\in E , \\
  &\phantom{\quad :=\Big\{ (\bx,w,\pi)\in \RR^d \times C(\RR_+, \RR^r) \times M(\RR_+ \times U)}
          Q(\bx,w,\pi,F) = Q_t(\bx,w,\pi,F) \Big\}\\
  &= \Big\{ (\bx,w,\pi)\in \RR^d \times C(\RR_+, \RR^r) \times M(\RR_+ \times U)
         : Q_t(\bx,w,\pi,F)\in E \Big\} \\
  &\quad \cap \Big\{ (\bx,w,\pi)\in \RR^d \times C(\RR_+, \RR^r) \times M(\RR_+ \times U)
         :  Q(\bx,w,\pi,F) = Q_t(\bx,w,\pi,F) \Big\},
 \end{align*}
 and
 \begin{align*}
  A_2&:= \Big\{ (\bx,w,\pi)\in \RR^d \times C(\RR_+, \RR^r) \times M(\RR_+ \times U)
         : Q(\bx,w,\pi,F)\in E,\\
     &\phantom{\quad :=\Big\{ (\bx,w,\pi)\in \RR^d \times C(\RR_+, \RR^r) \times M(\RR_+ \times U)}
          Q(\bx,w,\pi,F) \ne Q_t(\bx,w,\pi,F) \Big\} .
 \end{align*}
Here \ $A_1 \in \widehat{\cB}_t$, \ since, by \ $(\widetilde{\mathrm R2})$,
 \ the set
 \begin{align*}
  \Big\{ (\bx,w,\pi)\in \RR^d \times C(\RR_+, \RR^r) \times M(\RR_+ \times U)
          : Q_t(\bx,w,\pi,F)\in E \Big\}
 \end{align*}
 is in \ $\cB(\RR^d) \otimes \cC_t(\RR_+, \RR^r)
             \otimes \cM_t(\RR_+ \times U)
         \subset \widehat{\cB}_t$,
\ and
 \begin{align*}
   \Big\{ (\bx,w,\pi) &\in \RR^d \times C(\RR_+, \RR^r) \times M(\RR_+ \times U)
          :  Q(\bx,w,\pi,F) = Q_t(\bx,w,\pi,F) \Big\}\\
         & = \RR^d \times C(\RR_+, \RR^r) \times M(\RR_+ \times U) \setminus N
   \in \widehat{\cB}_t .
 \end{align*}
Further,
 \ $A_2\subset N \in \cB(\RR^d) \otimes \cC(\RR_+, \RR^r)
                     \otimes \cM(\RR_+ \times U)$
 \ and \ $(n \times P_{\bW,r} \times P_{U,m})(N) = 0$ \ imply
 \ $A_2 \in \cN \subset \widehat{\cB}_t$.

Unfortunately, \eqref{QiQit} does not follow from the comparison of (R3) with
 \ $(\widetilde{\mathrm R3})$, \ since still we do not know weather the
 function \ $(\bx, w, \pi) \mapsto Q(\bx, w, \pi, F)$ \ is
 \ $\cB(\RR^d) \otimes \cC_t(\RR_+, \RR^r)
    \otimes \cM_t(\RR_+\times U)/\cB([0,1])$-measurable.
In order to show \eqref{QiQit}, it suffices to check that
 $(\widetilde{\mathrm R3})$ is valid for every
 \ $G \in \cB(\RR^d) \otimes \cC(\RR_+, \RR^r) \otimes \cM(\RR_+\times U)$.
\ Indeed, then, by (R3),
 \begin{align*}
  \int_G Q(\bx, w, \pi, F) \, n(\dd\bx) \, P_{\bW,r}(\dd w) \, P_{U,m}(\dd\pi)
  = \int_G Q_t(\bx, w, \pi, F) \, n(\dd\bx) \, P_{\bW,r}(\dd w) \, P_{U,m}(\dd\pi)
 \end{align*}
 for all
 \ $G \in \cB(\RR^d) \otimes \cC(\RR_+, \RR^r) \otimes \cM(\RR_+\times U)$
 \ and \ $F\in\cD_t(\RR_+,\RR^d)$, \ and hence, using also that the function
 \ $(\bx, w, \pi) \mapsto Q_t(\bx, w, \pi, F)$ \ is
 \ $\cB(\RR^d) \otimes \cC(\RR_+, \RR^r)
    \otimes \cM(\RR_+\times U)/\cB([0,1])$-measurable,
 by the uniqueness part of the Radon-Nikod\'ym theorem, we have \eqref{QiQit}.

The class \ $\cG$ \ of sets \ $G$ \ satisfying $(\widetilde{\mathrm R3})$ is a
 Dynkin system, i.e.,
 \begin{itemize}
  \item
   $\RR^d \times C(\RR_+, \RR^r) \times M(\RR_+ \times U) \in \cG$, \ since
    \ $\RR^d \times C(\RR_+, \RR^r) \times M(\RR_+ \times U) \in \cB(\RR^d)
       \otimes \cC_t(\RR_+, \RR^r) \otimes \cM_t(\RR_+ \times U)$
    \ and one can apply $(\widetilde{\mathrm R3})$,
  \item
   if \ $G_1, G_2 \in \cG$ \ and \ $G_1 \subset G_2$, \ then
    \ $G_2 \setminus G_1 \in \cG$.
   \ Indeed,
    \begin{align*}
     P((G_2 \setminus G_1) \times F)
     &= P(G_2 \times F) - P(G_1 \times F) \\
     &= \int_{G_2}
         Q_t(\bx, w, \pi, F)
         \, n(\dd\bx) \, P_{\bW,r}(\dd w) \, P_{U,m}(\dd\pi) \\
     &\phantom{=\;}
        - \int_{G_1}
           Q_t(\bx, w, \pi, F)
           \, n(\dd\bx) \, P_{\bW,r}(\dd w) \, P_{U,m}(\dd\pi) \\
     &= \int_{G_2\setminus G_1}
         Q_t(\bx, w, \pi, F)
         \, n(\dd\bx) \, P_{\bW,r}(\dd w) \, P_{U,m}(\dd\pi) .
    \end{align*}
  \item
   if \ $(G_n)_{n\in\NN} \subset \cG$ \ and \ $G_1 \subset G_2 \subset \cdots$,
    \ then \ $\bigcup_{n=1}^\infty G_n \in \cG$.
   \ Indeed, by the continuity of probability and dominated convergence
    theorem,
    \begin{align*}
     P&\left( \left(\bigcup_{n=1}^\infty G_n\right) \times F \right)
     = \lim_{n\to\infty} P(G_n \times F) \\
     &= \lim_{n\to\infty}
         \int_{G_n}
          Q_t(\bx, w, \pi, F)
          \, n(\dd\bx) \, P_{\bW,r}(\dd w) \, P_{U,m}(\dd\pi) \\
     &= \lim_{n\to\infty}
         \int_{\RR^d\times C(\RR_+,\RR^r)\times M(\RR_+\times U)}
          Q_t(\bx, w, \pi, F) \bbone_{G_n}(\bx, w, \pi)
          \, n(\dd\bx) \, P_{\bW,r}(\dd w) \, P_{U,m}(\dd\pi) \\
     &= \int_{\bigcup_{n=1}^\infty G_n}
         Q_t(\bx, w, \pi, F)
         \, n(\dd\bx) \, P_{\bW,r}(\dd w) \, P_{U,m}(\dd\pi) .
    \end{align*}
 \end{itemize}
Consider the collection of sets of the form
 \begin{equation}\label{Ds}
  G = G_1 \times (\varphi_t^{-1}(G_2) \cap \widetilde{\varphi}_t^{-1}(G_3))
      \times (\psi_t^{-1}(G_4) \cap \widetilde{\psi}_t^{-1}(G_5))
 \end{equation}
 for \ $G_1 \in \cB(\RR^d)$, \ $G_2, G_3 \in \cC(\RR_+, \RR^r)$,
 \ $G_4, G_5 \in \cM(\RR_+\times U)$, \ where, for each \ $t \in \RR_+$,
 \ $\varphi_t$ \ and \ $\psi_t$ \ are defined earlier,
 \ $\widetilde{\varphi}_t : C(\RR_+, \RR^r) \to C(\RR_+, \RR^r)$ \ denotes the
 increment mapping \ $(\widetilde{\varphi}_t(w))(s) := w(t + s) - w(t)$,
 \ $w \in C(\RR_+, \RR^r)$, \ $s \in \RR_+$, \ and
 \ $\widetilde{\psi}_t : M(\RR_+\times U) \to M(\RR_+\times U)$ \ denotes the
 increment mapping given by
 \ $\widetilde{\psi}_t(\pi)([0,s]\times B)
    := \pi( [0,t+s]\times B) - \pi([0,t]\times B)$, $s\in\RR_+$, $B\in\cB(U)$.
\ This collection of sets is closed under pairwise intersection and generates
 the $\sigma$-algebra
 \ $\cB(\RR^d) \otimes \cC(\RR_+, \RR^r) \otimes \cM(\RR_+\times U)$,
 \ since the collection of sets of the form
 \ $(\varphi_t^{-1}(G_2) \cap \widetilde{\varphi}_t^{-1}(G_3))$ \ with
 \ $G_2 = \{ w \in C(\RR_+, \RR^r) : (w(t_1), \ldots, w(t_n)) \in A\}$ \ for
 \ $n \in \NN$, \ $t \in \RR_+$, \ $t_1, \ldots, t_n \in [0, t]$,
 \ $A \in \cB(\RR^{rn})$, \ and \ $G_3 = C(\RR_+, \RR^r)$ \ generates
 \ $\cC(\RR_+, \RR^r)$ \ by \eqref{C_infty}, and the collection of sets
 of the form \ $(\psi_t^{-1}(G_4) \cap \widetilde{\psi}_t^{-1}(G_5))$ \ with
 \[
     G_4 = \{ \pi \in M(\RR_+ \times U) : \pi([0,t]\times B)\in A\}
 \]
 for \ $t \in \RR_+$, \ $B\in\cB(U)$, \ $A\in \cB([0,\infty])$, \ and
 \ $G_5 = M(\RR_+ \times U)$ \ generates
 \ $\cM(\RR_+ \times U)$ \ by \eqref{M_infty}.
By the Dynkin system theorem (see, e.g., Karatzas and Shreve
 \cite[Theorem 2.1.3]{KarShr}),
 \ $\cB(\RR^d) \otimes \cC(\RR_+, \RR^r) \otimes \cM(\RR_+ \times U)
    \subset \cG$
 \ provided that we prove \ $(\widetilde{\mathrm R3})$ \ for \ $G$ \ of the
 form \eqref{Ds}.
For such a \ $G$, \ by Fubini theorem, we have
 \begin{align*}
  &\int_G
    Q_t(\bx, w, \pi, F) \, n(\dd\bx) \, P_{\bW,r}(\dd w) \, P_{U,m}(\dd\pi) \\
  &= \int_{\psi_t^{-1}(G_4)\cap\widetilde{\psi}_t^{-1}(G_5)}
      \left( \int_{\varphi_t^{-1}(G_2)\cap\widetilde{\varphi}_t^{-1}(G_3)}
      \left( \int_{G_1}  Q_t(\bx, w, \pi, F) \, n(\dd\bx) \right)
      P_{\bW,r}(\dd w) \right) P_{U,m}(\dd\pi) \\
  &= \EE_{P_{\bW,r}\times P_{U,m}}
     \left[ \int_{G_1}  Q_t(\bx, w, \pi, F) \, n(\dd\bx)
            \bbone_{\varphi_t^{-1}(G_2)\cap\widetilde{\varphi}_t^{-1}(G_3)}(w)
            \bbone_{\psi_t^{-1}(G_4)\cap\widetilde{\psi}_t^{-1}(G_5)}(\pi) \right] \\
  &= \EE_{P_{\bW,r}\times P_{U,m}}
     \Big[ \EE_{P_{\bW,r}\times P_{U,m}} \Big[ \int_{G_1}  Q_t(\bx, w, \pi, F)
           \, n(\dd\bx) \\
  &\phantom{ = \EE_{P_{\bW,r}\times P_{U,m}} \Big[}
      \times \bbone_{\varphi_t^{-1}(G_2)}(w)
             \bbone_{\widetilde{\varphi}_t^{-1}(G_3)}(w)
             \bbone_{\psi_t^{-1}(G_4)}(\pi) \bbone_{\widetilde{\psi}_t^{-1}(G_5)}(\pi)
             \;\Big|\; \cC_t(\RR_+,\RR^r) \otimes \cM_t(\RR_+\times U)\Big] \Big] \\
  &= \EE_{P_{\bW,r}\times P_{U,m}}
       \biggl[ \int_{G_1}
                Q_t(\bx, w, \pi, F) \, n(\dd\bx)
                \bbone_{\varphi_t^{-1}(G_2)}(w) \bbone_{\psi_t^{-1}(G_4)}(\pi) \\
   &\phantom{= \EE_{P_{\bW,r}\times P_{U,m}} \biggl[ \int_{G_1}}
                (P_{\bW,r}\times P_{U,m})
                 \left( \widetilde{\varphi}_t^{-1}(G_3)
                        \times \widetilde{\psi}_t^{-1}(G_5)
                        \;\Big|\; \cC_t(\RR_+, \RR^r) \otimes \cM_t(\RR_+\times U) \right)
       \biggr] \\
   &= \EE_{P_{\bW,r}\times P_{U,m}}
       \biggl[ \int_{G_1}
               Q_t(\bx, w, \pi, F) \, n(\dd\bx)
               \bbone_{\varphi_t^{-1}(G_2)}(w)
               \bbone_{\psi_t^{-1}(G_4)}(\pi)  \\
   &\phantom{= \EE_{P_{\bW,r}\times P_{U,m}} \biggl[ \int_{G_1}}
      \times(P_{\bW,r}\times P_{U,m})
       \left( \widetilde{\varphi}_t^{-1}(G_3)
              \times \widetilde{\psi}_t^{-1}(G_5) \right)\biggr]\\
   &= \int_{G_1 \times \varphi_t^{-1}(G_2) \times \psi_t^{-1}(G_4)}
       Q_t(\bx, w, \pi, F) \, n(\dd\bx) \, P_{\bW,r}(\dd w) \, P_{U,m}(\dd\pi) \\
   &\quad
      \times(P_{\bW,r}\times P_{U,m})
      \left( \widetilde{\varphi}_t^{-1}(G_3)
              \times \widetilde{\psi}_t^{-1}(G_5) \right) \\
   &= P[G_1 \times \varphi_t^{-1}(G_2) \times \psi_t^{-1}(G_4) \times F] \,
      (P_{\bW,r}\times P_{U,m})
       \left( \widetilde{\varphi}_t^{-1}(G_3)
              \times \widetilde{\psi}_t^{-1}(G_5) \right) .
 \end{align*}
The fourth equality above follows from the
 \ $\cC_t(\RR_+,\RR^r) \otimes \cM_t(\RR_+\times U)/\cB([0,1])$-measurability
 of the function
 \[
   C(\RR_+,\RR^r)\times M(\RR_+\times U)\ni (w,\pi)
   \mapsto \int_{G_1}  Q_t(\bx, w, \pi, F) \, n(\dd\bx) ,
 \]
 which is a consequence of $(\widetilde{\mathrm R2})$ and Fubini theorem.
The fifth equality above follows from the independence of
 \ $\widetilde{\varphi}_t^{-1}(G_3) \times \widetilde{\psi}_t^{-1}(G_5)$ \ and
 \ $\cC_t(\RR_+, \RR^r) \otimes \cM_t(\RR_+ \times U)$ \ under the measure
 \ $P_{\bW,r}\times P_{U,m}$, \ see, e.g., Ikeda and Watanabe
 \cite[Chapter 2, Theorems 6.4 and 6.5]{IkeWat}.
For the last equality above we used \ $(\widetilde{\mathrm R3})$ \ and
 \[
   G_1 \times \varphi_t^{-1}(G_2) \times \psi_t^{-1}(G_4) \times F
   \in \cB(\RR^d) \otimes \cC_t(\RR_+, \RR^r) \otimes \cM_t(\RR_+ \times U)
       \otimes \cD_t(\RR_+, \RR^d) .
  \]
By \eqref{help10},
 \begin{align*}
  (P_{\bW,r} \times P_{U,m})
  \left( \widetilde{\varphi}_t^{-1}(G_3)
         \times \widetilde{\psi}_t^{-1}(G_5) \right)
  &= P[(\bx, w, \pi, y) \in \Theta
         : \widetilde{\varphi}_t(w) \in G_3 , \,
           \widetilde{\psi}_t(\pi) \in G_5 ] \\
  &= \PP[\widetilde{\varphi}_t(\bW) \in G_3 , \,
              \widetilde{\psi}_t(p) \in G_5 ] ,
 \end{align*}
 \[
   P[G_1 \times \varphi_t^{-1}(G_2) \times \psi_t^{-1}(G_4) \times F]
   = \PP[\bX_0 \in G_1, \, \varphi_t(\bW) \in G_2 , \,
              \psi_t(p) \in G_4 , \, \bY \in F] .
 \]
Therefore, if \ $G$ \ is of the form \eqref{Ds}, then
 \begin{align*}
  &\int_G
    Q_t(\bx, w, \pi, F) \, n(\dd\bx) \, P_{\bW,r}(\dd w) \, P_{U,m}(\dd\pi) \\
  &= \PP[\bX_0 \in G_1 , \, \varphi_t(\bW) \in G_2 , \,
              \psi_t(p) \in G_4 , \, \bY \in F] \,
     \PP[\widetilde{\varphi}_t(\bW) \in G_3 , \,
               \widetilde{\psi}_t(p) \in G_5 ] \\
  &= \PP[\bX_0 \in G_1 , \, \varphi_t(\bW) \in G_2 , \,
              \widetilde{\varphi}_t(\bW) \in G_3 , \,
              \psi_t(p) \in G_4 , \,
              \widetilde{\psi}_t(p) \in G_5 , \, \bY \in F] \\
  &= \PP[(\bX_0, \bW, p) \in G , \, \bY \in F] \\
  &= P[G \times F] .
 \end{align*}
The second equality above follows from the independence of
 \ $\{\bX_0 \in G_1 , \, \varphi_t(\bW) \in G_2 , \,
      \psi_t(p) \in G_4 , \, \bY \in F\}$
 \ and
 \ $\{\widetilde{\varphi}_t(\bW) \in G_3 , \,
      \widetilde{\psi}_t(p) \in G_5\}$
 \ under the probability measure \ $\PP$.
\ This independence holds because
 \begin{align}\label{meas}
  \begin{aligned}
   \{\bX_0 \in G_1 , & \, \varphi_t(\bW) \in G_2 , \,
      \psi_t(p) \in G_4 , \, \bY \in F\}\\
   & = \{\bX_0 \in G_1 , \, \varphi_t(\varphi_t(\bW)) \in G_2 , \,
        \psi_t(\psi_t(p)) \in G_4 , \, \bY \in F\}\\
   & = \{\bX_0 \in G_1 , \, \varphi_t(\bW) \in \varphi_t^{-1}(G_2), \,
         \psi_t(p) \in \psi_t^{-1}(G_4) , \, \bY \in F\}
     \in \cF_t
  \end{aligned}
 \end{align}
 and
 \ $\{\widetilde{\varphi}_t(\bW) \in G_3 , \,
      \widetilde{\psi}_t(p) \in G_5\}$ \ is independent of \ $\cF_t$
 \ under the probability measure \ $\PP$, \ see, e.g., Ikeda and Watanabe
 \cite[Chapter II, Theorems 6.4 and 6.5]{IkeWat}.
The relationship \eqref{meas} is valid since
 \ $\varphi_t^{-1}(G_2)\in\cC_t(\RR_+,\RR^r)$,
 \ $\psi_t^{-1}(G_4)\in\cM_t(\RR_+\times U)$
 \ and \ $F\in \cD_t(\RR_+,\RR^d)$, \ the mapping
 \ $\Omega \ni \omega \mapsto \varphi_t(\bW(\omega))$ \ is
 \ $\cF_t / \cC_t(\RR_+,\RR^r)$-measurable, and the mapping
 \ $\Omega \ni \omega \mapsto \psi_t(p(\omega))$ \ is
 \ $\cF_t / \cM_t(\RR_+ \times U)$-measurable, because the processes
 \ $\bW$ \ and \ $p$ \ are $(\cF_t)_{t\in\RR_+}$-adapted.
\proofend

\begin{Rem}
The filtration \ $(\widehat{\cB}_t)_{t\in\RR_+}$ \ defined in Lemma
 \ref{measurability} is the augmentated filtration generated by the
 coordinate processes on the canonical probability space
 \ $(\RR^d \times C(\RR_+, \RR^r) \times M(\RR_+\times U),
     \cB(\RR^d) \otimes \cC(\RR_+, \RR^r) \otimes \cM(\RR_+\times U) ,
     n \times P_{\bW,r} \times P_{U,m})$.
\ This is the counterpart of the augmentated filtration
 \ $(\cF_t^{\bxi,\bW,p})_{t\in\RR_+}$.
\proofend
\end{Rem}

The next lemma is a generalization of Corollary 1 in Yamada and Watanabe
 \cite{YamWat} (see also Problem 5.3.22 in Karatzas and Shreve \cite{KarShr})
 for stochastic differential equations with jumps.

\begin{Lem}\label{k}
Let us suppose that pathwise uniqueness holds for the SDE \eqref{SDE_X_YW}.
If
 \ $\bigl( \Omega^{(i)}, \cF^{(i)}, (\cF_t^{(i)})_{t\in\RR_+}, \PP^{(i)}, \bW^{(i)},
           p^{(i)}, \bX^{(i)} \bigr)$, $i \in \{1, 2\}$,
 \ are two weak solutions of the SDE \eqref{SDE_X_YW} with the same initial
 distribution \ $n$ \ on \ $(\RR^d, \cB(\RR^d))$, \ then there exists a
 function
 \ $k : \RR^d \times C(\RR_+, \RR^r) \times M(\RR_+\times U)
        \to D(\RR_+, \RR^d)$
 \ such that
 \begin{equation}\label{Q1Q2}
  Q_i(\bx, w, \pi, \{ k(\bx, w, \pi) \} ) = 1, \quad i\in\{1,2\}
 \end{equation}
 holds for \ $n \times P_{\bW,r} \times P_{U,m}$-almost every
 \ $(\bx, w, \pi) \in \RR^d \times C(\RR_+, \RR^r) \times M(\RR_+ \times U)$,
 \ where \ $Q_i$, $i\in\{1,2\}$, \ is given in \eqref{help16}.
This function \ $k$ \ is
 \ $\cB(\RR^d) \otimes \cC(\RR_+, \RR^r)\otimes \cM(\RR_+ \times U)
    /\cD(\RR_+, \RR^d)$-measurable,
 \ $\widehat{\cB}_t / \cD_t(\RR_+, \RR^d)$-measurable for every fixed
 \ $t \in \RR_+$, \ and
 \begin{align}\label{help13}
  \PP^{(i)}( k(\bX_0^{(i)}, \bW^{(i)}, p^{(i)}) = \bY^{(i)}) = 1, \quad i\in\{1,2\}.
 \end{align}
\end{Lem}

\noindent
\textbf{Proof.}
Fix \ $(\bx, w, \pi) \in \RR^d \times C(\RR_+, \RR^r) \times M(\RR_+ \times U)$
 \ and define the measure
 \ $Q_{1,2}(\bx, w, \pi, \dd y^{(1)}, \dd y^{(2)})
    := Q_1(\bx, w, \pi, \dd y^{(1)}) Q_2(\bx, w, \pi, \dd y^{(2)})$
 \ on the space \ $S := D(\RR_+, \RR^d) \times D(\RR_+, \RR^d)$ \ equipped with
 the $\sigma$-algebra
 \ $\cS := \cD(\RR_+, \RR^d) \otimes \cD(\RR_+, \RR^d)$.
\ By \eqref{help12} and Fubini theorem,
 \begin{equation}\label{P(GB)}
  \PP_{1,2}[G \times B]
  = \int_G
     Q_{1,2}(\bx, w, \pi, B) \, n(\dd\bx) \, P_{\bW,r}(\dd w) \, P_{U,m}(\dd\pi)
 \end{equation}
 for all
 \ $G \in \cB(\RR^d) \otimes \cC(\RR_+, \RR^r)
          \otimes \cM(\RR_+ \times U)$
 \ and \ $B \in \cS$.
\ With the choice
 \ $G = \RR^d \times C(\RR_+, \RR^r) \times M(\RR_+ \times U)$ \ and
 \ $B = \{ (y^{(1)}, y^{(2)}) \in S : y^{(1)} = y^{(2)} \}$,
 \ using that pathwise uniqueness holds for the SDE \eqref{SDE_X_YW}, relation
 \eqref{Pu} yields \ $\PP_{1,2}[G \times B] = 1$.
\ Since \ $Q_{1,2}(\bx,w,\pi,B)\leq 1$ \ for all
 \ $(\bx, w, \pi) \in \RR^d \times C(\RR_+,\RR^r) \times M(\RR_+ \times U)$,
 \ \eqref{P(GB)} yields the existence of a set
 \ $N \in \cB(\RR^d) \otimes \cC(\RR_+, \RR^r)
          \otimes \cM(\RR_+ \times U)$
 \ with \ $(n \times P_{\bW,r} \times P_{U,m})(N) = 0$ \ such that
 \[
   Q_{1,2}\big(\bx, w, \pi, \{ (y^{(1)}, y^{(2)}) \in S : y^{(1)} = y^{(2)} \}\big) = 1,
   \qquad (\bx, w, \pi) \notin N.
 \]
Again, by Fubini theorem,
 \begin{align}\label{help14}
  \begin{split}
   1 &= Q_{1,2}(\bx, w, \pi, \{ (y^{(1)}, y^{(2)}) \in S : y^{(1)} = y^{(2)} \} )\\
     &= \int_{D(\RR_+, \RR^d)} Q_1(\bx, w, \pi, \{y\}) \, Q_2(\bx, w, \pi, \dd y),
     \qquad (\bx, w, \pi) \notin N,
  \end{split}
 \end{align}
 which can occur only if for some \ $y_0 \in D(\RR_+,\RR^d)$, \ call it
 \ $\tk(\bx, w, \pi)$, \ we have
 \begin{align}\label{help15}
  Q_i(\bx, w, \pi, \{\tk(\bx, w, \pi)\}) = 1, \qquad i \in \{1, 2\},
  \quad (\bx, w, \pi) \notin N.
 \end{align}
Indeed, since for all
 \ $(\bx, w, \pi, y) \in \RR^d \times C(\RR_+, \RR^r) \times M(\RR_+ \times U)
                         \times D(\RR_+, \RR^d)$,
 \ $Q_1(\bx, w, \pi, \{y\})\in[0,1]$, \ we have
 \[
   Q_2(\bx, w, \pi, \{ y\in D(\RR_+,\RR^d) : Q_1(\bx,w,\pi,\{y\}) = 1\})  = 1 ,
   \qquad (\bx, w, \pi) \notin N .
 \]
Since for all
 \ $(\bx, w, \pi) \in \RR^d \times C(\RR_+, \RR^r) \times M(\RR_+ \times U)$,
 \ by (R1), the set function
 \ $\cD(\RR_+, \RR^d) \ni F \mapsto  Q_i(\bx, w, \pi, F)$ \ is a probability
 measure on \ $(D(\RR_+, \RR^d), \cD(\RR_+, \RR^d))$, \ $i \in \{1, 2\}$, \ we
 get the unique existence of \ $\tk(\bx, w, \pi)$ \ for all
 \ $(\bx, w, \pi) \notin N$ \ satisfying \eqref{help15}.
Then we have \eqref{Q1Q2} for \ $\widetilde k$.

For \ $(\bx, w, \pi) \notin N$ \ and any \ $B \in \cD(\RR_+, \RR^d)$,
 \ we have \ $\tk(\bx, w, \pi) \in B$ \ if and only if
 \ $Q_i(\bx, w, \pi, B) = 1$, \ $i \in \{1, 2\}$.

The aim of the following discussion is to show the
 \ $\widehat{\cB}_t / \cD_t(\RR_+, \RR^d)$-measurability of \ $\tk$ \ for all
 \ $t \in \RR_+$.
\ For all \ $t \in \RR_+$ \ and \ $B \in \cD_t(\RR_+, \RR^d)$, \ we have
 \begin{align*}
  \tk^{-1}(B)
  &= \{ (\bx, w, \pi) \in \RR^d \times C(\RR_+, \RR^r) \times M(\RR_+ \times U)
        : \tk(\bx, w, \pi) \in B \}
   =: A_1 \cup A_2 ,
 \end{align*}
 where
 \begin{align}\label{help_A1}
   A_1 := \{ (\bx, w, \pi)
             \in \RR^d \times C(\RR_+, \RR^r) \times M(\RR_+ \times U)
             : \tk(\bx, w, \pi) \in B , \, (\bx, w, \pi) \in N \}
 \end{align}
 and
 \begin{align}\label{help_A2}
 \begin{split}
  &A_2 := \{ (\bx, w, \pi)
             \in \RR^d \times C(\RR_+, \RR^r) \times M(\RR_+ \times U)
             : \tk(\bx, w, \pi) \in B , \, (\bx, w, \pi) \notin N \} \\
      &= \{ (\bx, w, \pi)
             \in \RR^d \times C(\RR_+, \RR^r) \times M(\RR_+ \times U)
             : (\bx, w, \pi) \notin N \}
         \cap Q_i(\cdot, \cdot, \cdot, B)^{-1}(\{1\})
 \end{split}
 \end{align}
 for \ $i\in\{1,2\}$.
\ Lemma \ref{measurability} implies
 \ $Q_i(\cdot, \cdot, \cdot, B)^{-1}(\{1\}) \in \widehat{\cB}_t$,
 \ $i \in \{1, 2\}$.
\ Moreover, \ $N \in \widehat{\cB}_t$ \ (due to the definition of
 \ $\widehat{\cB}_t$, \ for more details, see the proof of Lemma
 \ref{measurability}), hence \ $A_2 \in \widehat{\cB}_t$.
\ Using that \ $A_1 \subset N$, \ $(n \times P_{\bW,r} \times P_{U,m})(N) = 0$
 \ and the definition of the augmented $\sigma$-algebra \ $\widehat{\cB}_t$
 \ (see Lemma \ref{measurability}), we obtain \ $A_1 \in \widehat{\cB}_t$.
\ Hence \ $\tk^{-1}(B) = A_1 \cup A_2 \in \widehat\cB_t$, \ as desired.

The aim of the following discussion is to show that \ $\tk$ \ is
 \[
   \overline{\cB(\RR^d)\otimes\cC(\RR_+,\RR^r)\otimes\cM(\RR_+\times U)}
   ^{\,n \times P_{\bW,r} \times P_{U,m}}/\cD(\RR_+,\RR^d)\text{-measurable} ,
 \]
 where
 \ $\overline{\cB(\RR^d)\otimes\cC(\RR_+,\RR^r)\otimes\cM(\RR_+\times U)}
    ^{\,n \times P_{\bW,r} \times P_{U,m}}$
 \ denotes the completion of
 \ $\cB(\RR^d)\otimes\cC(\RR_+,\RR^r)\otimes\cM(\RR_+\times U)$ \ with
 respect to the measure \ $n \times P_{\bW,r} \times P_{U,m}$.
\ For all \ $B \in \cD(\RR_+, \RR^d)$, \ we have
 \ $\tk^{-1}(B) = A_1 \cup A_2$, \ where \ $A_1$ \ and \ $A_2$ \ are defined in
 \eqref{help_A1} and \eqref{help_A2}.
Property (R2) implies
 \ $Q_i(\cdot, \cdot, \cdot, B)^{-1}(\{1\})
    \in \cB(\RR^d)\otimes\cC(\RR_+,\RR^r)\otimes\cM(\RR_+\times U)$,
 \ $i \in \{1, 2\}$.
\ Moreover, by definition of completion (see, e.g., Definition 2.7.2 in
 Karatzas and Shreve \cite{KarShr}),
 \[
   N
   \in\overline{\cB(\RR^d)\otimes\cC(\RR_+,\RR^r)\otimes\cM(\RR_+\times U)}
      ^{\,n \times P_{\bW,r} \times P_{U,m}} ,
 \]
 hence
 \[
   A_2
   \in\overline{\cB(\RR^d)\otimes\cC(\RR_+,\RR^r)\otimes\cM(\RR_+\times U)}
      ^{\,n \times P_{\bW,r} \times P_{U,m}} .
 \]
Using that \ $A_1 \subset N$, \ $(n \times P_{\bW,r} \times P_{U,m})(N) = 0$,
 \ by definition of completion, we obtain
 \[
   A_1
   \in\overline{\cB(\RR^d)\otimes\cC(\RR_+,\RR^r)\otimes\cM(\RR_+\times U)}
      ^{\,n \times P_{\bW,r} \times P_{U,m}} .
 \]
Hence
 \ $\tk^{-1}(B)
    = A_1 \cup A_2 \in
      \overline{\cB(\RR^d)\otimes\cC(\RR_+,\RR^r)\otimes\cM(\RR_+\times U)}
      ^{\,n \times P_{\bW,r} \times P_{U,m}}$, \ as desired.

Next we check  \eqref{help13} for \ $\tk$.
\ For \ $i \in \{1, 2\}$, \ by \eqref{PP}, \eqref{help12}, (R1) and
 \eqref{help15},
 \begin{align*}
  \PP^{(i)}(\tk(\bX_0^{(i)},&\bW^{(i)},p^{(i)}) = \bY^{(i)})
  = \PP_{1,2}\big(\omega = (\bx, w, \pi, y^{(1)}, y^{(2)}) \in \Omega :
                   \tk(\bx, w, \pi) = y^{(i)} \big) \\
  &= \int_{\RR^d \times C(\RR_+, \RR^r) \times M(\RR_+ \times U)}
         Q_i(\bx, w, \pi, \{ \tk(\bx,w,\pi)\})
                     \, n(\dd\bx) \, P_{\bW\!,\,r}(\dd w) \, P_{U,m}(\dd\pi)
   = 1 ,
 \end{align*}
 as desired.

It remains to check that one can choose a version of \ $\tk$ \ which
 is
 \ $\cB(\RR^d) \otimes \cC(\RR_+, \RR^r)\otimes \cM(\RR_+ \times U)
    /\cD(\RR_+, \RR^d)$-measurable,
 \ $\widehat{\cB}_t / \cD_t(\RR_+, \RR^d)$-measurable for every fixed
 \ $t \in \RR_+$, \ and \eqref{Q1Q2} and \eqref{help13} remain hold for \ $k$.
\ Since \ $\tk$ \ is
 \[
   \overline{\cB(\RR^d)\otimes\cC(\RR_+,\RR^r)\otimes\cM(\RR_+\times U)}
   ^{\,n \times P_{\bW,r} \times P_{U,m}}/\cD(\RR_+,\RR^d)\text{-measurable} ,
 \]
 there exists a function
 \ $k : \RR^d \times C(\RR_+, \RR^r) \times M(\RR_+\times U)
        \to D(\RR_+, \RR^d)$
 \ which is
 \ $\cB(\RR^d) \otimes \cC(\RR_+,\RR^r) \otimes \cM(\RR_+\times U)
    / \cD(\RR_+,\RR^d)$-measurable
 and
 \[
   (n \times P_{\bW,r} \times P_{U,m})
   \Big(\{ (\bx,w,\pi)\in \RR^d\times C(\RR_+,\RR^r)\times M(\RR_+\times U)
           : \tk(\bx,w,\pi) \ne k(\bx,w,\pi) \}\Big) = 0
 \]
 see, e.g., Cohn \cite[Proposition 2.2.5]{Coh}.
First we check that \ $k$ \ is
 \ $\widehat{\cB}_t / \cD_t(\RR_+, \RR^d)$-measurable for every fixed
 \ $t \in \RR_+$.
\ For all \ $t \in \RR_+$ \ and \ $B \in \cD_t(\RR_+, \RR^d)$, \ we have
 \begin{align*}
  k^{-1}(B) & = (k^{-1}(B) \cap \{\tk=k\})
                         \cup (k^{-1}(B) \cap \{\tk\ne k\})\\
           & = (\tk^{-1}(B) \cap \{\tk=k\})
                           \cup (k^{-1}(B) \cap \{\tk\ne k\}),
 \end{align*}
 where \ $\tk^{-1}(B) \in \widehat{\cB}_t$ \ (since \ $\tk$ \ is
 \ $\widehat{\cB}_t / \cD_t(\RR_+, \RR^d)$-measurable),
 \ $ \{\tk \ne k\} \in \widehat{\cB}_t$ \ (due to the definition of
 completion, since \ $(n \times P_{\bW,r} \times P_{U,m})(\tk\ne k)=0$),
 \ $\{\tk = k\} \in \widehat{\cB}_t$ \ (since \ $ \widehat{\cB}_t$ \ is a
 \ $\sigma$-algebra), and
 \ $k^{-1}(B) \cap \{\tk\ne k\} \in \widehat{\cB}_t$ \ (due to the definition
 of completion, since \ $k^{-1}(B) \cap \{\tk\ne k\} \subset \{\tk\ne k\}$).
\ Hence \ $k^{-1}(B) \in \widehat{\cB}_t$.

Next we check \eqref{Q1Q2} for \ $k$.
\ Using that \eqref{Q1Q2} holds for \ $\tk$ \ and
 \ $(n \times P_{\bW,r} \times P_{U,m})(\tk\ne k)=0$, \ we have
 \begin{align*}
  &(n \times P_{\bW,r} \times P_{U,m})( H_1 \cup H_2 ) \\
  &= (n \times P_{\bW,r} \times P_{U,m})
     \big( (H_1 \cup H_2) \cap \{k = \tk\} \big)
     + (n \times P_{\bW,r} \times P_{U,m})
       \big( (H_1 \cup H_2) \cap \{k \ne \tk\} \big) \\
  &\leq (n \times P_{\bW,r} \times P_{U,m})( \widetilde H_1 \cup \widetilde H_2 )
        + (n \times P_{\bW,r} \times P_{U,m})(k \ne \tk)
   = 0 + 0 = 0 ,
 \end{align*}
 where
 \begin{align*}
  &\widetilde H_i
   := \{ (\bx,w,\pi)\in \RR^d\times C(\RR_+,\RR^r)\times M(\RR_+\times U)
         : Q_i(\bx,w,\pi, \{\widetilde k(\bx,w,\pi)\})\ne 1\}, \\
  &H_i := \{ (\bx,w,\pi)\in \RR^d\times C(\RR_+,\RR^r)\times M(\RR_+\times U)
             : Q_i(\bx,w,\pi, \{k(\bx,w,\pi)\})\ne 1\}
 \end{align*}
 for \ $i\in\{1, 2\}$.
\ This implies \eqref{Q1Q2} for \ $k$.

Finally, we check \eqref{help13} for \ $k$.
\ First observe that \ $\PP_{1,2}(\tk = k) = 1$, \ since, by \eqref{P(GB)},
 \begin{align*}
  &\PP_{1,2}(\tk = k) = 1 - \PP_{1,2}(\tk \ne k) \\
  & = 1 - \int_{\{\tk \ne k\}}
           Q_{1,2}(\bx,w,\pi, D(\RR_+,\RR^d), D(\RR_+,\RR^d)) \,
                     n(\dd\bx) \, P_{\bW\!,\,r}(\dd w) \, P_{U,m}(\dd\pi) \\
  & = 1 - \int_{\{\tk \ne k\}}
           Q_1(\bx,w,\pi, D(\RR_+,\RR^d)) Q_2(\bx,w,\pi, D(\RR_+,\RR^d)) \,
                     n(\dd\bx) \, P_{\bW\!,\,r}(\dd w) \, P_{U,m}(\dd\pi) \\
  & = 1 - (n\times P_{\bW\!,\,r}\times P_{U,m})(\tk \ne k)
    = 1 - 0 = 1 ,
 \end{align*}
 where we used (R1) as well.
Then, by \eqref{PP} and \eqref{help12}, for \ $i \in \{1, 2\}$, \ we obtain
 \begin{align*}
  \PP^{(i)}( k(\bX_0^{(i)},& \bW^{(i)}, p^{(i)}) = \bY^{(i)})
  = \PP_{1,2}\big(\omega = (\bx, w, \pi, y^{(1)}, y^{(2)}) \in \Omega :
                   k(\bx, w, \pi) = y^{(i)} \big) \\
    & =  \PP_{1,2}\big(\{\omega = (\bx, w, \pi, y^{(1)}, y^{(2)}) \in \Omega :
                   k(\bx, w, \pi) = y^{(i)} \} \cap \{ \tk = k\} \big) \\
    & =  \PP_{1,2}\big(\{\omega = (\bx, w, \pi, y^{(1)}, y^{(2)}) \in \Omega :
                    \tk(\bx, w, \pi) = y^{(i)} \} \cap \{ \tk = k\} \big) \\
    & =  \PP_{1,2}\big(\omega = (\bx, w, \pi, y^{(1)}, y^{(2)}) \in \Omega :
                        \tk(\bx, w, \pi) = y^{(i)}  \big) \\
    & = \PP^{(i)}( \tk(\bX_0^{(i)}, \bW^{(i)}, p^{(i)}) = \bY^{(i)})
      = 1 ,
 \end{align*}
 where, for the last equality, we applied that \eqref{help13} holds for
 \ $\tk$.
\proofend

\begin{Rem}
Note that the function \ $k$ \ in Lemma \ref{k} and the
 \ $n \times P_{\bW,r} \times P_{U,m}$-null set on which \eqref{Q1Q2} does not
 hold depend on the two weak solutions in question.
\proofend
\end{Rem}

Applying Lemma \ref{k} for weak solutions
 \ $\bigl( \Omega^{(i)}, \cF^{(i)}, (\cF_t^{(i)})_{t\in\RR_+}, \PP^{(i)}, \bW^{(i)},
           p^{(i)}, \bX^{(i)} \bigr)
    = \bigl( \Omega, \cF, (\cF_t)_{t\in\RR_+}, \PP, \bW, p, \bX \bigr)$,
 \ $i \in \{1, 2\}$, \ of the SDE \eqref{SDE_X_YW} with the same initial
 distribution \ $n$ \ on \ $(\RR^d, \cB(\RR^d))$, \ we obtain the following
 corollary.

\begin{Cor}\label{Corrolary1}
If pathwise uniqueness holds for the SDE \eqref{SDE_X_YW} and
 $\bigl( \Omega, \cF, (\cF_t)_{t\in\RR_+}, \PP, \bW, p, \bX \bigr)$ is a weak
 solution of the SDE \eqref{SDE_X_YW} with initial distribution \ $n$ \ on
 \ $(\RR^d, \cB(\RR^d))$, \ then there exists a function
 \ $k : \RR^d \times C(\RR_+, \RR^r) \times M(\RR_+\times U)
        \to D(\RR_+, \RR^d)$ \  such that
 \ $Q(\bx, w, \pi, \{ k(\bx, w, \pi) \} ) = 1$ \ holds for
 \ $n \times P_{\bW,r} \times P_{U,m}$-almost every
 \ $(\bx, w, \pi) \in \RR^d \times C(\RR_+, \RR^r) \times M(\RR_+ \times U)$,
 \ where \ $Q$ \ is given in \eqref{help11}.
This function \ $k$ \ is
 \ $\cB(\RR^d) \otimes \cC(\RR_+, \RR^r)\otimes \cM(\RR_+ \times U)
    /\cD(\RR_+, \RR^d)$-measurable,
 \ $\widehat{\cB}_t / \cD_t(\RR_+, \RR^d)$-measurable for every fixed
 \ $t \in \RR_+$, \ and \ $\PP( k(\bX_0, \bW, p) = \bY) = 1$.
\end{Cor}

Next we give the precise formulation of Theorem \ref{Thm_pathwise_2}.

\begin{Thm}\label{YW}
Let us suppose that pathwise uniqueness holds for the SDE \eqref{SDE_X_YW} and
 there exists a weak solution
 \ $\bigl( \Omega', \cF', (\cF'_t)_{t\in\RR_+}, \PP', \bW', p', \bX' \bigr)$
 \ of the SDE \eqref{SDE_X_YW} with initial distribution \ $n'$.
\ Then there exists a function
 \ $h' : \RR^d \times C(\RR_+, \RR^r) \times M(\RR_+\times U)
         \to D(\RR_+, \RR^d)$
 \ which is
 \ $\cB(\RR^d) \otimes \cC(\RR_+, \RR^r)\otimes \cM(\RR_+ \times U)
    /\cD(\RR_+, \RR^d)$-measurable,
 \ $\hcB_t / \cD_t(\RR_+, \RR^d)$-measurable for every fixed \ $t \in \RR_+$,
 \ and
 \begin{equation}\label{th}
  \bX' = h'(\bX'_0, \bW', p')  \qquad \text{$\PP'$-almost surely.}
 \end{equation}
Moreover, if objects \textup{(E1)--(E4)} are given such that the distribution
 of \ $\bxi$ \ is \ $n'$, \ then the process
 \[
   \bX := h'(\bxi, \bW, p)
 \]
 is a strong solution of the SDE \eqref{SDE_X_YW} with initial value \ $\bxi$.
\end{Thm}

\noindent
\textbf{Proof.}
Let \ $h'(\bx, w, \pi) := \bx + k'(\bx, w, \pi)$ \ for \ $\bx \in \RR^d$,
 \ $w \in C(\RR_+, \RR^r)$, \ $\pi \in M(\RR_+ \times U)$, \ where \ $k'$ \ is
 as in Corollary \ref{Corrolary1}.
By Corollary \ref{Corrolary1}, for the function \ $h'$, \ the desired
 measurability properties hold.
Using Corollary \ref{Corrolary1} and \ $\bX' = \bX'_0 + \bY'$, \ we have
 \begin{align*}
  \PP'\left( \bX' = h'(\bX'_0, \bW', p') \right)
  & = \PP'\left( \bX'_0 + \bY' = \bX'_0 + k'(\bX'_0, \bW', p') \right) \\
  & = \PP'\left( \bY' = k'(\bX'_0, \bW', p') \right)
    = 1,
 \end{align*}
 implying \eqref{th}.

Note that, for \ $\bxi$, \ $\bW$ \ and \ $p$ \ as described in (E1)--(E4), the
 triplets \ $(\bX'_0, \bW', p')$ \ and \ $(\bxi, \bW, p)$ \ induce the same
 probability measure \ $n' \times P_{\bW,r} \times P_{U,m}$ \ on the measurable
 space
 \[
   \Big( \RR^d \times C(\RR_+,\RR^r) \times M(\RR_+\times U) ,
         \cB(\RR^d) \otimes \cC(\RR_+, \RR^r)
         \otimes \cM(\RR_+\times U) \Big)
 \]
 with respect to the probability measure \ $\PP'$ \ and \ $\PP$,
 \ respectively, where \ $\PP$ \ denotes the probability measure appears in
 (E1), since \ $\bX'_0$, $\bW'$, $p'$ \ are $\PP'$-independent and
 \ $\bxi$, $\bW$, $p$ \ are $\PP$-independent, see Remarks \ref{dRM_YW} and
 \ref{dRM_strong}.

Observe also that the mappings
 \begin{align}\label{mapping_1}
  \Omega' \ni \omega'
  \mapsto (\bX'_0(\omega'), (\bW_t'(\omega'))_{t\in\RR_+}, N_{p'(\omega')})
  \in \RR^d \times C(\RR_+, \RR^r) \times M(\RR_+ \times U)
 \end{align}
 and
 \begin{align}\label{mapping_2}
  \Omega \ni \omega
  \mapsto (\bxi(\omega), (\bW_t(\omega))_{t\in\RR_+}, N_{p(\omega)} )
  \in \RR^d \times C(\RR_+, \RR^r) \times M(\RR_+ \times U)
 \end{align}
 are
 \ $\cF'/ \cB(\RR^d) \otimes \cC(\RR_+, \RR^r)
          \otimes \cM(\RR_+ \times U)$-measurable and
 \[
   \sigma(\bxi,\bW_s,s\in\RR_+, p(s), s\in\RR_{++}\cap D(p))
   / \cB(\RR^d) \otimes \cC(\RR_+, \RR^r)
     \otimes \cM(\RR_+ \times U)\text{-measurable} ,
 \]
 respectively.
Further, they are
 \ $\cF'_t / \cB(\RR^d) \otimes \cC_t(\RR_+, \RR^r)
             \otimes \cM_t(\RR_+ \times U)$-measurable
 and
 \[
   \sigma(\bxi,\bW_s,s\in[0,t], p(s), s\in(0,t]\cap D(p))
   / \cB(\RR^d) \otimes \cC_t(\RR_+, \RR^r)
     \otimes \cM_t(\RR_+ \times U)\text{-measurable}
 \]
  for all \ $t \in \RR_+$, \ respectively.
Indeed, since \ $\bX_0'$ \ and \ $\bxi$ \ are \ $\cF'/\cB(\RR^d)$-measurable
 and \ $\sigma(\bxi)/\cB(\RR^d)$-measurable, respectively, by \eqref{C_infty}
 and \eqref{M_infty}, it is enough to check that for all \ $t \in \RR_+$,
 \ $n \in \NN$, \ $A_1 \in \cB(\RR^{nr})$, \ $t_1, \ldots, t_n \in [0, t]$,
 \ $s \in [0, t]$, \ $B \in \cB(U)$, \ $A_2 \in \cB([0, \infty])$,
 \begin{align*}
  &\Big\{ \omega'\in\Omega' : (\bW'_{t_1}(\omega'), \ldots,\bW'_{t_n}(\omega'))
   \in A_1 \Big\} \in \cF' , \\
  &\Big\{ \omega\in\Omega : (\bW_{t_1}(\omega), \ldots,\bW_{t_n}(\omega)) \in A_1 \Big\}
      \in \sigma(\bW_s,s\in\RR_+) , \\
  &\Big\{ \omega'\in\Omega' :  N_{p'(\omega')}([0,s]\times B)\in A_2 \Big\}
   \in \cF' , \\
  &\Big\{ \omega\in\Omega :  N_{p(\omega)}([0,s]\times B)\in A_2 \Big\}\in\sigma(p(s), s\in\RR_{++}\cap D(p)) .
 \end{align*}
These relations hold since \ $\bW'_{t_i}$, \ $i \in \{1, \ldots, n\}$, \ and
 \ $\bW_{t_i}$, \ $i \in \{1, \ldots, n\}$, \ are
 \ $\cF'/\cB(\RR^r)$-measurable and
 $\sigma(\bW_s, s \in \RR_+) / \cB(\RR^r)$-measurable, and \ $p'$ \ and \ $p$
 \ are \ $\cF'/\cM(\RR_+ \times U)$-measurable and $\sigma(p(s),
 \ s \in \RR_{++} \cap D(p)) / \cM(\RR_+ \times U)$-measurable, respectively.
Similarly, one can argue that the functions in question are
 \ $\cF'_t/ \cB(\RR^d) \otimes \cC_t(\RR_+, \RR^r)
            \otimes \cM_t(\RR_+ \times U)$-measurable
 and
 \ $\sigma(\bxi, \bW_s, s \in [0, t], p(s), s \in (0, t] \cap D(p))
    / \cB(\RR^d) \otimes \cC_t(\RR_+, \RR^r)
      \otimes \cM_t(\RR_+ \times U)$-measurable
 for all \ $t \in \RR_+$, \ respectively.

Next, we check that the process \ $\bX$ \ is adapted to the augmented
 filtration \ $(\cF_t^{\bxi,\bW,p})_{t\in\RR_+}$.
\ First, note that the process \ $\bX$ \ is adapted to
 \ $(\cF_t^{\bxi,\bW,p})_{t\in\RR_+}$ \ if and only if \ $\varphi_t(\bX)$ \ is
 \ $\cF_t^{\bxi,\bW,p}/\cD_t(\RR_+,\RR^d)$-measurable for all \ $t \in \RR_+$,
 \ where \ $\varphi_t$ \ is given in \eqref{varphi_t}.
Indeed,
 \begin{align*}
  &\text{$(\bX_t)_{t\in\RR_+}$ \ is \ $(\cF_t^{\bxi,\bW,p})_{t\in\RR_+}$-adapted}
     \qquad \Longleftrightarrow \qquad
  \text{$\sigma(\bX_t)\subset\cF_t^{\bxi,\bW,p}$ \ for all \ $t\in\RR_+$} \\
  &\qquad \qquad \Longleftrightarrow \qquad
  \text{$\sigma(\bX_s : s\in[0,t])\subset\cF_t^{\bxi,\bW,p}$ \ for all \ $t\in\RR_+$}\\
  &\qquad \qquad \Longleftrightarrow \qquad
  \text{$\varphi_t(\bX)$ \ is \ $\cF_t^{\bxi,\bW,p}/\cD_t(\RR_+,\RR^d)$-measurable for all \ $t\in\RR_+$},
 \end{align*}
 where the last equivalence can be checked as follows.
Since \ $\cD_t(\RR_+, \RR^d)$ \ coincides with the smallest \ $\sigma$-algebra
 containing the finite-dimensional cylinder sets of the form
 \[
   \big\{ y \in D(\RR_+,\RR^d) : (y(t_1), \ldots, y(t_n)) \in A \big\} ,
   \qquad n \in \NN , \quad A \in \cB(\RR^{nd}) , \quad
   t_1, \ldots, t_n \in [0,t] ,
 \]
 it is enough to check that
 \ $\sigma(\bX_s : s \in [0,t]) \subset \cF_t^{\bxi,\bW,p}$ \ for all
 \ $t \in \RR_+$ \ is equivalent with
 \begin{align*}
  \{ \omega \in \Omega
     : \big( (\varphi_t(\bX))_{t_1}(\omega), \cdots,
             (\varphi_t(\bX))_{t_n}(\omega) \big)
       \in A \}
  \in \cF_t^{\bxi,\bW,p}
 \end{align*}
 for all \ $n \in \NN$, \ $A \in \cB(\RR^{nd})$, \ $t_1, \ldots, t_n \in [0,t]$,
 $t \in \RR_+$, \ which readily follows from
 \begin{align*}
  \{ \omega \in \Omega
     : \big( (\varphi_t(\bX))_{t_1}(\omega), \cdots,
             (\varphi_t(\bX))_{t_n}(\omega) \big)
       \in A \}
  = \{ \omega \in \Omega
       : \big( \bX_{t_1}(\omega), \cdots, \bX_{t_n}(\omega) \big) \in A \} .
 \end{align*}
Since \ $\varphi_t(\bX) = \varphi_t \circ h' \circ (\bxi, \bW,p)$,
 \ $t \in \RR_+$, \ the mapping \ $\varphi_t$ \ is
 \ $\cD_t(\RR_+,\RR^d)/\cD_t(\RR_+,\RR^d)$-measurable for all \ $t \in \RR_+$,
 \ $h'$ \ is \ $\hcB_t / \cD_t(\RR_+, \RR^d)$-measurable for all
 \ $t \in \RR_+$, \ it remains to check that the mapping \eqref{mapping_2} is
 \ $\cF_t^{\bxi,\bW,p}/\widehat\cB_t$-measurable for all \ $t \in \RR_+$.
\ Recall that
 \begin{align*}
  &\hcB_t = \sigma\Big(\cB(\RR^d)\otimes \cC_t(\RR_+,\RR^r)\otimes\cM_t(\RR_+\times U) \cup \cN\Big),\qquad t\in\RR_+,\\
  &\cF_t^{\bxi,\bW,p} = \sigma\Big(  \sigma(\bxi,\bW_s,s\in[0,t], p(s), s\in(0,t]\cap D(p)) \cup \cN^{\bxi,\bW,p}\Big) ,
  \qquad t \in \RR_+ ,
 \end{align*}
 where
 \begin{align*}
  \cN = \Big\{ A \subset \RR^d \times C(\RR_+, \RR^r) \times M(\RR_+ \times U)
               : \, & \exists
             B\in\cB(\RR^d)\otimes \cC(\RR_+,\RR^r)\otimes\cM(\RR_+\times U) \\
               &\text{with} \; A\subset B, \; (n'\times P_{\bW,r}\times P_{U,m})(B)=0
          \Big\} ,
 \end{align*}
 and
 \begin{align*}
  \cN^{\bxi,\bW,p} := \Big\{ A \subset \Omega : \, & \exists  B\in \sigma(\bxi,\bW_s,s\in\RR_+, p(s), s\in\RR_{++}\cap D(p)) \\
                                     &\text{with} \; A \subset B, \; \PP(B)=0
          \Big\} .
 \end{align*}
Since a generator system of
 \ $\cB(\RR^d) \otimes \cC_t(\RR_+,\RR^r) \otimes\cM_t(\RR_+\times U)$
 \ together with \ $\cN$ \ is a generator system of \ $\hcB_t$, \ and we have
 already checked that the mapping \eqref{mapping_2} is
 \[
   \sigma(\bxi,\bW_s,s\in\RR_+, p(s), s\in\RR_{++}\cap D(p))/\cB(\RR^d) \otimes \cC(\RR_+, \RR^r)\otimes \cM(\RR_+ \times U)
   \text{-measurable,}
 \]
 it remains to verify that \ $(\bxi,\bW,p)^{-1}(A) \in \cF_t^{\bxi,\bW,p}$ \ for
 all \ $A \in \cN$ \ and \ $t \in \RR_+$.
\ We show that \ $(\bxi, \bW,p)^{-1}(A) \in \cN^{\bxi,\bW,p}$ \ for all
 \ $A \in \cN$, \ implying \ $(\bxi, \bW,p)^{-1}(A) \in \cF_t^{\bxi,\bW,p}$ \ for
 all \ $t \in \RR_+$, \ as desired.
If \ $A \in \cN$, \ then there exists
 \ $B \in \cB(\RR^d) \otimes \cC(\RR_+,\RR^r)
    \otimes \cM(\RR_+\times U)$
 \ such that \ $A \subset B$ \ and
 \ $(n' \times P_{\bW,r} \times P_{U,m})(B) = 0$.
\ Hence
 \[
   (\bxi,\bW,p)^{-1}(A)\subseteq (\bxi,\bW,p)^{-1}(B)
   \in \sigma(\bxi,\bW_s,s\in\RR_+, p(s), s\in\RR_{++}\cap D(p))
 \]
 and
 \begin{align*}
  \PP((\bxi,\bW,p)^{-1}(B)) = \PP( (\bxi,\bW,p) \in B)
                           = (n'\times P_{\bW,r} \times P_{U,m})(B)
                           = 0 ,
 \end{align*}
 where, for the last but one equality, we used that the distribution of
 \ $(\bxi, \bW,p)$ \ under \ $\PP$ \ is \ $n'\times P_{\bW,r} \times P_{U,m}$
 \ (as it was explained at the beginning of the proof).
By definition, this means that \ $(\bxi, \bW,p)^{-1}(A) \in \cN^{\bxi,\bW,p}$.

Next we check that \ $(\bX_t)_{t\in\RR_+}$ \ satisfies the SDE \eqref{SDE_X_YW}
 \ $\PP$-almost surely.
Since \ $h'$ \ is
 \ $\cB(\RR^d) \otimes \cC(\RR_+, \RR^r)\otimes \cM(\RR_+ \times U)
    /\cD(\RR_+, \RR^d)$-measurable,
 and the triplets \ $(\bX'_0, \bW', p')$ \ and \ $(\bxi, \bW, p)$ \ induce the
 same probability measure \ $n' \times P_{\bW,r} \times P_{U,m}$ \ on the
 measurable space
 \[
   \Big(\RR^d\times C(\RR_+,\RR^r) \times M(\RR_+\times U),
         \cB(\RR^d) \otimes \cC(\RR_+, \RR^r) \otimes \cM(\RR_+\times U)\Big)
 \]
 with respect to the probability measure \ $\PP'$ \ and \ $\PP$,
 \ respectively, the triplets \ $(\bX', \bW', p')$ \ and \ $(\bX, \bW, p)$
 \ induce the same probability measure on the measurable space
 \[
   \Big( D(\RR_+,\RR^d) \times C(\RR_+,\RR^r) \times M(\RR_+\times U),
         \cD(\RR_+, \RR^d) \otimes \cC(\RR_+, \RR^r) \otimes \cM(\RR_+\times U)\Big)
  \]
 with respect to the probability measure \ $\PP'$ \ and \ $\PP$,
 \ respectively.
Let us apply Lemma \ref{Lemma_distribution2} with the following choices
 \[
   \bigl( \Omega^{(1)}, \cF^{(1)}, (\cF_t^{(1)})_{t\in\RR_+}, \PP^{(1)}, \bW^{(1)},
          p^{(1)}, \bX^{(1)} \bigr)
   := \bigl( \Omega', \cF', (\cF_t')_{t\in\RR_+}, \PP', \bW', p', \bX' \bigr)
 \]
 and
 \[
   \bigl( \Omega^{(2)}, \cF^{(2)}, (\cF_t^{(2)})_{t\in\RR_+}, \PP^{(2)}, \bW^{(2)},
          p^{(2)}, \bX^{(2)} \bigr)
   := \bigl( \Omega, \cF, (\cF_t^{\bxi,\bW,p})_{t\in\RR_+}, \PP, \bW, p, \bX \bigr).
 \]
Since
 \ $\bigl( \Omega', \cF', (\cF'_t)_{t\in\RR_+}, \PP', \bW', p', \bX' \bigr)$
 \ is a weak solution of the SDE \eqref{SDE_X_YW} with initial distribution
 \ $n'$, \ the tuple
 \ $\bigl( \Omega^{(1)}, \cF^{(1)}, (\cF_t^{(1)})_{t\in\RR_+}, \PP^{(1)}, \bW^{(1)},
           p^{(1)}, \bX^{(1)} \bigr)$
 \ satisfies \textup{(D1)}, \textup{(D2)}, \textup{(D3)} and
 \textup{(D4)(b)--(e)}.
Further, as it was explained before Definition \ref{Def_strong_solution2}, the
 tuple
 \ $\bigl( \Omega^{(2)}, \cF^{(2)}, (\cF_t^{(2)})_{t\in\RR_+}, \PP^{(2)}, \bW^{(2)},
           p^{(2)}, \bX^{(2)} \bigr)$
 \ satisfies \textup{(D1)}, \textup{(D2)} and \textup{(D3)}, and we have
 already checked that \ $\bX$ \ is adapted to the augmented filtration
 \ $(\cF_t^{\bxi,\bW,p})_{t\in\RR_+}$.
\ Then Lemma  \ref{Lemma_distribution2} yields that the tuple
 \ $\bigl( \Omega^{(2)}, \cF^{(2)}, (\cF_t^{(2)})_{t\in\RR_+}, \PP^{(2)}, \bW^{(2)},
           p^{(2)}, \bX^{(2)} \bigr)$
 \ satisfies \textup{(D4)(b)--(d)} and the distribution of
 \begin{align*}
  \biggl(&\bX'_t - \bX'_0
          - \int_0^t b(s, \bX'_s) \, \dd s
          - \int_0^t \sigma(s, \bX'_s) \, \dd \bW'_s \\
         &- \int_0^t \int_{U_0} f(s, \bX'_{s-}, u) \, \tN'(\dd s, \dd u)
          - \int_0^t \int_{U_1}
             g(s, \bX'_{s-}, u) \, N'(\dd s, \dd u) \biggr)_{t\in\RR_+}
 \end{align*}
 on \ $(D(\RR_+, \RR^d), \cD(\RR_+, \RR^d))$ \ under \ $\PP'$ \ is the
 same as the distribution of
 \begin{align*}
  \biggl(&\bX_t - \bX_0
          - \int_0^t b(s, \bX_s) \, \dd s
          - \int_0^t \sigma(s, \bX_s) \, \dd \bW_s \\
         &- \int_0^t \int_{U_0} f(s, \bX_{s-}, u) \, \tN(\dd s, \dd u)
          - \int_0^t \int_{U_1}
             g(s, \bX_{s-}, u) \, N(\dd s, \dd u) \biggr)_{t\in\RR_+}
 \end{align*}
 on \ $(D(\RR_+, \RR^d), \cD(\RR_+, \RR^d))$ \ under \ $\PP$, \ where
 \ $N'(\dd s, \dd u)$ \ and \ $N(\dd s, \dd u)$ \ is the counting measure of
 \ $p'$ \ and \ $p$ \ on \ $\RR_+ \times U$, \ respectively, and
 \ $\tN'(\dd s, \dd u) := N'(\dd s, \dd u) - \dd s \, m(\dd u)$ \ and
 \ $\tN(\dd s, \dd u) := N(\dd s, \dd u) - \dd s \, m(\dd u)$.
\ Using that the first process and the identically 0 process are
 indistinguishable (since the SDE \eqref{SDE_X_YW} holds \ $\PP'$-a.s.\ for
 \ $(\bX'_t)_{t\in\RR_+}$), \ we obtain that the SDE \eqref{SDE_X_YW} holds
 \ $\PP$-a.s. for \ $(\bX_t)_{t\in\RR_+}$ \ as well, i.e., \textup{(D4)(e)}
 holds.

Finally, we show that \ $\PP(\bX_0 = \bxi) = 1$.
\ Since, as it was checked that the distribution of \ $\bX'$ \ and \ $\bX$
 \ coincide, especially, the distribution of \ $\bX'_0$ \ and \ $\bX_0$
 \ coincide, and consequently, the distribution of \ $\bX_0$ \ and \ $\bxi$
 \ coincide (both are equal to \ $n'$).
\ Using Corollary \ref{Corrolary1} for
 \ $\bigl( \Omega, \cF, (\cF_t^{\bxi,\bW,p})_{t\in\RR_+}, \PP, \bW, p, \bX \bigr)$
 \ (which is especially a weak solution of the SDE \eqref{SDE_X_YW} with
 initial distribution \ $n'$) \ we get
 \begin{align*}
  \PP(\bX_0 = \bxi)
   = \PP(\bxi + k'(\bxi,\bW,p)_0 = \bxi)
   = \PP(k'(\bxi,\bW,p)_0 = \bzero)
   = \PP(k'(\bxi,\bW,p)_0=\bY_0)
   = 1,
 \end{align*}
 as desired.

Summarizing, \ $(\bX_t)_{t\in\RR_+}$ \ is a strong solution of the SDE
 \eqref{SDE_X_YW} with initial value \ $\bxi$.
\proofend

\appendix

\section{Appendix}
\label{appendix}

Let \ $(\Omega, \cF, (\cF_t)_{t\in\RR_+}, \PP)$ \ be a filtered probability
 space.
First we recall the notion of \ $(\cF_t)_{t\in\RR_+}$-predictability, see, e.g.,
 Ikeda and Watanabe \cite[Chapter II, Definition 3.3]{IkeWat}.
The predictable \ $\sigma$-algebra \ $\cP$ \ on
 \ $\RR_+ \times \Omega \times U$ \ is given by
 \begin{align*}
  \cP := \sigma( h : \RR_+ \times \Omega \times U \to \RR
                 \mid
                 &\text{$h(t, \cdot, \cdot)$ \ is
                        \ $\cF_t \otimes \cB(U) / \cB(\RR)$-measurable for all
                        \ $t \in \RR_{++}$,} \\
                 &\text{$h(\cdot, \omega, u)$ \ is left continuous for all
                        \ $(\omega, u) \in \Omega \times U$} ) .
 \end{align*}
A function \ $H : \RR_+ \times \Omega \times U \to \RR^d$ \ is called
 \ $(\cF_t)_{t\in\RR_+}$-predictable if it is \ $\cP / \cB(\RR^d)$-measurable.

\begin{Lem}\label{Lemma_pred}
Let \ $(\Omega, \cF, (\cF_t)_{t\in\RR_+}, \PP)$ \ be a filtered probability
 space.
Let \ $(\bX_t)_{t\in\RR_+}$ \ be an \ $(\cF_t)_{t\in\RR_+}$-adapted c\`adl\`ag
 process with values in \ $\RR^d$.
\begin{enumerate}
 \item[\textup{(i)}]
  If \ $w : \RR^d \to \RR$ \ is a continuous function, then for each
   \ $T \in \RR_+$ \ and \ $B \in \cB(U)$, \ the function
   \ $h(t, \omega, u) := w(\bX_{t-}(\omega)) \bbone_{[0, T]}(t) \bbone_B(u)$,
   \ $(t, \omega, u) \in \RR_+ \times \Omega \times U$, \ is
   \ $(\cF_t)_{t\in\RR_+}$-predictable.
 \item[\textup{(ii)}]
  If \ $T \in \RR_+$, \ $A \in \cB(\RR^d)$ \ is an open set and
   \ $B \in \cB(U)$, \ then
   \[
     \{ (t, \omega, u) \in \RR_+ \times \Omega \times U
        : t \in [0, T] , \, \bX_{t-}(\omega) \in A , \, u \in B \}
     \in \cP .
   \]
 \item[\textup{(iii)}]
  If \ $f : \RR_+ \times \RR^d \times U \to \RR^d$ \ is
   \ $\cB(\RR_+) \otimes \cB(\RR^d) \otimes \cB(U) / \cB(\RR^d)$-measurable,
   then the function \ $H(t, \omega, u) := f(t, \bX_{t-}(\omega), u)$,
   \ $(t, \omega, u) \in \RR_+ \times \Omega \times U$, \ is
   \ $(\cF_t)_{t\in\RR_+}$-predictable.
\end{enumerate}
\end{Lem}

\noindent
\textbf{Proof.}
(i) \ The function \ $h$ \ is $(\cF_t)_{t\in\RR_+}$-predictable, since it
 belongs to the generator system of \ $\cP$.
\ Indeed, for each \ $t \in \RR_+$, \ the mapping
 \ $\Omega \times U \ni (\omega, u) \mapsto h(t, \omega, u)$
 \ is \ $\cF_t \otimes \cB(U) / \cB(\RR)$-measurable, because
 \ $\bX_s$ \ is \ $\cF_s / \cB(\RR^d)$-measurable and \ $\cF_s \subset \cF_t$
 \ for all \ $s < t$, \ and hence \ $\bX_{t-} := \lim_{s\uparrow t} \bX_s$ \ is
 \ $\cF_t / \cB(\RR^d)$-measurable, and \ $w$ \ is
 \ $\cB(\RR^d) / \cB(\RR)$-measurable.
Moreover, for each \ $(\omega, u) \in \Omega \times U$, \ the function
 \ $\RR_+ \ni t \mapsto h(t, \omega, u)$ \ is left continuous, because the
 functions \ $\RR_+ \ni t \mapsto \bbone_{[0, T]}(t)$ \ and
 \ $\RR_+ \ni t \mapsto \bX_{t-}(\omega)$ \ are left continuous and \ $w$ \ is
 continuous.

(ii) \ Consider the function \ $w_A : \RR^d \to \RR_+$ \ given by
 \ $w_A(\bx) := \varrho(\bx, \RR^d \setminus A)$, \ $\bx \in \RR^d$, \ where
 \ $\varrho$ \ denotes the Euclidean distance of \ $\bx$ \ and
 \ $\RR^d \setminus A$.
\ Then \ $w_A$ \ is continuous and \ $A = w_A^{-1}(\RR_{++})$.
\ Put
 \ $h_A(t, \omega, u) := w_A(\bX_{t-}(\omega)) \bbone_{[0, T]}(t) \bbone_B(u)$,
 \ $(t, \omega, u) \in \RR_+ \times \Omega \times U$.
\ Then, by (i), we obtain
 \begin{align*}
  &\{(t, \omega, u) \in \RR_+ \times \Omega \times U
     : t \in [0, T] , \, \bX_{t-}(\omega) \in A , \, u \in B \} \\
  &=\{(t, \omega, u) \in \RR_+ \times \Omega \times U
      : t \in [0, T] , \, w_A(\bX_{t-}(\omega)) \in \RR_{++} ,
         \, u \in B \} \\
  &=\{(t, \omega, u) \in \RR_+ \times \Omega \times U
       : h_A(t, \omega, u) \in \RR_{++} \}
  \in \cP .
 \end{align*}
(iii) \ We have \ $H = f \circ G$, \ where
 \ $G(t, \omega, u) := (t, \bX_{t-}(\omega), u)$,
 \ $(t, \omega, u) \in \RR_+ \times \Omega \times U$,
 \ thus it suffices to show that \ $G$ \ is
 \ $\cP / \cB(\RR_+) \otimes \cB(\RR^d) \otimes \cB(U)$-measurable.
The \ $\sigma$-algebra \ $\cB(\RR_+) \otimes \cB(\RR^d) \otimes \cB(U)$ \ is
 generated by the sets \ $[0, T] \times A \times B$ \ with \ $T \in \RR_+$,
 \ open sets \ $A \in \cB(\RR^d)$ \ and \ $B \in \cB(U)$, \ hence it suffices
 to show that
 \[
   \{(t, \omega, u) \in \RR_+ \times \Omega \times U
      : t \in [0, T] , \, \bX_{t-}(\omega) \in A , \, u \in B \} \in \cP .
 \]
This holds by (ii).
\proofend

Note that using Lemma \ref{Lemma_pred}, one can relax Assumption 6.2.8
 in Applebaum \cite{App}.

The next lemma plays a similar role as Lemma 139 in Situ \cite{Sit}.

\begin{Lem}\label{Lemma_distribution}
Let
 \ $\bigl( \Omega^{(i)}, \cF^{(i)}, (\cF_t^{(i)})_{t\in\RR_+}, \PP^{(i)}, \bW^{(i)},
           p^{(i)}, \bX^{(i)} \bigr)$,
 \ $i \in \{1, 2\}$, \ be tuples satisfying
 \textup{(D1)}, \textup{(D2)}, \textup{(D3)} and \textup{(D4)(b)--(d)}.
Suppose that \ $(\bW^{(1)}, p^{(1)}, \bX^{(1)})$ \ and
 \ $(\bW^{(2)}, p^{(2)}, \bX^{(2)})$ \ have the same distribution on
 \ $C(\RR_+, \RR^r) \times M(\RR_+ \times U) \times D(\RR_+, \RR^d)$.
\ Then
 \begin{align}\label{process1}
  \begin{split}
  \biggl(&\bX_t^{(1)} ,
          \int_0^t b(s, \bX_s^{(1)}) \, \dd s ,
          \int_0^t \sigma(s, \bX_s^{(1)}) \, \dd \bW_s^{(1)} , \\
  &\int_0^t \int_{U_0} f(s, \bX_{s-}^{(1)}, u) \, \tN^{(1)}(\dd s, \dd u) ,
   \int_0^t \int_{U_1}
    g(s, \bX_{s-}^{(1)}, u) \, N^{(1)}(\dd s, \dd u) \biggr)_{t\in\RR_+}
  \end{split}
 \end{align}
 and
 \begin{align}\label{process2}
 \begin{split}
  \biggl(&\bX_t^{(2)} ,
          \int_0^t b(s, \bX_s^{(2)}) \, \dd s ,
          \int_0^t \sigma(s, \bX_s^{(2)}) \, \dd \bW_s^{(2)} , \\
  &\int_0^t \int_{U_0} f(s, \bX_{s-}^{(2)}, u) \, \tN^{(2)}(\dd s, \dd u) ,
   \int_0^t \int_{U_1}
    g(s, \bX_{s-}^{(2)}, u) \, N^{(2)}(\dd s, \dd u) \biggr)_{t\in\RR_+}
    \end{split}
 \end{align}
 have the same distribution on \ $(D(\RR_+, \RR^d))^5$, \ where, for each
 \ $i \in \{1, 2\}$, \ $N^{(i)}(\dd s, \dd u)$ \ is the counting measure of
 \ $p^{(i)}$ \ on \ $\RR_{++} \times U$, \ and
 \ $\tN^{(i)}(\dd s, \dd u) := N^{(i)}(\dd s, \dd u) - \dd s \, m(\dd u)$.
\end{Lem}

\noindent
\textbf{Proof.}
By Remark \ref{intexist}, the above processes have c\`adl\`ag modifications.
According to Lemma VI.3.19 in Jacod and Shiryaev \cite{JacShi}, it suffices to
 show that the finite dimensional distributions of the above processes
 coincide.

By Proposition I.4.44 in Jacod and Shiryaev \cite{JacShi}, for each
 \ $i \in \{1, 2\}$ \ and \ $t \in \RR_+$,
 \ $I_{1,n}^{(i)}(t) \stochi \int_0^t b(s, \bX_s^{(i)}) \, \dd s$ \ and
 \ $I_{2,n}^{(i)}(t) \stochi \int_0^t \sigma(s, \bX_s^{(i)}) \, \dd \bW_s^{(i)}$
 \ as \ $n \to \infty$, \ where
 \[
   I_{1,n}^{(i)}(t)
   := \frac{1}{n}
      \sum_{k=1}^\nt
       b\left( \frac{k-1}{n}, \bX_{\frac{k-1}{n}}^{(i)} \right) , \qquad
   I_{2,n}^{(i)}(t)
   := \sum_{k=1}^\nt
       \sigma\left( \frac{k-1}{n}, \bX_{\frac{k-1}{n}}^{(i)} \right)
       \left( \bW_{\frac{k}{n}}^{(i)} - \bW_{\frac{k-1}{n}}^{(i)} \right) .
 \]
Let \ $U_{1,j} \in \cB(U)$, \ $j \in \NN$, \ be such that they are disjoint,
 \ $m(U_{1,j}) < \infty$, \ $j \in \NN$, \ and
 \ $U_1 = \bigcup_{j=1}^\infty U_{1,j}$ \ (such a sequence exists since \ $m$ \ is
 \ $\sigma$-finite, see, e.g., Cohn \cite[page 9]{Coh}).
Then for each
 \ $i \in \{1, 2\}$ \ and \ $t \in \RR_+$,
 \ $I_{3,n}^{(i)}(t)
    \to \int_0^t \int_{U_1} g(s, \bX_{s-}^{(i)}, u) \, N^{(i)}(\dd s, \dd u)$
 \ as \ $n \to \infty$ \ $\PP^{(i)}$-almost surely, where
 \begin{align*}
  I_{3,n}^{(i)}(t)
  &:= \sum_{j=1}^n
       \int_0^t \int_{U_{1,j}} g(s, \bX_{s-}^{(i)}, u) \, N^{(i)}(\dd s, \dd u)
   = \sum_{j=1}^n
      \sum_{s\in(0,t]\cap D(p_{1,j}^{(i)})} g(s, \bX_{s-}^{(i)}, p_{1,j}^{(i)}(s)) ,
 \end{align*}
 where \ $p_{1,j}^{(i)}$ \ denotes the thinning of \ $p^{(i)}$ \ onto \ $U_{1,j}$,
 \ see, e.g., Ikeda and Watanabe \cite[page 62]{IkeWat}.
Since \ $m(U_{1,j}) < \infty$, \ by Remark \ref{Rem_sum_finite}, the set \ $(0,t] \cap D(p_{1,j}^{(i)})$ \ is
 finite \ $\PP^{(i)}$-almost surely for all \ $t \in \RR_+$ \ and
 \ $i \in \{1, 2\}$, \ $j \in \NN$, \ and hence one can order the set
 \ $D(p_{1,j}^{(i)})$ \ according to magnitude, say
 \ $0 < \zeta^{(i)}_{1,j,1} < \zeta^{(i)}_{1,j,2} < \cdots$, \ $j \in \NN$,
 \ $i \in \{1, 2\}$.
\ Namely,
 \begin{align}\label{def_zeta}
   \zeta^{(i)}_{1,j,\ell}
   = \inf\{t\in\RR_+ : N^{(i)}((0,t]\times U_{1,j}) \geq \ell\},
   \qquad \ell \in \NN , \quad j \in \NN , \quad i \in \{1, 2\}
 \end{align}
 on the event
 \[
     \Omega^{(i)}_{1,j}
     :=\bigcap_{k=1}^\infty
        \Big\{ \omega \in \Omega^{(i)}
               : N_{p^{(i)}_{1,j}(\omega)}((0,k]\times U_{1,j}) < \infty \Big\} ,
     \qquad j \in \NN ,\quad i\in\{1,2\},
 \]
 having \ $\PP^{(i)}$-probability 1, where we used that the point measure
 corresponding to the point function \ $p^{(i)}_{1,j}(\omega)$ \ is its counting
 measure \ $N_{p^{(i)}_{1,j}(\omega)}$, \ see Section \ref{section_preliminaries}.
Then we can write \ $I_{3,n}^{(i)}(t)$ \ in the form
 \begin{align*}
  I_{3,n}^{(i)}(t)
  = \sum_{j=1}^n
     \sum_{\ell=1}^\infty
      g\Bigl(\zeta^{(i)}_{1,j,\ell}, \bX_{\zeta^{(i)}_{1,j,\ell}-}^{(i)},
             p_{1,j}^{(i)}(\zeta^{(i)}_{1,j,\ell})\Bigr)
      \bbone_{(0,t]}(\zeta^{(i)}_{1,j,\ell}) ,
  \qquad t \in \RR_+ , \quad n \in \NN , \quad i \in \{1, 2\} ,
 \end{align*}
 where
 \ $\sum_{\ell=1}^\infty
     g\Bigl(\zeta^{(i)}_{1,j,\ell}, \bX_{\zeta^{(i)}_{1,j,\ell}-}^{(i)},
            p_{1,j}^{(i)}(\zeta^{(i)}_{1,j,\ell})\Bigr)
     \bbone_{(0,t]}(\zeta^{(i)}_{1,j,\ell})$
 \ is a finite sum \ $\PP^{(i)}$-almost surely.
Furthermore, by Remark \ref{intexist}, for \ $i \in \{1, 2\}$ \ and
 \ $t \in \RR_+$,
 \ $I_{4,n}^{(i)}(t)
    \to \int_0^t \int_{U_0} f(s, \bX_{s-}^{(i)}, u) \, \tN^{(i)}(\dd s, \dd u)$
 \ as \ $n \to \infty$ \ $\PP^{(i)}$-almost surely, where
 \begin{align*}
  I_{4,n}^{(i)}(t)
  := \int_0^t \int_{U_0}
      \bbone_{[0,\tau_n^{(i)}]}(s)f(s, \bX_{s-}^{(i)}, u) \, \tN^{(i)}(\dd s, \dd u)
 \end{align*}
 with
 \[
   \tau_n^{(i)}
   := \inf\left\{t \in \RR_+
                 : \int_0^t \int_{U_0}
                    \|f(s, \bX_s^{(i)}, u)\|^2 \, \dd s \, m(\dd u)
                   \geq n\right\}
      \land n ,
   \qquad n \in \NN , \quad i \in \{1, 2\} ,
 \]
 satisfying  \ $\tau_n^{(i)} \uparrow \infty$ \ $\PP^{(i)}$-almost surely as
 \ $n \to \infty$.
\ Let \ $U_{0,j} \in \cB(U)$, \ $j \in \NN$, \ be such that they are disjoint,
 \ $m(U_{0,j}) < \infty$, \ $j \in \NN$, \ and
 \ $U_0 = \bigcup_{j=1}^\infty U_{0,j}$ \ (such a sequence exists since \ $m$ \ is
 \ $\sigma$-finite, see, e.g., Cohn \cite[page 9]{Coh}).
Then, by pages 47 and 63 in Ikeda and Watanabe \cite{IkeWat},
 for all \ $t \in \RR_+$, \ $i \in \{1, 2\}$ \ and \ $n \in \NN$,
 \ $I_{4,n,j}^{(i)}(t) \stochi I_{4,n}^{(i)}(t)$ \ as \ $j \to \infty$, \ where
 \begin{align*}
  I_{4,n,j}^{(i)}(t)
  &:=\int_0^t \int_{U_0}
      \bbone_{(-j,j)}\left(\bbone_{[0,\tau_n^{(i)}]}(s)f(s, \bX_{s-}^{(i)}, u)\right)
      \bbone_{U_{0,j}}(u) \bbone_{[0,\tau_n^{(i)}]}(s)f(s, \bX_{s-}^{(i)}, u)
      \, \tN^{(i)}(\dd s, \dd u) \\
  &=\int_0^t \int_{U_{0,j}}
     \bbone_{(-j,j)}\left(f(s, \bX_{s-}^{(i)}, u)\right)
     \bbone_{[0,\tau_n^{(i)}]}(s)f(s, \bX_{s-}^{(i)}, u) \, \tN^{(i)}(\dd s, \dd u) .
 \end{align*}
By page 62 in Ikeda and Watanabe \cite{IkeWat}, for all \ $t \in \RR_+$,
 \ $i \in \{1, 2\}$, \ $n \in \NN$, \ and \ $j \in \NN$,
 \ $I_{4,n,j}^{(i)}(t) = I_{4,n,j}^{(i),a}(t) - I_{4,n,j}^{(i),b}(t)$, \ where
 \begin{align*}
  I_{4,n,j}^{(i),a}(t)
  &:=\int_0^t \int_{U_{0,j}}
      \bbone_{(-j,j)}\left(f(s, \bX_{s-}^{(i)}, u)\right)
      \bbone_{[0,\tau_n^{(i)}]}(s)f(s, \bX_{s-}^{(i)}, u)
      \, N^{(i)}(\dd s, \dd u) , \\
  I_{4,n,j}^{(i),b}(t)
  &:=\int_0^t
      \left(\int_{U_{0,j}} \bbone_{(-j,j)}\left(f(s, \bX_{s-}^{(i)}, u)\right)
      \bbone_{[0,\tau_n^{(i)}]}(s)f(s, \bX_{s-}^{(i)}, u)
      \, m(\dd u)\right) \dd s .
 \end{align*}
Similarly as for the integrals
 \ $ \int_0^t \int_{U_1} g(s, \bX_{s-}^{(i)}, u) \, N^{(i)}(\dd s, \dd u)$ \ and
 \ $ \int_0^t b(s, \bX_s^{(i)}) \, \dd s$, \ there exist sequences of random
 variables  \ $(I_{4,n,j,\ell}^{(i),a}(t))_{\ell\in\NN}$ \ and
 \ $(I_{4,n,j,\ell}^{(i),b}(t))_{\ell\in\NN}$ \ such that
 \ $I_{4,n,j,\ell}^{(i),a}(t) \stochi I_{4,n,j}^{(i),a}(t)$ \ and
 \ $I_{4,n,j,\ell}^{(i),b}(t) \stochi I_{4,n,j}^{(i),b}(t)$ \ as \ $\ell\to\infty$,
 \ respectively.
Then, for all \ $t \in \RR_+$ \ and \ $i \in \{1, 2\}$,
 \ $I_{4,n,j,\ell}^{(i),a}(t) - I_{4,n,j,\ell}^{(i),b}(t)
    \stochi \int_0^t \int_{U_0} f(s, \bX_{s-}^{(i)}, u) \, \tN^{(i)}(\dd s, \dd u)$
 \ as \ $\ell \to \infty$, \ then \ $j \to \infty$, \ and, finally,
 \ $n \to \infty$.
\ Using part (vi) of Theorem 2.7 in van der Vaart \cite{Vaart}, we get for all
 \ $K \in \NN$, \ $t_1, \ldots, t_K \in \RR_+$ \ and \ $i \in \{1, 2\}$,
 \begin{align*}
  &\left( \bX_{t_k}^{(i)}, I_{1,n}^{(i)}(t_k), I_{2,n}^{(i)}(t_k), I_{3,n}^{(i)}(t_k),
          I_{4,n,j,\ell}^{(i),a}(t_k)
          - I_{4,n,j,\ell}^{(i),b}(t_k) \right)_{k\in\{1,\ldots,K\}} \\
  &\stochi
   \Bigl( \bX_{t_k}^{(i)}, \int_0^{t_k} b(s, \bX_s^{(i)}) \, \dd s,
          \int_0^{t_k} \sigma(s, \bX_s^{(i)}) \, \dd \bW_s^{(i)}, \\
  &\phantom{\stochi\Bigl(}
          \int_0^{t_k} \int_{U_1} g(s, \bX_{s-}^{(i)}, u) \, N^{(i)}(\dd s, \dd u),
          \int_0^{t_k} \int_{U_0}
           f(s, \bX_{s-}^{(i)}, u) \,
           \tN^{(i)}(\dd s, \dd u) \Bigr)_{k\in\{1,\ldots,K\}}
 \end{align*}
 as \ $\ell, j, n \to \infty$.
\ Since \ $(\bW^{(1)}, p^{(1)}, \bX^{(1)})$ \ and \ $(\bW^{(2)}, p^{(2)}, \bX^{(2)})$
 \ have the same distribution, the random vectors
 \[
   \left( \bX_{t_k}^{(1)}, I_{1,n}^{(1)}(t_k), I_{2,n}^{(1)}(t_k), I_{3,n}^{(1)}(t_k),
          I_{4,n,j,\ell}^{(1),a}(t_k)
          - I_{4,n,j,\ell}^{(1),b}(t_k) \right)_{k\in\{1,\ldots,K\}}
 \]
 and
 \[
   \left( \bX_{t_k}^{(2)}, I_{1,n}^{(2)}(t_k), I_{2,n}^{(2)}(t_k), I_{3,n}^{(2)}(t_k),
          I_{4,n,j,\ell}^{(2),a}(t_k)
          - I_{4,n,j,\ell}^{(2),b}(t_k) \right)_{k\in\{1,\ldots,K\}}
 \]
 have the same distribution for all \ $\ell, j, n \in \NN$, \ as well.
Indeed, the random vectors above can be considered as some appropriate
 measurable function of \ $(\bW^{(1)}, p^{(1)}, \bX^{(1)})$ \ and
 \ $(\bW^{(2)}, p^{(2)}, \bX^{(2)})$, \ respectively.
For this, it is enough to verify that each coordinate of the above random
 vectors can be considered as some appropriate measurable function of
 \ $(\bW^{(1)}, p^{(1)}, \bX^{(1)})$ \ and \ $(\bW^{(2)}, p^{(2)}, \bX^{(2)})$,
 \ respectively, hence we fix \ $k \in \{1, \ldots, K\}$.
\begin{itemize}
 \item
  First observe, that \ $\bX_{t_k}^{(i)}$ \ is a
   \ $\cD(\RR_+, \RR^d) / \cB(\RR^d)$-measurable function of
   \ $\bX^{(i)}$, \ namely, \ $\bX_{t_k}^{(i)} = \Psi_0(\bX^{(i)})$, \ where
   \ $\Psi_0 : D(\RR_+, \RR^d) \to \RR^d$ \ is given by
   \ $\Psi_0(y) := y(t_k)$, \ $y \in D(\RR_+, \RR^d)$.

 \item
  Next, \ $I_{1,n}^{(i)}(t_k)$ \ is a
   \ $\cD(\RR_+, \RR^d) / \cB(\RR^d)$-measurable function of
   \ $\bX^{(i)}$ \ as well, namely,
   \ $I_{1,n}^{(i)}(t_k) = \Psi_1(\bX^{(i)})$, \ where
   \ $\Psi_1 : D(\RR_+, \RR^d) \to \RR^d$ \ is given by
   \ $\Psi_1(y)
      := \frac{1}{n}
         \sum_{k=1}^{\lfloor n t_k \rfloor}
          b\left( \frac{k-1}{n}, y\left(\frac{k-1}{n}\right) \right)$,
   \ $y \in D(\RR_+, \RR^d)$.

 \item
  In a similar way, \ $I_{2,n}^{(i)}(t_k)$ \ is a
   \ $\cD(\RR_+,\RR^d)\times\cC(\RR_+,\RR^r)/\cB(\RR^d)$-measurable function of
   \ $(\bX^{(i)}, \bW^{(i)})$, \ namely,
   \ $I_{2,n}^{(i)}(t_k) = \Psi_2(\bX^{(i)}, \bW^{(i)})$, \ where
   \ $\Psi_2 : D(\RR_+, \RR^d) \times C(\RR_+, \RR^r) \to \RR^d$ \ is given by
   \ $\Psi_2(y, w)
      := \sum_{k=1}^{\lfloor n t_k \rfloor}
          \sigma\left( \frac{k-1}{n}, y\left(\frac{k-1}{n}\right) \right)
          \left( w\left(\frac{k}{n}\right)
                 - w\left(\frac{k-1}{n}\right)\right)$,
   \ $y \in D(\RR_+, \RR^d)$, \ $w \in C(\RR_+, \RR^r)$.

 \item
  Now we show that \ $I_{3,n}^{(i)}(t_k)$ \ is a
   \ $\cD(\RR_+, \RR^d) \otimes \cM(\RR_+\times U) / \cB(\RR^d)$-measurable
   function of \ $(\bX^{(i)},p^{(i)})$.
  \ As a first step, we show that for each \ $j, \ell \in \NN$ \ there exist
   functions \ $\Phi_{j,\ell} : M(\RR_+\times U) \to \RR_+$ \ and
   \ $\Xi_{j,\ell} : M(\RR_+\times U) \to U$ \ such that
   \ $\Phi_{j,\ell}$ \ is \ $\cM(\RR_+\times U) / \cB(\RR_+)$-measurable,
   \ $\Xi_{j,\ell}$ \ is \ $\cM(\RR_+\times U) / \cB(U)$-measurable, and
   \ $(\zeta^{(i)}_{1,j,\ell}, p^{(i)}_{1,j}(\zeta^{(i)}_{1,j,\ell}))
      = (\Phi_{j,\ell}(N_{p^{(i)}_{1,j}}), \Xi_{j,\ell}(N_{p^{(i)}_{1,j}}))$
   \ holds \ $\PP^{(i)}$-almost surely.
  Then it will follow that \ $I_{3,n}^{(i)}(t_k) = \Psi_3(\bX^{(i)}, p^{(i)})$,
   \ where \ $\Psi_3 : D(\RR_+,\RR^d) \times M(\RR_+\times U) \to \RR^d$
   \ given by
   \begin{align*}
    \Psi_3(y,\pi)
    := \sum_{j=1}^n
       \sum_{\ell=1}^\infty
        g\Bigl( \Phi_{j,\ell}(\pi), y(\Phi_{j,\ell}(\pi)-) ,
                 \Xi_{j,\ell}(\pi) \Bigr)
      \bbone_{(0,t_k]}(\Phi_{j,\ell}(\pi))
   \end{align*}
   for \ $(y,\pi)\in D(\RR_+,\RR^d) \times M(\RR_+\times U)$ \ is
   \ $\cD(\RR_+,\RR^d) \otimes \cM(\RR_+\times U) / \cB(\RR^d)$-measurable.
  To prove the existence of \ $\Phi_{j,\ell}$ \ and \ $\Xi_{j,\ell}$, \ first we
   verify that \ $(\zeta^{(i)}_{1,j,\ell}, p^{(i)}_{1,j}(\zeta^{(i)}_{1,j,\ell}))$ \ is
   measurable with respect to the \ $\sigma$-algebra
   \ $\sigma\big( N_{p^{(i)}_{1,j}} \big) \cap \Omega^{(i)}_{1,j}$ \ having the form
   \begin{align}\label{help_szigma_algebra}
        \sigma\Big( \Big\{ \omega\in\Omega^{(i)}_{1,j}
                         : N_{p^{(i)}_{1,j}(\omega)}((0,t]\times B) =k \Big\} \,
                  \Big| \, t \in \RR_+ , \; B \in \cB(U_{1,j}) , \;
                              k \in \NN \Big) .
   \end{align}
   We have
   \begin{align}\label{help_zeta}
   \begin{split}
    &\bigl\{ \omega \in \Omega^{(i)}_{1,j}
             : \bigl( \zeta^{(i)}_{1,j,\ell}(\omega) ,
                      p^{(i)}_{1,j}(\omega)(\zeta^{(i)}_{1,j,\ell}(\omega)) \bigr)
               \in (0, t] \times B \bigr\} \\
    &=\bigcap_{n=1}^\infty \bigcup_{k=1}^n
       \Bigl\{ \omega \in \Omega^{(i)}_{1,j}
                : N_{p^{(i)}_{1,j}(\omega)}((0, (k-1)t/n] \times U_{1,j})
                  \leq \ell - 1 , \\[-7mm]
    &\phantom{=\bigcap_{n=1}^\infty \bigcup_{k=1}^n
                \Bigl\{ \omega \in \Omega^{(i)}_{1,j} : \;\;}
                 N_{p^{(i)}_{1,j}(\omega)}(((k-1)t/n, kt/n] \times B)
                 \geq 1 , \\[-7mm]
    &\phantom{=\bigcap_{n=1}^\infty \bigcup_{k=1}^n
                \Bigl\{ \omega \in \Omega^{(i)}_{1,j} : \;\;}
                 N_{p^{(i)}_{1,j}(\omega)}((0, kt/n] \times U_{1,j})
                 \geq \ell \Bigr\}
    \end{split}
   \end{align}
  for \ $t \in \RR_{++}$, \ $j, \ell \in \NN$, \ $B \in \cB(U_{1,j})$,
  \ $i \in \{1, 2\}$.
 \ Indeed, on the one hand, if \ $\omega \in \Omega^{(i)}_{1,j}$ \ is such that
  \ $\zeta^{(i)}_{1,j,\ell}(\omega) \in (0, t]$ \ and
  \ $p^{(i)}_{1,j}(\omega)(\zeta^{(i)}_{1,j,\ell}(\omega))\in  B$, \ then for each
  \ $n \in \NN$, \ there exists a unique \ $k \in \{1, \ldots, n\}$ \ with
  \ $\zeta^{(i)}_{1,j,\ell}(\omega) \in ((k-1)t/n, kt/n]$, \ and hence
  \ $N_{p^{(i)}_{1,j}(\omega)}((0, (k-1)t/n] \times U_{1,j}) \leq \ell - 1$,
  \ $N_{p^{(i)}_{1,j}(\omega)}(((k-1)t/n, kt/n] \times B) \geq 1$ \ and
  \ $N_{p^{(i)}_{1,j}(\omega)}((0, kt/n] \times U_{1,j}) \geq \ell$.
 \ On the other hand,
 \begin{align*}
  \{\omega\in \Omega^{(i)}_{1,j} : \;\; \zeta^{(i)}_{1,j,\ell}(\omega) \notin (0, t] \}
   & = \{\omega\in \Omega^{(i)}_{1,j} : \;\; N_{p^{(i)}_{1,j}(\omega)}((0, t] \times U_{1,j}) \leq \ell - 1 \} \\
   & \subset
    \bigcup_{n=1}^\infty \bigcap_{k=1}^n
       \Bigl\{ \omega \in \Omega^{(i)}_{1,j} : \;\;
         N_{p^{(i)}_{1,j}(\omega)}((0, kt/n] \times U_{1,j})
                 \leq \ell - 1 \Bigr\},
 \end{align*}
 and
 \begin{align*}
  &\bigl\{ \omega \in \Omega^{(i)}_{1,j}
             :  \zeta^{(i)}_{1,j,\ell}(\omega) \in (0, t] ,\;
                      p^{(i)}_{1,j}(\omega)(\zeta^{(i)}_{1,j,\ell}(\omega))
                     \notin  B \bigr\} \\
  &\qquad\subset
    \bigcup_{n=1}^\infty \bigcap_{k=1}^n
      \biggl( \Bigl\{ \omega \in \Omega^{(i)}_{1,j} : \;\;
    N_{p^{(i)}_{1,j}(\omega)}((0, (k-1)t/n] \times U_{1,j}) \geq \ell \Bigl\}\\[-7mm]
  &\phantom{\qquad\subset \bigcup_{n=1}^\infty \bigcap_{k=1}^n \Biggl( }
      \;\cup \Bigl\{ \omega \in \Omega^{(i)}_{1,j} : \;\; N_{p^{(i)}_{1,j}(\omega)}(((k-1)t/n, kt/n] \times B) = 0 \Bigl\}\\[-7mm]
  &\phantom{\qquad\subset \bigcup_{n=1}^\infty \bigcap_{k=1}^n \Biggl( }
      \;\cup \Bigl\{ \omega \in \Omega^{(i)}_{1,j} : \;\; N_{p^{(i)}_{1,j}(\omega)}((0, kt/n] \times U_{1,j}) \leq \ell - 1
   \Bigr\} \biggr).
 \end{align*}
For the second inclusion, for each \ $\omega \in \Omega^{(i)}_{1,j}$, \ let us
 choose \ $n(\omega) \in \NN$ \ such that
 \[
   n(\omega)
     > \max\left( \frac{1}{\zeta^{(i)}_{1,j,\ell}(\omega) - \zeta^{(i)}_{1,j,\ell-1}(\omega)} ,
                  \frac{1}{\zeta^{(i)}_{1,j,\ell+1}(\omega) - \zeta^{(i)}_{1,j,\ell}(\omega)} \right) .
 \]
If \ $\omega \in \Omega^{(i)}_{1,j}$ \ is such that
 \ $\zeta^{(i)}_{1,j,\ell}(\omega) \in (0, t]$ \ and
  \ $p^{(i)}_{1,j}(\omega)(\zeta^{(i)}_{1,j,\ell}(\omega)) \notin B$, \ then there
 exists a unique \ $k^* \in \{1, \ldots, n\}$ \ with
  \ $\zeta^{(i)}_{1,j,\ell}(\omega) \in ((k^*-1)t/n, k^*t/n]$, \ and hence we have
 \ $N_{p^{(i)}_{1,j}(\omega)}((0, kt/n] \times U_{1,j}) \leq \ell-1$ \ for
 \ $k \in \{1, \ldots, k^*-1\}$,
 \ $N_{p^{(i)}_{1,j}(\omega)}(((k^*-1)t/n, k^*t/n] \times B) = 0$, \ and
 \ $N_{p^{(i)}_{1,j}(\omega)}((0, (k-1)t/n] \times U_{1,j}) \geq \ell$ \ for
 \ $k \in \{k^*+1, \ldots, n\}$.

Since the set on right hand side of \eqref{help_zeta} is in the
 $\sigma$-algebra given in \eqref{help_szigma_algebra} and
 \ $\{(0,t]  \times B : t\in\RR_+ ,\; B\in \cB(U_{1,j})\}$ \ is a generator system
  of \ $\cB(\RR_+)\otimes \cB(U_{1,j})$, \ we readily get that the random
  variable  \ $(\zeta^{(i)}_{1,j,\ell}, p^{(i)}_{1,j}(\zeta^{(i)}_{1,j,\ell}))$ \
  is measurable with respect to the \ $\sigma$-algebra given in \eqref{help_szigma_algebra}.
  Let us apply Theorem 4.2.8 in Dudley \cite{Dud} with the following choices:
   \begin{itemize}
     \item[$\circ$] $X:=\Omega^{(i)}_{1,j}$, \ $Y:=M(\RR_+\times U)$,
     \item[$\circ$] $T:\Omega^{(i)}_{1,j} \to M(\RR_+\times U)$,
                    \ $T(\omega):=N_{p^{(i)}_{1,j}(\omega)}$,
           \ $\omega\in\Omega^{(i)}_{1,j}$,
     \item[$\circ$] $f:\Omega^{(i)}_{1,j} \to \RR_+\times U$, \ $f(\omega):= (\zeta^{(i)}_{1,j,\ell}(\omega),
             p^{(i)}_{1,j}(\omega)(\zeta^{(i)}_{1,j,\ell}(\omega)))$, \ $\omega\in\Omega^{(i)}_{1,j}$.
   \end{itemize}
  Then there exist functions
  \ $\Phi_{j,\ell}: M(\RR_+\times U) \to \RR_+$ \ and
  \ $\Xi_{j,\ell}: M(\RR_+\times U) \to U$ \ such that
  \ $\Phi_{j,\ell}$ \ is \ $\cM(\RR_+\times U)/\cB(\RR_+)$-measurable,
  \ $\Xi_{j,\ell}$ \ is \ $\cM(\RR_+\times U)/\cB(U)$-measurable, and
  \ $(\zeta^{(i)}_{1,j,\ell}, p^{(i)}_{1,j}(\zeta^{(i)}_{1,j,\ell}))
      = (\Phi_{j,\ell}(N_{p^{(i)}_{1,j}}), \Xi_{j,\ell}(N_{p^{(i)}_{1,j}}))$
  \ holds on \ $\Omega^{(i)}_{1,j}$.
  \ Since \ $\PP^{(i)}(\Omega^{(i)}_{1,j})=1$, \ we have
  \ $(\zeta^{(i)}_{1,j,\ell}, p^{(i)}_{1,j}(\zeta^{(i)}_{1,j,\ell}))
      = (\Phi_{j,\ell}(N_{p^{(i)}_{1,j}}), \Xi_{j,\ell}(N_{p^{(i)}_{1,j}}))$
  \ $\PP^{(i)}$-almost surely, as desired.

  In what follows we provide an alternative argument for verifying that \ $\zeta_{1,j,\ell}^{(i)}$ \ is an
   \ $\cM(\RR_+ \times U) / \cB(\RR)$-measurable function of
   \ $p^{(i)}$ \ with the advantage that the measurable function in question shows up explicitly.
   We have
   \ $\zeta_{1,j,1}^{(i)} = \inf\{ t \in \RR_{++} : |\Delta y_{i,j}(t)| > 1/2 \}$,
   \ where \ $y_{i,j}(t) := N^{(i)}((0,t]\times U_{1,j})$ \ and
   \ $\Delta y_{i,j}(t) := y_{i,j}(t) - y_{i,j}(t-) = N^{(i)}(\{t\}\times U_{1,j})$
   \ for \ $t \in \RR_{++}$.
  \ Further,
   \ $\zeta_{1,j,\ell+1}^{(i)}
      = \inf\{ t \in (\zeta_{1,j,\ell}^{(i)}, \infty)
               : |\Delta y_{i,j}(t)| > 1/2 \}$
   \ for all \ $\ell \in \NN$.
  \ Consider the mappings \ $\Psi_{3,\ell} : D(\RR_+, \RR) \to \RR_+$,
   \ $\ell \in \NN$, \ defined by
   \ $\Psi_{3,1}(y) := \inf\{ t \in \RR_{++} : |\Delta y(t)| > 1/2 \}$ \ and
   \ $\Psi_{3,\ell+1}(y)
      := \inf\{ t \in (\Psi_{3,\ell}(y), \infty) : |\Delta y(t)| > 1/2 \}$,
   \ $y \in D(\RR_+, \RR)$, \ $\ell \in \NN$.
  \ By Proposition VI.2.7 in Jacod and Shiryaev \cite{JacShi}, the mappings
   \ $\Psi_{3,\ell}$, \ $\ell \in \NN$, \ are continuous at each point
   \ $y \in D(\RR_+, \RR)$ \ such that \ $|\Delta y(t)| \ne 1/2$ \ for all
   \ $t \in \RR_+$.
  \ Moreover, we have
   \ $\zeta_{1,j,\ell}^{(i)} = \Psi_{3,\ell}(\Psi_{4,j}(p^{(i)}))$, \ where the
   mappings \ $\Psi_{4,j} : M(\RR_+ \times U) \to D(\RR_+, \RR)$,
   \ $j \in \NN$, \ are given by
   \ $\Psi_{4,j}(\pi) := (\pi((0,t]\times U_{1,j}))_{t\in\RR_+}$,
   \ $\pi \in M(\RR_+ \times U)$.
  \ Observe, that for each \ $\pi \in M(\RR_+ \times U)$, \ we have
   \ $|\Delta \Psi_{4,j}(\pi)(t)| \ne 1/2$ \ for all \ $t \in \RR_+$
   \ (since \ $|\Delta \Psi_{4,j}(\pi)(t)| \in\ZZ_+$ \ for all \ $t \in \RR_+$),
   \ hence, it remains to check that the mappings \ $\Psi_{4,j}$, \ $j \in \NN$, \ are
   \ $\cM(\RR_+ \times U) / \cD(\RR_+, \RR)$-measurable.
  This follows from
   \ $\{ \pi \in  M(\RR_+ \times U)
        : (\pi((0,t]\times U_{1,j}))_{t\in\{t_1,\ldots,t_L\}} \in B \}
     \in \cM(\RR_+ \times U)$
  \ for all \ $L \in \NN$, \ $t_1, \ldots, t_L \in \RR_+$ \ and \ $B \in \RR^L$,
   \ which is a consequence of the definition of \ $\cM(\RR_+ \times U)$.

 \item
  Finally, we verify that \ $I_{4,n,j,\ell}^{(i),a}(t_k) - I_{4,n,j,\ell}^{(i),b}(t_k)$ \ is a
 \ $\cD(\RR_+, \RR^d) \otimes \cM(\RR_+\times U) / \cB(\RR^d)$-measurable function of
 \ $(\bX^{(i)},p^{(i)})$.
 \ Based on the findings for \ $I_{1,n}^{(i)}(t_k)$ \ and \ $I_{3,n}^{(i)}(t_k)$, \
 it is enough to check that
 \begin{align}\label{help17}
  \sigma\bigl( \zeta_{0,j,\ell}^{(i)}, p^{(i)}_{0,j}(\zeta^{(i)}_{0,j\ell}),
               \tau_n^{(i)} \bigr)
  \cap \Omega^{(i)}_{0,j}
  \subset \sigma\bigl(\bX^{(i)}, p_{0,j}^{(i)}\bigr) \cap \Omega^{(i)}_{0,j}
  \subset \sigma\bigl(\bX^{(i)}, p^{(i)}\bigr) \cap \Omega^{(i)}_{0,j},
 \end{align}
 where \ $\zeta_{0,j,\ell}^{(i)}$ \ and \ $\Omega^{(i)}_{0,j}$ \ can be defined
 similarly as \ $\zeta_{1,j,\ell}^{(i)}$ \ and \ $\Omega^{(i)}_{0,j}$ \ for all
 \ $i\in\{1,2\}$ \ and \ $j,\ell\in\NN$, \ respectively (replacing in the
 definitions \ $U_{1,j}$ \ and \ $p^{(i)}_{1,j}$ \ by \ $U_{0,j}$ \ and
 \ $p^{(i)}_{0,j}$, \ respectively).
Note that
 \begin{align*}
  &\Big\{ \omega\in \Omega^{(i)}_{0,j} : \zeta^{(i)}_{0,j,\ell}(\omega) \in(0,t],\,
        p^{(i)}_{0,j}(\omega)(\zeta^{(i)}_{0,j,\ell}(\omega))\in B,\,
        \tau_n^{(i)}(\omega) \in[0,T]  \Big\} \\
  &\qquad = \bigcap_{n=1}^\infty \bigcup_{k=1}^n
            \Bigl\{ \omega \in \Omega^{(i)}_{0,j}
                   : N_{p^{(i)}_{0,j}(\omega)}((0, (k-1)t/n] \times U_{0,j})
                    \leq \ell - 1 , \\[-5mm]
          &\phantom{\qquad=\bigcap_{n=1}^\infty \bigcup_{k=1}^n \Bigl\{}
                    \Bigl\{ \omega \in \Omega^{(i)}_{0,j} : \;\;
                    N_{p^{(i)}_{0,j}(\omega)}(((k-1)t/n, kt/n] \times B)
                    \geq 1 , \\[-5mm]
         &\phantom{\qquad=\bigcap_{n=1}^\infty \bigcup_{k=1}^n \Bigl\{}
                   \Bigl\{ \omega \in \Omega^{(i)}_{0,j} : \;\;
                   N_{p^{(i)}_{0,j}(\omega)}((0, kt/n] \times U_{0,j})
                   \geq \ell \Bigr\}\\[-2mm]
  &\phantom{\qquad=\bigcap_{n=1}^\infty\,}
    \bigcap \Big\{ \omega\in \Omega^{(i)}_{0,j} :
             \int_0^T \int_{U_0}
              \|f(s, \bX_s^{(i)}(\omega), u)\|^2 \, \dd s \, m(\dd u)
             \geq n \Big\}
 \end{align*}
 for \ $t \in \RR_{++}$, \ $T\in\RR_+$, \ $j , \ell \in \NN,$
 \ $B\in\cB(U_{0,j})$, \ $i \in \{1, 2\}$.
\ Similarly, as it was explained in case of \ $I^{(i)}_{n,1}(t)$, \ one can
 approximate
 \ $\int_0^T \int_{U_0} \|f(s, \bX_s^{(i)}, u)\|^2 \, \dd s \, m(\dd u)$ \ by
 \ $\cD(\RR_+,\RR^d)/\cB(\RR_+)$-measurable functions of \ $\bX^{(i)}$, \ which yields
 \eqref{help17}.
\end{itemize}
Hence we obtain the statement.
\proofend

\begin{Rem}
In case of \ $f = 0$ \ and \ $g = 0$, \ the statement of Lemma
 \ref{Lemma_distribution} basically follows by Exercise (5.16) in Chapter IV
 in Revuz and Yor \cite{RevYor}, see also Lemma 12.4.5 in Weizs\"acker and
 Winkler \cite{WeizWin}.
\proofend
\end{Rem}

Next we formulate a corollary of Lemma \ref{Lemma_distribution}.

\begin{Lem}\label{Lemma_distribution2}
Let
 \ $\bigl( \Omega^{(1)}, \cF^{(1)}, (\cF_t^{(1)})_{t\in\RR_+}, \PP^{(1)}, \bW^{(1)},
           p^{(1)}, \bX^{(1)} \bigr)$ \ be a tuple satisfying
 \textup{(D1)}, \textup{(D2)}, \textup{(D3)} and \textup{(D4)(b)--(d)}
 and let
 \ $\bigl( \Omega^{(2)}, \cF^{(2)}, (\cF_t^{(2)})_{t\in\RR_+}, \PP^{(2)}, \bW^{(2)},
           p^{(2)}, \bX^{(2)} \bigr)$ \ be another tuple satisfying
 \textup{(D1)}, \textup{(D2)}, \textup{(D3)} such that
 \ $(\bX^{(2)}_t)_{t\in\RR_+}$ \ is an \ $\RR^d$-valued \ $(\cF_t^{(2)})_{t\in\RR_+}$-adapted
 c\`{a}dl\`{a}g process.
Suppose that \ $(\bW^{(1)}, p^{(1)}, \bX^{(1)})$ \ and
 \ $(\bW^{(2)}, p^{(2)}, \bX^{(2)})$ \ have the same distribution on
 \ $C(\RR_+, \RR^r) \times M(\RR_+ \times U) \times D(\RR_+, \RR^d)$.
\ Then \textup{(D4)(b)--(d)} hold for the tuple
 \ $\bigl( \Omega^{(2)}, \cF^{(2)}, (\cF_t^{(2)})_{t\in\RR_+}, \PP^{(2)}, \bW^{(2)},
           p^{(2)}, \bX^{(2)} \bigr)$ \ as well, and the processes
 \eqref{process1} and \eqref{process2} have the same distribution on
 \ $(D(\RR_+, \RR^d))^5$.
\end{Lem}

\noindent{\bf Proof.}
First we check that \
 $\PP^{(2)}\left(\int_0^t \|b(s, \bX_s^{(2)})\|\, \dd s < \infty\right) = 1$ \ for all
 \ $t \in \RR_+$.
\ Since \ $b$ \ is \ $\cB(\RR_+)\otimes\cB(\RR^d)\otimes\cB(U)/\cB(\RR^d)$-measurable and
 \ $\bX^{(1)}$ \ and \ $\bX^{(2)}$ \ have the same law,
 the processes \ $(b(s,\bX^{(1)}_s))_{s\in\RR_+}$ \ and \ $(b(s,\bX^{(2)}_s))_{s\in\RR_+}$ \
 have the same law as well.
Since the mapping
 \ $D(\RR_+,\RR^d) \ni f \mapsto \big(\int_0^t f(s)\,\dd s\big)_{t\in\RR_+}\in D(\RR_+,\RR^d)$
 \ is continuous
 (see, e.g.,  Ethier and Kurtz \cite[Chapter III, Section 11, Exercise 26]{EthKur},
 or Barczy et al.\ \cite[Proof of Lemma B.3]{BarIspPap}), and consequently
 \ $\cD(\RR_+,\RR^d)/\cD(\RR_+,\RR^d)$-measurable,
 the processes \ $\Big(\int_0^t \|b(s, \bX_s^{(1)})\|\, \dd s\Big)_{t\in\RR_+}$ \ and
 \ $\Big(\int_0^t \|b(s, \bX_s^{(2)})\|\, \dd s\Big)_{t\in\RR_+}$ \ have the same
 distribution with respect to \ $\PP^{(1)}$ \ and \ $\PP^{(2)}$, \ respectively.
Since \ $\PP^{(1)}\left(\int_0^t \|b(s, \bX_s^{(1)})\|\, \dd s < \infty\right) = 1$
 \ for all \ $t\in\RR_+$, \ this yields
 \ $\PP^{(2)}\left(\int_0^t \|b(s, \bX_s^{(2)})\|\, \dd s < \infty\right) = 1$ \ for
 all \ $t \in \RR_+$, \ as desired.

Similarly, one can check that
 \ $\PP^{(2)}\left(\int_0^t \|\sigma(s, \bX_s^{(2)})\|^2 \, \dd s < \infty\right) = 1$
 \ for all \ $t \in \RR_+$, \ and
 \begin{equation*}
   \PP^{(2)}\left( \int_0^t \int_{U_0} \|f(s, \bX^{(2)}_s, u)\|^2 \, \dd s \, m(\dd u)
                    < \infty \right)
          = 1 , \qquad t \in \RR_+ .
 \end{equation*}

It remains to check that
 \begin{equation}\label{D4_D_X2}
  \PP^{(2)}\left( \int_0^t \int_{U_1} \|g(s, \bX^{(2)}_{s-}, u)\| \, N^{(2)}(\dd s, \dd u)
            < \infty \right)
   = 1 , \qquad t \in \RR_+ ,
 \end{equation}
 where \ $N^{(2)}(\dd s, \dd u)$ \ is the counting measure of \ $p^{(2)}$ \ on
 \ $\RR_{++} \times U$.
\  Recall that, in the proof of Lemma \ref{Lemma_distribution}, \ $U_{1,j} \in \cB(U)$, \ $j \in \NN$,
 \ have been chosen such that they are disjoint,
 \ $m(U_{1,j}) < \infty$, \ $j \in \NN$, \ and
 \ $U_1 = \bigcup_{j=1}^\infty U_{1,j}$.
\ Further, the set \ $D(p_{1,j}^{(i)})$ \ is ordered according to magnitude as
 \ $0 < \zeta^{(i)}_{1,j,1} < \zeta^{(i)}_{1,j,2} < \cdots$, \ $j \in \NN$,
 \ $i \in \{1, 2\}$, \ see \eqref{def_zeta}.
Then for each \ $i \in \{1, 2\}$ \ and \ $t \in \RR_+$,
  \ $K_n^{(i)}(t) \to \int_0^t \int_{U_1} \|g(s, \bX^{(i)}_{s-}, u)\|  \, N^{(i)}(\dd s, \dd u)$
 \ as \ $n \to \infty$ \ $\PP^{(i)}$-almost surely, where
 \begin{align*}
  K_n^{(i)}(t)
  &:= \sum_{j=1}^n
       \int_0^t \int_{U_{1,j}} \|g(s, \bX^{(i)}_{s-}, u)\| \, N^{(i)}(\dd s, \dd u)
   = \sum_{j=1}^n
      \sum_{s\in(0,t]\cap D(p_{1,j}^{(i)})} \|g(s, \bX_{s-}^{(i)}, p_{1,j}^{(i)}(s))\| ,
 \end{align*}
 where \ $p_{1,j}^{(i)}$ \ denotes the thinning of \ $p^{(i)}$ \ onto \ $U_{1,j}$.
\ Since \ $(p^{(1)},\bX^{(1)})$ \ and \ $(p^{(2)},\bX^{(2)})$ \ have the
 same distribution with
 respect to \ $\PP^{(1)}$ \ and \ $\PP^{(2)}$, \ respectively, \ $K^{(1)}_n(t)$
 \ and \ $K^{(2)}_n(t)$ \ have the same distribution with respect to
 \ $\PP^{(1)}$ \ and \ $\PP^{(2)}$, \ respectively for all \ $n \in \NN$ \ and
 \ $t \in \RR_+$ \ (which can be checked in the same way as in the proof of
 Lemma \ref{Lemma_distribution} by replacing \ $g$ \ with \ $\Vert g\Vert$).
\ Consequently,
 \ $\int_0^t \int_{U_1} \|g(s, \bX^{(1)}_{s-}, u)\| \, N^{(1)}(\dd s, \dd u)$ \ and
 \ $\int_0^t \int_{U_1} \|g(s, \bX^{(2)}_{s-}, u)\| \, N^{(2)}(\dd s, \dd u)$
 \ have the same distribution with respect to \ $\PP^{(1)}$ \ and \ $\PP^{(2)}$,
 \ respectively for all \ $t \in \RR_+$.
\ Since
 \begin{equation*}
  \PP^{(1)}\left( \int_0^t \int_{U_1} \|g(s, \bX^{(1)}_{s-}, u)\| \, N^{(1)}(\dd s, \dd u)
            < \infty \right)
  = 1 , \qquad t \in \RR_+ ,
 \end{equation*}
 we have \eqref{D4_D_X2}.
All in all, the tuple
 \ $\bigl( \Omega^{(2)}, \cF^{(2)}, (\cF_t^{(2)})_{t\in\RR_+}, \PP^{(2)}, \bW^{(2)},
           p^{(2)}, \bX^{(2)} \bigr)$ \ satisfies \textup{(D4)(b)--(d)},
and then Lemma \ref{Lemma_distribution} yields that the processes
 \eqref{process1} and \eqref{process2} have the same distribution on
 \ $(D(\RR_+, \RR^d))^5$.
\proofend

The next lemma corresponds to Fact B on page 107 in Situ \cite{Sit}.

\begin{Lem}\label{Lemma_filtration}
Let us consider the filtered probability space
 \ $\bigl( \Omega, \cF, (\cF_t)_{t\in\RR_+}, \PP_{1,2}\bigr)$ \ given in the proof
 of Theorem \ref{Thm_pathwise_1}.
The process \ $\Omega \ni (\bx, w, \pi, y^{(1)}, y^{(2)}) \mapsto w_t \in \RR^r$,
 \ $t \in \RR_+$, \ is an $r$-dimensional standard $(\cF_t)_{t\in\RR_+}$-Brownian
 motion, and the process
 \ $\Omega \ni (\bx, w, \pi, y^{(1)}, y^{(2)})
    \mapsto N_{p_\pi}\vert_{(0,t]\times U} \in M(\RR_+ \times U)$,
 \ $t \in \RR_+$, \ is a stationary $(\cF_t)_{t\in\RR_+}$-Poisson point
 process on \ $U$ \ with characteristic measure \ $m$ \ under the measure
 \ $\PP_{1,2}$.
\end{Lem}

\noindent
\textbf{Proof.}
Using that the \ $w$-coordinate process is an \ $r$-dimensional standard
 \ $(\cG_t)_{t\in\RR_+}$-Brownian motion under \ $\PP_{1,2}$, \ for the first
 statement, it is enough to prove the independence of \ $w_t - w_s$ \ and
 \ $\cF_s$ \ for every \ $s, t \in \RR_+$ \ with \ $s < t$.
\ For this, it is sufficient to show
 \begin{align}\label{help_jobbrol_folytonossag}
  \EE_{\PP_{1,2}}( \ee^{\ii\langle\by,w_t-w_s\rangle} \bbone_G )
  = \ee^{-(t-s) \|\by\|^2/2} \, \PP_{1,2}(G), \qquad
    \by \in \RR^r , \quad G \in \cG_s , \quad 0 \leq s < t .
 \end{align}
Indeed, if \ $A \in \widetilde \cG_s$, \ then there exists some \ $G \in \cG_s$
 \ such that \ $A\Delta G = (A\setminus G) \cup(G\setminus A) \in \cN$,
 \ and consequently \ $\PP_{1,2}(A \Delta G) = 0$.
\ Then,
 \begin{align*}
  \EE_{\PP_{1,2}}( \ee^{\ii\langle\by,w_t-w_s\rangle} \bbone_A )
  & = \EE_{\PP_{1,2}}( \ee^{\ii\langle\by,w_t-w_s\rangle} \bbone_{A\cap G} )
    = \EE_{\PP_{1,2}}( \ee^{\ii\langle\by,w_t-w_s\rangle} \bbone_G ) \\
  & = \ee^{-(t-s) \|\by\|^2/2} \, \PP_{1,2}(G)
    = \ee^{-(t-s) \|\by\|^2/2} \, \PP_{1,2}(A) ,
  \qquad A \in \widetilde\cG_s , \quad 0 \leq s < t .
 \end{align*}
Moreover, if \ $A \in \cF_s$, \ then \ $A \in \widetilde\cG_{s+\varepsilon}$ \ for
 all \ $\vare > 0$, \ and hence
 \[
   \EE_{\PP_{1,2}}( \ee^{\ii\langle\by,w_t-w_{s+\varepsilon}\rangle} \bbone_{A} )
   = \ee^{-(t-s-\varepsilon) \|\by\|^2/2} \, \PP_{1,2}(A) ,
   \qquad A \in \cF_s, \quad 0 \leq s < t , \quad \vare > 0 .
 \]
By dominated convergence theorem, using that \ $w$ \ has continuous sample
 paths \ $\PP_{1,2}$-almost surely, we get
 \[
   \EE_{\PP_{1,2}}( \ee^{\ii\langle\by,w_t-w_{s}\rangle} \bbone_{A} )
   = \ee^{-(t-s) \|\by\|^2/2} \, \PP_{1,2}(A) ,
   \qquad A \in \cF_s , \quad 0 \leq s < t ,
 \]
 i.e.,
 \[
   \EE_{\PP_{1,2}}\left[ \ee^{\ii\langle\by,w_t-w_s\rangle}  \mid \cF_s \right]
   = \ee^{-(t-s) \|\by\|^2/2}, \qquad 0 \leq s < t.
 \]
Thus, in the light of Lemma 2.6.13 of Karatzas and Shreve \cite{KarShr}, we
 get the independence of \ $w_t - w_s$ \ and \ $\cF_s$ \ for every
 \ $s, t \in \RR_+$ \ with \ $s < t$.

Using that \ $w_t - w_s$ \ is independent of \ $\cG_s$ \ under \ $\PP_{1,2}$,
 \ we obtain
 \begin{align*}
  &\EE_{\PP_{1,2}}\left[ \ee^{\ii\langle\by,w_t-w_s\rangle}
                \bbone_G \right]
   =  \EE_{\PP_{1,2}}\left[ \EE_{\PP_{1,2}}\left[ \ee^{\ii\langle\by,w_t-w_s\rangle}
                \bbone_G
                \;\big\vert\;  \cG_s \right] \right]
    = \EE_{\PP_{1,2}} \left[\bbone_G
      \EE_{\PP_{1,2}} \left[ \ee^{\ii\langle\by,w_t-w_s\rangle}
                \;\big\vert\;  \cG_s \right] \right] \\
  &\qquad
    = \EE_{\PP_{1,2}} \left[\bbone_G \EE_{\PP_{1,2}}
                      \left[ \ee^{\ii\langle\by,w_t-w_s\rangle} \right] \right]
    = \EE_{\PP_{1,2}} \left[\bbone_G
                        \ee^{-(t-s)\|\by\|^2/2}  \right]
   = \ee^{-(t-s)\|\by\|^2/2}
      \, \PP_{1,2}(G)
 \end{align*}
 for all \ $\by \in \RR^r$ \ and \ $G \in \cG_s$, \ hence we conclude
 \eqref{help_jobbrol_folytonossag} and then the first statement.

Using that the process \ $p_\pi$ \ is a stationary $(\cG_t)_{t\in\RR_+}$-Poisson
 point process on \ $U$ \ with characteristic measure \ $m$, \ as it was
 explained in the proof of the first statement, for the second statement, it
 is enough to show that for every \ $s, t \in \RR_+$ \ with \ $s < t$, \ every
 \ $n \in \NN$, \  every disjoint subsets \ $B_1, \ldots, B_n \in \cB(U)$ \ and
 \ $\lambda_1, \ldots, \lambda_n \in \RR_+$,
 \[
   \EE_{\PP_{1,2}}\left[\ee^{- \sum_{j=1}^n \lambda_j N_{p_\pi}((s,t] \times B_j)}
                \bbone_G\right]
   = \ee^{(t-s) \sum_{j=1}^n (\ee^{-\lambda_j}-1) m(B_j)}
     \, \PP_{1,2}(G) , \qquad G \in \cG_s .
 \]
Using that \ $N_{p_\pi}((s,t] \times B_j)$, \ $j \in \{1, \ldots, n\}$, \ are
 independent of each other and from \ $\cG_s$ \  under \ $\PP_{1,2}$, \ we get
 \begin{align*}
  &\EE_{\PP_{1,2}}\left[\ee^{- \sum_{j=1}^n \lambda_j N_{p_\pi}((s,t] \times B_j)}
               \bbone_G\right]
    =  \EE_{\PP_{1,2}}\left[ \EE_{\PP_{1,2}}\left[ \ee^{- \sum_{j=1}^n \lambda_j N_{p_\pi}((s,t] \times B_j)}
                \bbone_G
                \;\Big\vert\;  \cG_s \right] \right] \\
  &\qquad
    = \EE_{\PP_{1,2}} \left[\bbone_G
      \EE_{\PP_{1,2}} \left[ \ee^{- \sum_{j=1}^n \lambda_j N_{p_\pi}((s,t] \times B_j)}
                \;\Big\vert\;  \cG_s \right] \right]
    = \EE_{\PP_{1,2}} \left[\bbone_G \EE_{\PP_{1,2}}
                      \left[ \ee^{- \sum_{j=1}^n \lambda_j N_{p_\pi}((s,t] \times B_j)} \right] \right]\\
  &\qquad
    = \EE_{\PP_{1,2}} \left[\bbone_G
                        \; \ee^{(t-s) \sum_{j=1}^n (\ee^{-\lambda_j}-1) m(B_j)} \right]
   = \ee^{(t-s) \sum_{j=1}^n (\ee^{-\lambda_j}-1) m(B_j)}
      \PP_{1,2}(G)
 \end{align*}
 for all \ $G \in \cG_s$.
\ The last but one equality above is a consequence that
 \ $N_{p_\pi}((s,t] \times B_j)$ \ is a Poisson distributed random variable with
 parameter \ $(t-s) m(B_j)$, \ $j \in \{1, \ldots, n\}$, \ under \ $\PP_{1,2}$.
\ Hence we conclude the second statement as well.
\proofend


\begin{thebibliography}{99}

\bibitem{App}
\textsc{Applebaum, D.} (2009).
\textit{L\'evy Processes and Stochastic Calculus, 2nd ed}.
Cambridge University Press, Cambridge.

\bibitem{BarIspPap}
\textsc{Barczy, M., Isp\'any, M.} and \textsc{Pap, G.} (2012).
Asymptotic behavior of CLS estimators for unstable INAR(2) models.\\
Available on the ArXiv: \texttt{http://arxiv.org/abs/1202.1617}.

\bibitem{Che2}
\textsc{Cherny, A. S.} (2000).
On strong and weak uniqueness for stochastic differential equations.
\textit{Theory of Probability and its Applications}
\textbf{46(3)} 406--419.

\bibitem{Coh}
\textsc{Cohn, D. L.} (2013).
\textit{Measure Theory, 2nd ed}.
Birkh\"auser.

\bibitem{DawLi}
\textsc{Dawson, D. A.} and \textsc{Li, Z.} (2012).
Stochastic equations, flows and measure-valued processes.
\textit{The Annals of Probability}
\textbf{40(2)} 813--857.

\bibitem{DorBar}
\textsc{D\"oring, L.} and \textsc{Barczy, M.} (2012).
A jump type SDE approach to positive self-similar Markov processes.
\textit{Electronic Journal in Probability}
\textbf{17} Paper no. 94, 1--39.

\bibitem{Dud}
\textsc{Dudley, R. M.} (1989).
\textit{Real Analysis and Probability}.
Wadsworth \& Brooks/Cole Advanced Books \& Software, Pacific Grove, California.

\bibitem{Eng}
\textsc{Engelbert, H. J.} (1991).
On the theorem of T. Yamada and S. Watanabe.
\textit{Stochastics and Stochastics Reports}
\textbf{36} 205--216.

\bibitem{EthKur}
\textsc{Ethier, S. N.} and \textsc{Kurtz, T. G.} (1986).
\textit{Markov Processes},
 John Wiley \& Sons, Inc., New York.

\bibitem{IkeWat}
\textsc{Ikeda, N.} and \textsc{Watanabe, S.} (1989).
\textit{Stochastic Differential Equations and Diffusion Processes}, 2nd ed.
North-Holland/Kodansha, Amsterdam/Tokyo.

\bibitem{Jac}
\textsc{Jacod, J.} (1980).
Weak and strong solutions of stochastic differential equations.
\textit{Stochastics}
\textbf{3}, 171--191.

\bibitem{JacShi}
\textsc{Jacod, J.} and \textsc{Shiryaev, A. N.} (2003).
\textit{Limit Theorems for Stochastic Processes}, 2nd ed.
Springer-Verlag, Berlin.

\bibitem{KarShr}
\textsc{Karatzas, I.} and \textsc{Shreve, S. E.} (1991).
\textit{Brownian Motion and Stochastic Calculus}, 2nd ed.
Springer, Berlin.

\bibitem{Kec}
\textsc{Kechris, A. S.} (1995).
\textit{Classical descriptive set theory}.
Graduate Texts in Mathematics, 156.
Springer-Verlag, New York.

\bibitem{Kur1}
\textsc{Kurtz, T. G.} (2007).
The Yamada--Watanabe--Engelbert theorem for general stochastic equations and
 inequalities.
\textit{Electronic Journal of Probability}
\textbf{12} Paper no. 33, 951--965.

\bibitem{Kur2}
\textsc{Kurtz, T. G.} (2014).
Weak and strong solutions of general stochastic models.
\textit{Electronic Communications in Probability}
\textbf{19} Paper no. 58, 1--16.

\bibitem{Li}
\textsc{Li, Z. H.} (2011).
\textit{Measure-Valued Branching Markov Processes}.
 Springer-Verlag, Heidelberg.

\bibitem{LiMyt}
\textsc{Li, Z.} and \textsc{Mytnik, L.} (2011).
Strong solutions for stochastic differential equations with jumps.
\textit{Annales de l'Institut Henri Poincar\'e (B) Probability and
         Statistics}
\textbf{47(4)}, 1055--1067.

\bibitem{LiPu}
\textsc{Li, Z.} and \textsc{Pu, F.} (2012).
Strong solutions of jump-type stochastic equations.
\textit{Electronic Communications in Probability}
\textbf{17} Paper no. 33, 1--13.

\bibitem{Ond}
\textsc{Ondrej\'at, M.} (2004).
Uniqueness for stochastic evolution equations in Banach spaces.
\textit{Dissertationes Math. (Rozprawy Mat.)}
\textbf{426}, 63 pp.

\bibitem{Res}
\textsc{Resnick, S. I.} (2008).
\textit{Extreme Values, Regular Variation, and Point Processes.}
Springer Science+Business Media, LLC.

\bibitem{RevYor}
\textsc{Revuz, D.} and \textsc{Yor, M.} (2001).
\textit{Continuous Martingales and Brownian Motion},
 3rd ed., corrected 2nd printing.
Springer-Verlag, Berlin.

\bibitem{RocSchZha}
\textsc{R\"ockner, M., Schmuland, B.} and \textsc{Zhang, X.} (2008).
Yamada--Watanabe theorem for stochastic evolution equations in infinite
 dimensions.
\textit{Condensed Matter Physics}
\textbf{11} No 2(54) 247--259.

\bibitem{Sit}
\textsc{Situ, R.} (2005).
\textit{Theory of Stochastic Differential Equations with Jumps and
        Applications}.
Springer, New York.

\bibitem{Tap}
\textsc{Tappe, S.} (2013).
The Yamada--Watanabe theorem for mild solutions to stochastic partial
 differential equations.
\textit{Electronic Communications in Probability}
\textbf{18} no. 24, 1--13.

\bibitem{Vaart}
\textsc{van der Vaart, A. W.} (1998).
\textit{Asymptotic Statistics},
Cambridge University Press.

\bibitem{WeizWin}
\textsc{von Weizs\"acker, H.} and \textsc{Winkler, G.} (1990).
\textit{Stochastic Integrals. An Introduction},
Advanced Lectures in Mathematics.
Friedr.\ Vieweg \& Sohn, Braunschweig.

\bibitem{YamWat}
\textsc{Yamada, T.} and \textsc{Watanabe, S.} (1971).
On the uniqueness of solutions of stochastic differential equations.
\textit{Journal of Mathematics of Kyoto University}
\textbf{11(1)} 155--167.

\bibitem{Zhao}
\textsc{Zhao, H.} (2014).
Yamada-Watanabe theorem for stochastic evolution equation driven by Poisson random measure.
\textit{International Scholarly Research Notices. Probability and Statistics}
\textbf{2014} Article ID 982190.

\end{thebibliography}
\end{document}